\documentclass{article}

\newenvironment{rcases}
  {\left.\begin{aligned}}
  {\end{aligned}\right\rbrace}

\usepackage{amssymb,amsmath,amsthm,mathrsfs}
\usepackage{multirow}
\usepackage{tabularx}
\usepackage{diagbox}
\usepackage{cite}
\usepackage{indentfirst}
\usepackage{graphicx,subcaption}
\usepackage[table]{xcolor}
\usepackage{standalone}
\usepackage{relsize}
\usepackage{structuralanalysis}
\usepackage{changes}
\graphicspath{ {images/} }
\usepackage{tikz,array,calc}
\usetikzlibrary{decorations.pathreplacing}
\definecolor{green1}{RGB}{217,240,211}

\usepackage{xfrac}

\DeclareMathOperator*{\A}{ \mathlarger{\mathlarger{\mathlarger{\boldsymbol{\mathsf{A}}}}} }

\newcommand\mi{\mathbf{I}}

\newcommand{\rn}{\mathbb{R}^d}

\renewcommand{\vec}[1]{\boldsymbol{\mathbf{#1}}}
\newcommand{\stkout}[1]{\ifmmode\text{\sout{\ensuremath{#1}}}\else\sout{#1}\fi}
\setdeletedmarkup{\stkout{#1}}
\usepackage{geometry}
\newcommand{\Bezier}{{B\'{e}zier} }
\renewenvironment{abstract}{%
\hfill\begin{minipage}{0.95\textwidth}
\rule{\textwidth}{1pt}}
{\par\noindent\rule{\textwidth}{1pt}\end{minipage}}
\providecommand{\keywords}[1]{\textbf{\textit{Keywords: }}#1}

\title{B\'{e}zier $\bar{B}$ projection}
\author{Di Miao$^{1,}$\thanks{miaodi1987@gmail.com} \and Michael J. Borden$^1$ \and Michael A. Scott$^1$ \and Derek C. Thomas$^2$}
\date{}

\begin{document}
\maketitle
\vspace{-2em}
%
\begin{center}

$^1$Department of Civil and Environmental Engineering\\
Brigham Young University\\
368 CB, Provo, UT 84602, USA\\ \vspace{1em}
 
$^2$Coreform LLC\\
P.O. Box 970336, Orem, UT 84097, USA \\ 
\end{center}


\begin{abstract}
  In this paper we demonstrate the use of B\'{e}zier projection to alleviate locking phenomena in structural mechanics applications of isogeometric analysis. {Interpreting the well-known $\bar{B}$ projection in two different ways we develop two formulations for locking problems in beams and nearly incompressible elastic solids.  One formulation leads to a sparse symmetric symmetric system and the other leads to a sparse non-symmetric system.} To demonstrate the utility of {\Bezier projection} for both geometry and material locking phenomena we focus on transverse shear locking in Timoshenko beams and volumetric locking in nearly compressible linear elasticity although the approach can be applied generally to other types of locking phenemona as well. B\'{e}zier projection is a local projection technique with optimal approximation properties, which in many cases produces solutions that are comparable to global $L^2$ projection. In the context of $\bar{B}$ methods, the use of B\'ezier projection produces sparse stiffness matrices with only a slight increase in bandwidth when compared to standard displacement-based methods. Of particular importance is that the approach is applicable to any spline representation that can be written in B\'ezier form like NURBS, T-splines, LR-splines, etc. We discuss in detail how to integrate this approach into an existing finite element framework with minimal disruption through the use of B\'ezier extraction operators and a newly introduced dual {basis for} the \Bezier {projection} operator. We then demonstrate the behavior of the two proposed formulations through several challenging benchmark problems.
\end{abstract}
\keywords{Isogeometric analysis, B\'ezier extraction, {\Bezier dual basis}, B\'ezier projection, $\bar{B}$-projection, locking}
\section{Introduction}
Isogeometric analysis (IGA), introduced by Hughes et al. \protect\cite{hughes_isogeometric_2005}, adopts the spline basis, which underlies the CAD geometry, as the basis for analysis. Of particular importance is the positive impact of smoothness on numerical solutions, where, in many application domains, IGA outperforms classical finite elements~\protect\cite{cottrell_isogeometric_2009,cottrell_studies_2007,cottrell2006isogeometric,hughes_duality_2008,bazilevs_isogeometric_2010,evans_n-widths_2009}. Initial investigations of IGA focused on non-uniform rational B-splines (NURBS) due to their dominance in commercial CAD packages. However, many advances are being made in analysis-suitable geometry representations that overcome the strict rectangular topological restrictions of NURBS. Examples include T-splines~\protect\cite{bazilevs_isogeometric_2010,sederberg_t-splines_2003} and their analysis-suitable restriction~\protect\cite{scott_local_2012, li_analysis-suitable_2013}, hierarchical B-splines~\protect\cite{bornemann_subdivision-based_2013,scott_isogeometric_2014,schillinger_isogeometric_2012,evans_hierarchical_2015,forsey_hierarchical_1988}, and locally refined B-splines~\protect\cite{dokken_polynomial_2013,johannessen_isogeometric_2014} among others.

The purpose of this paper is to demonstrate how B\'{e}zier projection~\protect\cite{thomas_bezier_2015} can be employed as the underlying local projection framework for a $\bar{B}$ approach to treat locking in isogeometric structural elements. B\'ezier projection is an element-based local projection methodology for B-splines, NURBS, and T-splines. It relies on the concept of B\'ezier extraction~\protect\cite{borden_isogeometric_2011, scott_isogeometric_2011} and an associated operation, spline reconstruction, which enables the use of B\'ezier projection in standard finite element codes.

B\'ezier projection exhibits provably optimal convergence and yields
projections that are virtually indistinguishable from global $L^2$ projection. For an isogeometric finite element code that leverages B\'ezier extraction, B\'ezier projection can be employed virtually for free. To simplify the implementation of the \Bezier $\bar{B}$ method in existing finite element codes we develop a \textit{dual} element \Bezier extraction operator that can be derived directly from the \Bezier extraction of a spline representation. It is worth noting that B\'ezier projection can also be used to develop a unified framework for spline operations including cell subdivision and merging, degree elevation and reduction, basis roughening and smoothing, and spline reparameterization and is applicable to any spline representation that can be written in B\'ezier form.

Numerical locking in structural finite elements includes geometric locking in thin curved structural members such as membrane and shear locking and also includes volumetric locking in incompressible and nearly incompressible elasticity. There is an immense body of literature on approaches to overcome locking in the finite element community and various approaches have emerged as dominant. These include reduced quadrature~\protect\cite{malkus_mixed_1978,zienkiewicz_reduced_1971}, $\bar{B}$ projection methods~\protect\cite{nagtegaal_numerically_1974,hughes_generalization_1980, NME:NME4328}, and mixed methods based on the Hu-Washizu variational principle~\protect\cite{dolbow_volumetric_1999,kasper_mixed-enhanced_2000,hughes_variational_1986, NME:NME3048}. It is important to mention that, although ameliorated at high polynomial degrees, smooth splines in the context of IGA still exhibit locking behavior~\protect\cite{echter_numerical_2010, bouclier_locking_2012}.

In IGA, there is a growing literature on the treatment of locking in structural elements. Leveraging higher-order smoothness, transverse shear locking can be eliminated at the theoretical level by employing Kirchhoff-Love~\protect\cite{kiendl_isogeometric_2009, kiendl_isogeometric_2015} and hierarchic Reissner-Mindlin~\protect\cite{oesterle_hierarchic_2017, oesterle_shear_2016,echter_hierarchic_2013} shell elements. Reduced quadrature schemes have been explored in~\protect\cite{adam_improved_2014, adam_improved_2015, adam_selective_2015} as a way to alleviate transverse shear locking. The extension of $\bar{B}$ projection to the isogeometric setting was initiated in~\protect\cite{elguedj:hal-00457010} for both elastic and plastic problems and was extended in~\protect\cite{bouclier_efficient_2013} to include local projection techniques~\protect\cite{mitchell_method_2011,govindjee_convergence_2012}.

{In this paper we introduce two methods that employ \Bezier projection to produce a localized approximation to the standard $\bar{B}$ method. The motivation behind these methods is the fact that $\bar{B}$ methods result in dense linear systems. The methods we introduce result in a sparse linear system irrespective of the choice of basis functions. We call these two methods the symmetric and non-symmetric \Bezier $\bar{B}$ projection methods, where the names indicate the symmetry of the resulting stiffness matrix. We show that both methods result in a sparse stiffness matrix and reduce locking. We also show that optimal convergence rates are achieved in the case of the non-symmetric method and near optimal convergence rates are achieved in the case of the symmetric method.  We also perform an inf-sup analysis of these methods.}

The outline of this paper is as follows: First, we briefly review spline basis functions in Section~\protect\ref{sec:preliminaries}. In Sections~\protect\ref{sec:extraction} and~\protect\ref{sec:bproject}, we describe B\'ezier extraction and projection. We then formulate and use B\'ezier $\bar{B}$ projection for the Timoshenko beam (to treat transverse shear locking) and nearly incompressible elasticity (to treat volumetric locking) in Sections~\protect\ref{sec:timoshenkobeam} and~\protect\ref{sec:nearlyincompressible}, respectively. We provide detailed element level operations in both settings. We also presents numerical tests to show the performance of the proposed strategy. 

\section{Preliminaries and notation}
\label{sec:preliminaries}
In this section a brief overview of univariate Bernstein, B-spline, and NURBS basis functions is provided. We also describe how these univariate basis functions are extended to higher dimensions.
\subsection{Univariate Bernstein basis functions}
The $i$th univariate Bernstein basis function of degree $p$ is defined by
\begin{equation}
B_{i,p}(\xi)=\binom {p}{i}\xi^i(1-\xi)^{n-i}
\end{equation}
where $\xi\in\left[ 0,1 \right]$ and $\binom {p}{i}=\dfrac{p!}{i!(n-i)!}$, $0\leq{i}\leq{p}$, is a binomial coefficient.
\subsection{Univariate spline basis functions}
A univariate B-spline basis of dimension $n$ is defined by a polynomial degree $p$ and a knot vector $\mathbf{\Xi}=\lbrace{\xi_0,\xi_1,\ldots, \xi_{n+p}}\rbrace$, which is a non-decreasing sequence of real numbers. The $A$th B-spline basis function can then be defined using the Cox-de Boor recursion formula:
\begin{gather}
N_{A,0}(\xi)=\begin{cases}1 & \xi_A\leq{\xi}\leq{\xi_{A+1}}\\0 & otherwise \end{cases} \\
N_{A,p}(\xi)=\dfrac{\xi-\xi_A}{\xi_{A+p}-\xi_A}N_{A,p-1}(\xi)+\dfrac{\xi_{A+p+1}-\xi}{\xi_{A+p+1}-\xi_{A+1}}N_{A+1,p-1}(\xi).
\end{gather}\par
For simplicity, we will always use open knot vectors defined over the interval $\left[ 0,1\right]$. An open knot vector satisfies the conditions $\xi_0=\xi_1=\dots=\xi_{p}=0$ and $\xi_{n}=\xi_{n+1}=\dots=\xi_{n+p}=1$ and creates interpolatory end conditions. B-spline basis functions can be used to represent piecewise polynomial functions but are not capable of representing conic sections (e.g. circles, ellipses and hyperbolas). NURBS overcome this shortcoming. A NURBS basis function can be written as
\begin{equation}
R_{A,p}(\xi)=\dfrac{N_{A,p}(\xi)w_A}{W{(\xi)}}
\end{equation}
where $w_A$ is called a weight and
\begin{align}
  \label{eq:weight}
W(\xi)=\sum_{A} N_{A,p}(\xi)w_A
\end{align}
is called the weight function. A $d$-dimensional rational curve $\mathbf{S}(\xi)\in{\rn}$ can then be defined as
\begin{equation}
\mathbf{S}(\xi)=\sum_A R_{A,p}(\xi)\mathbf{P}_A
\end{equation}
where $\mathbf{P}_A=(p_A^1,p_A^2,\ldots,p_A^d)^T$ is a $d$-dimensional control point. It is often more convenient to represent the $d$-dimensional NURBS in a $(d+1)$-dimensional homogeneous space by defining $\mathbf{P}_A^w=(p_A^1w_A,p_A^2w_A,\ldots,p_A^dw_A,w_A)^T$ and the corresponding $(d+1)$-dimensional B-spline curve as
\begin{align}
\mathbf{S}^w(\xi)&=\sum_A N_{A,p}(\xi)\mathbf{P}_A^w
\end{align}
such that each component of $\mathbf{S}^w$ can be written as
\begin{align}
S_i(\xi)&=\dfrac{{S}_i^w(\xi)}{{S}_{d+1}^w(\xi)}.
\end{align}
In the homogeneous form, NURBS can be manipulated with standard B-spline algorithms.

\subsection{Multivariate spline basis functions}
In higher dimensions, Bernstein, B-spline, and NURBS basis functions are formed by the Kronecker product of univariate basis functions. For example, two-dimensional B-spline basis functions of degree $\mathbf{p}=(p_\xi, p_\eta)$ are defined by
\begin{equation}
\mathbf{N}^\mathbf{p}(\xi,\eta)=\mathbf{N}^{p_\xi}(\xi)\otimes\mathbf{N}^{p_\eta}(\eta)
\end{equation}
where $\mathbf{N}^{p_\xi}(\xi)$ and $\mathbf{N}^{p_\eta}(\eta)$ are vectors of basis functions in the $\xi$ and $\eta$ directions, respectively. A particular multivariate basis function can be written as
\begin{equation}
{N}_{A(i,j)}^\mathbf{p}(\xi,\eta)={N}_{i,p_\xi}(\xi){N}_{j,p_\eta}(\eta)
\end{equation}
where the index mapping is defined as
\begin{equation}
A(i,j)=n_\eta{i}+j.
\end{equation}
The integer $n_\eta$ is the number of basis functions in $\eta$ direction.
\section{B\'{e}zier extraction}
\label{sec:extraction}
\begin{figure}
  \centering
  \includegraphics[width=.6\linewidth]{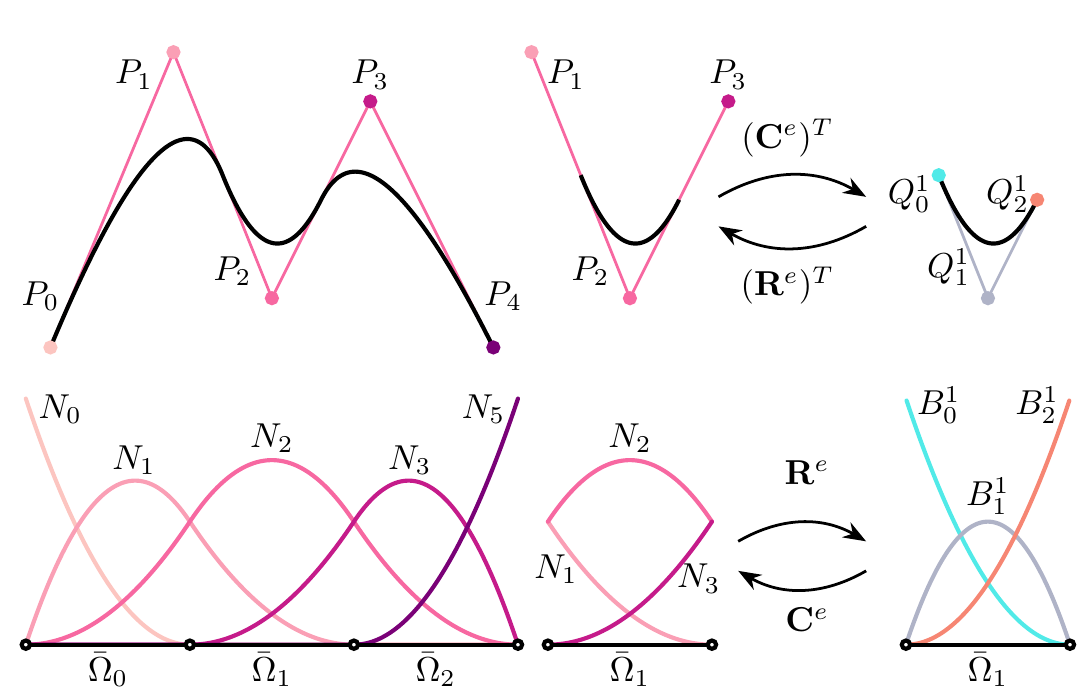}
  \caption{Illustration of B\'ezier extraction and projection in one dimension for a B-spline of degree 2 and knot vector $[0,0,0,1/3/2/3,1,1,1]$ (restricted to the second element for illustrative purposes).}
    \label{fig:extraction_and_projection}
\end{figure}

Given a spline basis $\mathbf{N}$ there exists a Bernstein basis $\mathbf{B}$ and a linear operator $\mathbf{C}$ (see~\protect\cite{borden_isogeometric_2011}) such that
\begin{equation}
\mathbf{N}(\xi)=\mathbf{C}\mathbf{B}(\xi).
\end{equation}
The localization of $\mathbf{C}$ to an element domain produces the element extraction operator $\mathbf{C}^e$.
Given control points $\mathbf{P}^e$, the corresponding B\'ezier control points $\mathbf{Q}^e$ can be computed directly as
\begin{equation}
\mathbf{Q}^e=(\mathbf{C}^e)^T\mathbf{P}^e.
\end{equation}
A graphical depiction of B\'{e}zier extraction is shown in Figure~\protect\ref{fig:extraction_and_projection}.

\section{B\'ezier projection}
\label{sec:bproject}

B\'{e}zier projection can be viewed as the inverse of extraction~\protect\cite{thomas_bezier_2015}. B\'ezier projection uses an element reconstruction operator $\mathbf{R}^e\equiv(\mathbf{C}^e)^{-1}$ such that the global control point values, corresponding to those basis functions defined over the support of an element $e$, can be determined directly from \Bezier control values as
\begin{equation}
\mathbf{P}^e=(\mathbf{R}^e)^T\mathbf{Q}^e
\end{equation}
where $\mathbf{Q}^e$ is any field in B\'ezier form. The action of the element reconstruction operator is depicted graphically in Figure~\protect\ref{fig:extraction_and_projection}. For example, given any function $u \in L^2$, we can compute $\mathbf{Q}^e$ as
\begin{equation}
\mathbf{Q}^e=(\mathbf{G}^e)^{-1}\mathbf{F}^e
\label{eq:element-Qi}
\end{equation}
where $\mathbf{G}^e$ is the Gramian matrix corresponding to the Bernstein basis with components
\begin{align}
  {G}_{ij}^e &= \int_{\Omega^e} B^e_i B^e_j \, d\Omega =\langle{B^e_{i},B^e_{j}}\rangle_{\Omega^e}
\end{align}
and
\begin{align}
  {F}^e_i &=  \int_{\Omega^e} B^e_i u \, d\Omega = \langle{B^e_{i,},u}\rangle_{\Omega^e}.
\end{align}
Note that efficiency gains can be had at the expense of accuracy by instead performing the integration in the parametric domain of the element~\protect\cite{thomas_bezier_2015}. 

\begin{figure}[htb]
  \centering
  \includegraphics[width=5in]{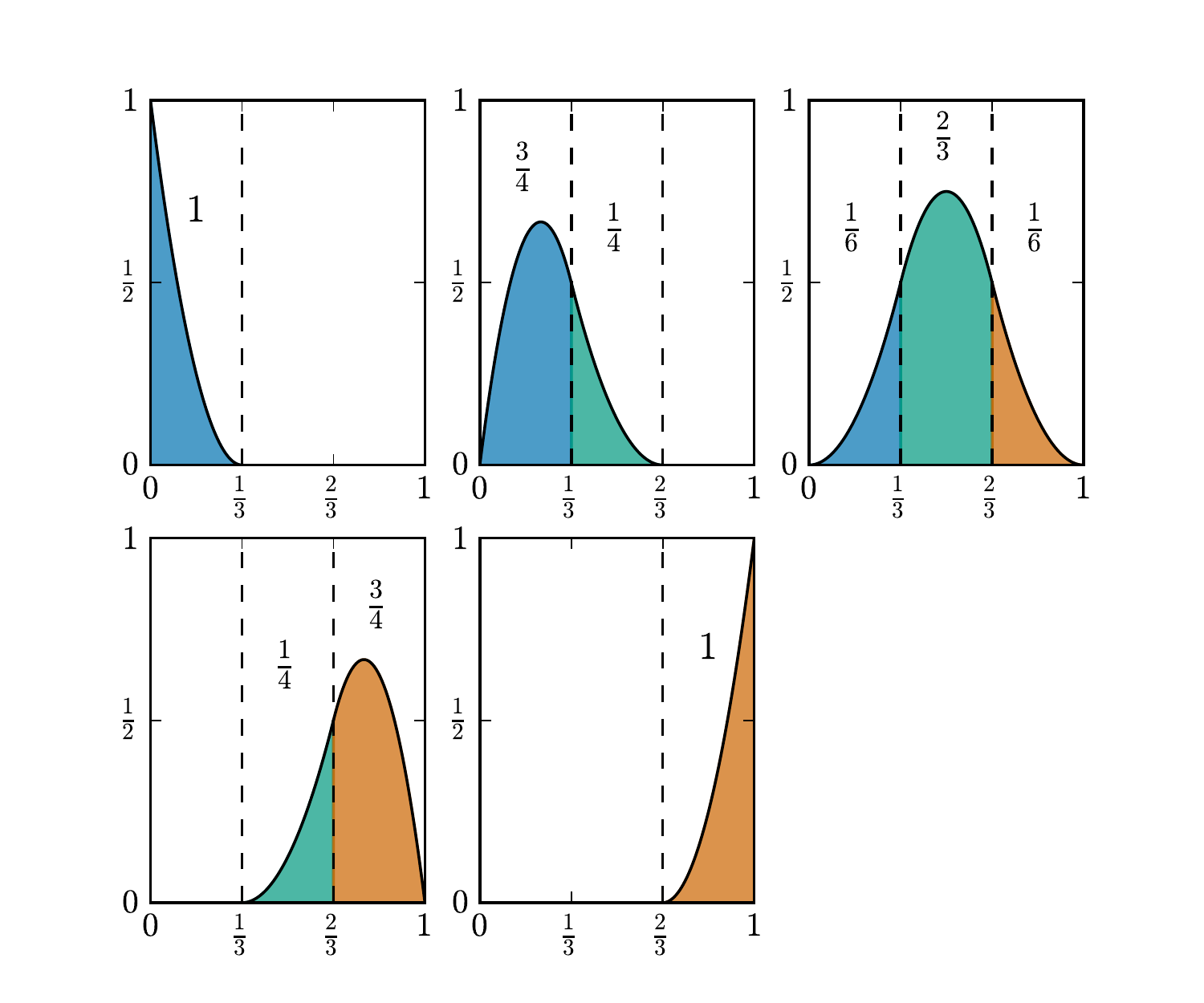}
  \caption{\label{fig:weights}Weights over each knot span associated
    with the basis function defined by the knot vector
    $[0,0,0,\sfrac{1}{3},\sfrac{2}{3},1,1,1]$.} 
\end{figure}

The element-wise projection produces one control value for each element in the support of the function.  These values must be combined in order to provide the final control value.  A core component of the B\'ezier projection algorithm is the definition of an appropriate averaging operation. {The process of computing the weights is illustrated in Figure~\protect\ref{fig:weights}.} A weighted average of the values is computed using the weighting
\begin{equation}\label{eqn:Bezier_weight}
\omega_a^e=\dfrac{\int_{\Omega^e} N_{a}^e \, d\Omega}{\int_{\Omega^A} N_{A(e,a)} \, d\Omega}
\end{equation}
where $\Omega^e$ corresponds to the physical domain of element $e$, $A(e,a)$ is a mapping from a local nodal index $a$ defined over element $e$ to a corresponding global node index $A$, and $\Omega^A$ corresponds to the physical support of $N_A$. The final averaged global control point is then calculated as
\begin{equation}
P_A=\sum_{\Omega^e\in \Omega^A } \omega_{A(e,a)} P_{A(e,a)}.
\end{equation}
B\'ezier projection onto NURBS functions can be defined in an analogous manner~\protect\cite{thomas_bezier_2015}.

\begin{figure}[htb]
  \centering
  \begin{tabular}{c p{4in} p{2in}}
    (0)&\raisebox{-.5\height}{\includegraphics[width=4in]{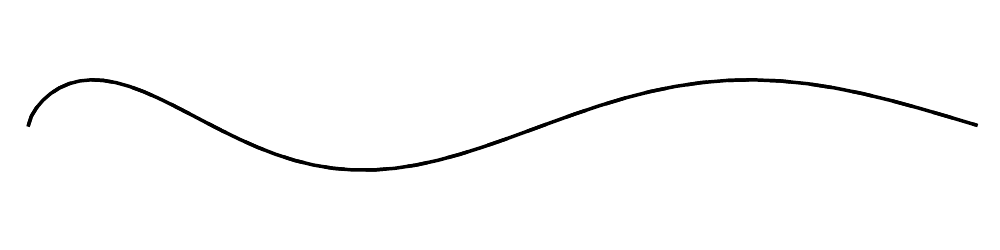}} & \begin{minipage}[t]{2in}Target function\end{minipage}\\
    (1)&\raisebox{-.5\height}{\includegraphics[width=4in]{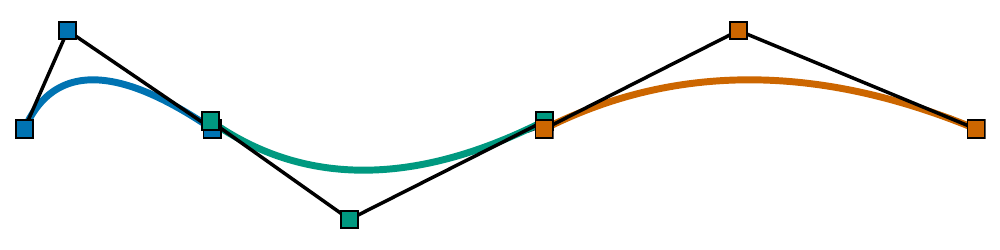}} & \begin{minipage}[t]{2in}Perform local projection to obtain \Bezier control points (represented by squares, colored to match elements)\end{minipage}\\
    (2)&\raisebox{-.5\height}{\includegraphics[width=4in]{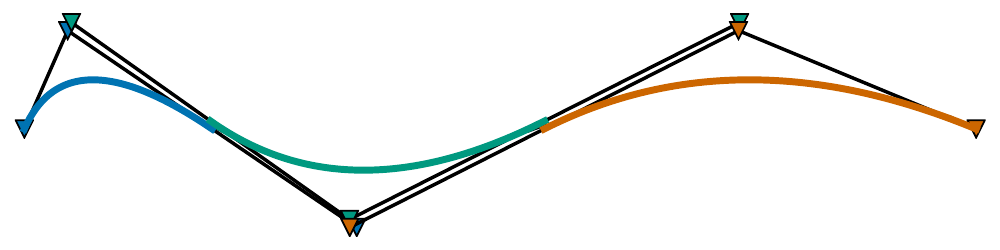}} & \begin{minipage}[t]{2in}Use element reconstruction operator to project \Bezier points to spline control points (represented by inverted triangles, colored to match elements)\end{minipage}\\
    (3)&\raisebox{-.5\height}{\includegraphics[width=4in]{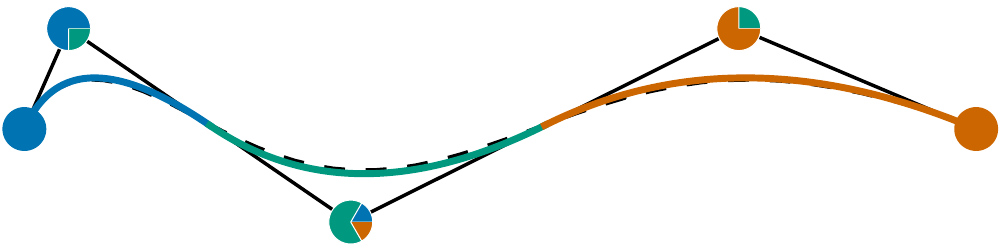}} & \begin{minipage}[t]{2in}Apply smoothing algorithm (contribution of each element to each control point shown by colored fraction)\end{minipage}\\
    (4)&\raisebox{-.5\height}{\includegraphics[width=4in]{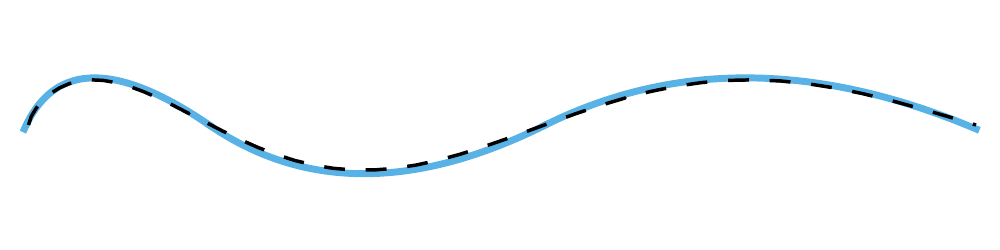}} & \begin{minipage}[t]{2in}Comparison of final function (light ) and target function (dashed)\end{minipage}
  \end{tabular}
  \caption{\label{fig:loc-proj-example}Steps of \Bezier projection.}
\end{figure}

The individual steps comprising the \Bezier projection algorithm are
illustrated in Figure~\protect\ref{fig:loc-proj-example} where
the curve defined by $\vec{f}(t)=\left( \frac{t}{3}
\right)^{3/2}\vec{e}_1+\frac{1}{10}\sin (\pi t )\,\vec{e}_2$,
$t\in[0,3]$ is projected onto the quadratic B-spline basis defined by
the knot vector $[0,0,0,\sfrac{1}{3},\sfrac{2}{3},1,1,1]$. For this example,
the algorithm proceeds as follows:

\begin{description}

\item{Step 1:} The function $\mathbf{f}$ is projected onto the Bernstein basis of each element. This results in a set of
  \Bezier coefficients that define an approximation to $\mathbf{f}$.
  The \Bezier coefficients are indicated in part (1) of Figure~\protect\ref{fig:loc-proj-example} by 
  square markers that have been colored to match the corresponding
  element. Each \Bezier segment is discontinuous.

\item{Step 2:} The element reconstruction operator $\mathbf{R}^e$ is used to convert the
  \Bezier control points into spline control points associated with the
  basis function segments over each element.
  The new control points are marked with inverted triangles and
  again colored to indicate the element with which the control point is
  associated. The control points occur in clusters.
  The clusters of control points represent the contributions from
  multiple elements to a single spline basis function control point.

\item{Step 3:} Each cluster of control points is averaged to obtain a
  single control point by weighting each point in the cluster according
  to the weighting given in (\protect\ref{eqn:Bezier_weight}). The resulting control
  points are shown as circles with the relative contribution from each
  element to each control point indicated by the colored fraction of the
  control point marker. Colors in Figures~\protect\ref{fig:weights} and~\protect\ref{fig:loc-proj-example} are coordinated
  to illustrate where the averaging weights come from and their values.
\end{description}


\subsection{Dual basis formulation of B\'{e}zier projection}
\label{sec:dual-funct-dual}
To integrate \Bezier projection into a standard finite element assembly algorithm, it is convenient to recast \Bezier projection in terms of a dual basis. 
A dual basis has the distinguishing property that
\begin{align}
  \int_{\Omega} \hat{N}_A N_B \, d\Omega = \delta_{AB}.
  \end{align}
Once a dual basis is defined it can be processed in much the same manner as standard basis functions are processed in a finite element code. A complete exposition on the subject of dual bases and the \Bezier projection framework can be found in~\protect\cite{thomas_bezier_2015}. We first define the dual element extraction operator
\begin{align}
 \hat{\mathbf{D}}^e=\operatorname{diag}({\boldsymbol{\omega}^e})\mathbf{R}^e(\mathbf{G}^e)^{-1}
\end{align}
where $\mathbf{G}^e$ is the Gramian matrix of the Bernstein basis functions over the element and $\operatorname{diag}({\boldsymbol{\omega}^e})$ is a diagonal matrix that contains the B\'ezier projection weights computed by (\protect\ref{eqn:Bezier_weight}). We can then define a dual basis function $\hat{N}_{A(e,a)}$ restricted to element $e$ as
\begin{align}
\hat{N}^e_{a} = \sum_{j}\hat{D}_{aj}^e B_{j}.
\label{eq:local-dual-basis}
\end{align}

The biorthogonality of the dual basis can be seen by noting that
\begin{align}
  \int_{\Omega^e } \hat{\mathbf{N}}^{e} (\mathbf{N}^{e})^T \, d\Omega &= \operatorname{diag}(\boldsymbol{\omega}^e)
\end{align}
and
\begin{align}
\A_{e} \left[\int_{\Omega^e } \hat{\mathbf{N}}^{e} (\mathbf{N}^{e})^T \, d\Omega \right] &= \mathbf{I}
\end{align}
where $\A$ is the standard finite element assembly operator \protect\cite{hughes_finite_2012}.

\begin{figure}
    \centering
    \begin{subfigure}[b]{0.32\linewidth}        
        \centering
        \includegraphics[width=\linewidth]{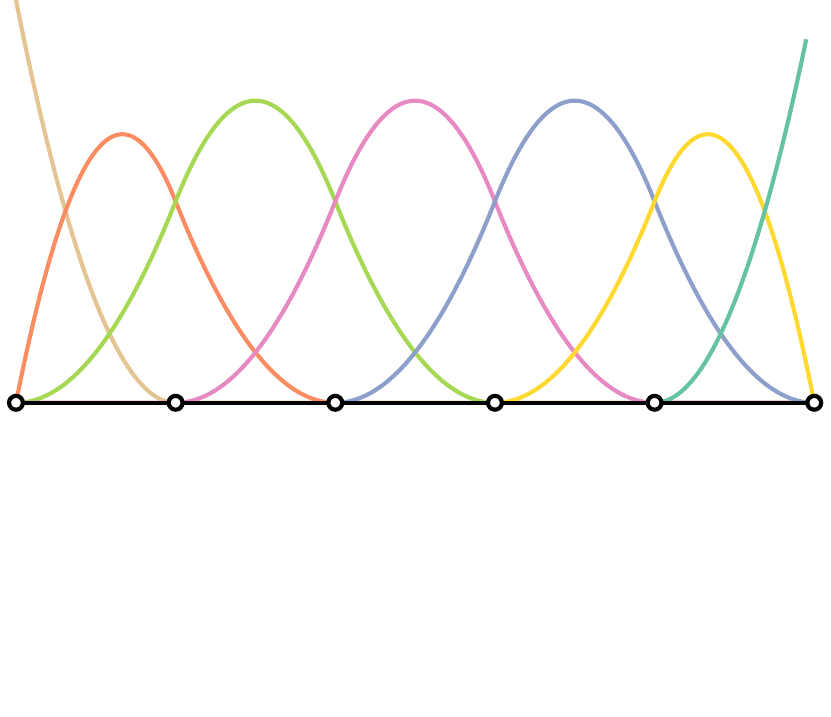}
        \caption{}
    \end{subfigure}
    \begin{subfigure}[b]{0.32\linewidth}        
        \centering
        \includegraphics[width=\linewidth]{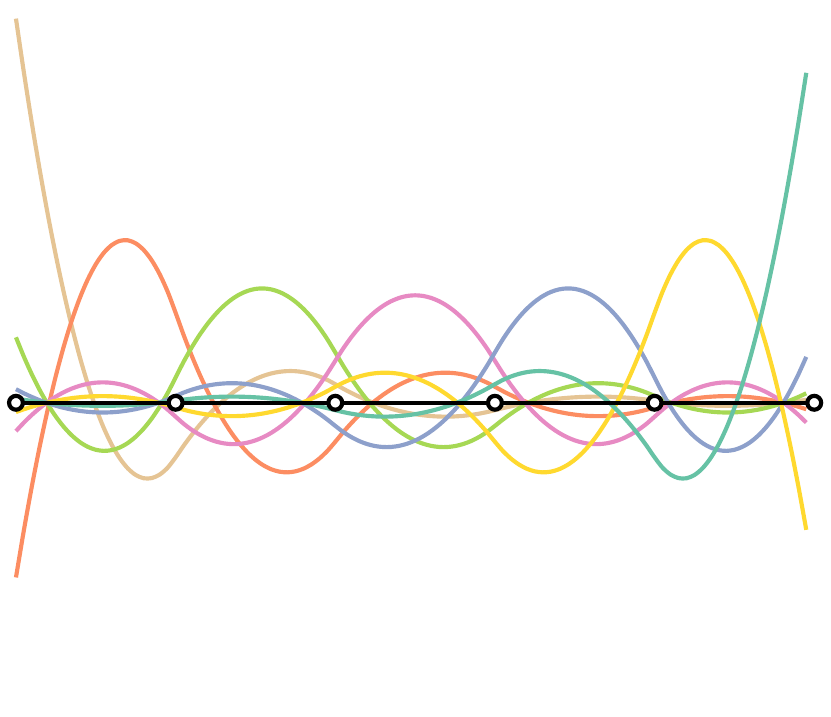}
        \caption{}
    \end{subfigure}
    \begin{subfigure}[b]{0.32\linewidth}        
        \centering
        \includegraphics[width=\linewidth]{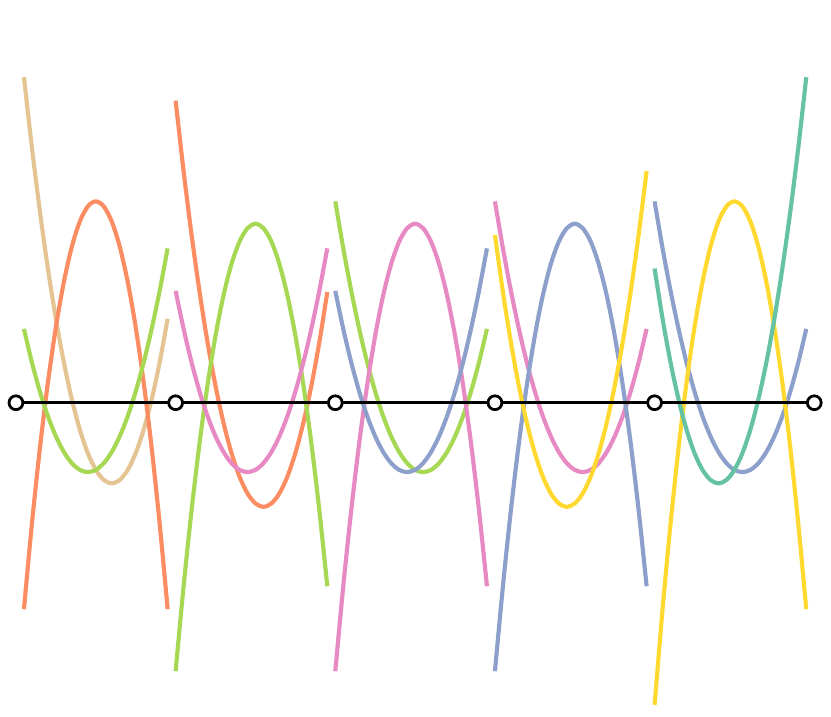}
        \caption{}
    \end{subfigure}
    \caption{A set of quadratic B-spline basis functions (a), and corresponding dual basis functions computed by a global $L^2$ projection (b) and the dual basis functions computed using B\'ezier projection (c). Note that the support of the dual basis functions in (b) is not compact. The dual basis functions shown in (c) have compact support but are not continuous.}
    \label{fig:dualbasis}
\end{figure}

Now, given any function $u \in L^2$ we can use the dual basis to find its representation in terms of the corresponding spline basis as
\begin{align}
  u &= \sum_A P_A N_A
\end{align}
where
\begin{align}
	P_A = \int_{\Omega^A} \hat{N}_{A} u \, d\Omega = \langle \hat{N}_A, u \rangle_{\Omega^A}.
\end{align}
A set of dual basis functions corresponding to the quadratic maximally smooth B-spline basis shown in Figure~\protect\ref{fig:dualbasis}a is shown in Figure~\protect\ref{fig:dualbasis}c. Note that these dual functions have compact support and discontinuities which coincide with the underlying knots in the knot vector. The compact support of the dual basis functions will be crucial for maintaining the sparsity of the stiffness matrix for the B\'{e}zier $\bar{B}$ formulations presented in this paper. For comparison, the dual basis corresponding to \textit{global} $L^2$ projection are shown in Figure~\protect\ref{fig:dualbasis}b. Each of these dual basis functions has global support which explains why the use of global $\bar{B}$ projections results in dense stiffness matrices. 

\subsubsection{Rational dual basis functions}
If rational basis functions are used, the construction of the dual basis must be modified slightly. A rational dual basis must satisfy the biorthogonality requirement
\begin{align}
  \int_{\Omega} \bar{R}_A R_B \, d\Omega= \delta_{AB}.
\end{align}
A simple way to achieve biorthogonality is to define
\begin{align}
  \bar{R}_A &= W \bar{N}_A
\end{align}
where $W$ is the rational weight given in (\protect\ref{eq:weight}). Now
\begin{align}
  \int_{\Omega} \bar{R}_A R_B \, d\Omega &= \int_{\Omega} \bar{N}_A N_B \, d\Omega = \delta_{AB}.
  \end{align}

\section{Geometric locking: Timoshenko beams}
\label{sec:timoshenkobeam}
\begin{figure}[ht]
  \centering
  \includegraphics[width=0.25\linewidth]{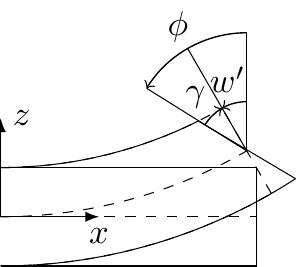}
  \caption{Deformation of a Timoshenko beam. The normal rotates by the angle $\phi$, which is not equal to $w'$, due to shear deformation.}
  \label{fig:Timoshenko_beam_cross}
\end{figure}

To illustrate the use of \Bezier $\bar{B}$ projection to overcome geometric locking effects we study transverse shear locking in Timoshenko beams. The Timoshenko beam problem provides a simple one dimensional setting in which to describe B\'ezier $\bar{B}$ projection. Note however, that the approach can be directly generalized to more complex settings like spatial beams and shells and other geometric locking mechanisms like membrane locking. We consider a planar cantilevered Timoshenko beam as shown in Figure~\protect\ref{fig:Timoshenko_beam_cross}. The strong form for this problem can be stated as
\begin{align}
\begin{rcases}
\begin{aligned}
  -s{GA}\gamma' &=f(x) \\
  -EI\kappa'-s{GA}\gamma &= 0\\
  \kappa &= \phi' \\
  \gamma &= \omega'-\phi
\end{aligned}
\end{rcases}
& \text{ in $\Omega$}\\
\begin{rcases}
\begin{aligned}
\omega&=0\\
\phi&=0
\end{aligned}
\end{rcases}
&\text{ at $x=0$}\\
\begin{rcases}
\begin{aligned}
s{AG}\gamma &= Q\\
-EI\kappa &= M
\end{aligned}
\end{rcases}
&\text{ at $x=L$}
\end{align}
where $\gamma$ is the shear strain, $\kappa$ is the bending strain, $\omega$ is the vertical displacement, $\phi$ is the angle of rotation of the normal to the mid-plane of the beam, $f$ is the distributed transverse load, $Q$ is a point load, $M$ is the moment, $E$ is the Young's modulus, $G$ is the shear modulus, $A$ is the cross-sectional area, $I$ is the second moment of inertia of the beam cross-section, $s$ is the shear correction factor, normally set to $5/6$ for rectangular cross-sections, and $\Omega = (0,L)$. When $\omega$ and $\phi$ are interpolated by basis functions of the same order the finite element solution to this problem exhibits shear locking as the beam becomes slender.

\subsection{Symmetric \Bezier $\bar{B}$ projection}
\label{sec:symmetric-projection}
{Locking is caused in Timoshenko beam problems when equal order interpolation is used for both midline displacements and rotations. To reduce the effects of locking for these problems the $\bar{B}$ method projects the shear strain, $\gamma$, onto a lower order function space. This process produces a projected shear strain, which we call $\bar{\gamma}$,  that is substituted into the weak form of the problem statement. In this section we use the \Bezier projection operator to compute an approximation of $\bar{\gamma}$, and refer to this approach as the symmetric \Bezier $\bar{B}$ projection method because the resulting stiffness matrix is symmetric.}

\subsubsection{The weak form}

\sloppy Given the function spaces $\mathcal{S}(\Omega)=\{{\mathbf{u} \, \vert \, {\mathbf{u}\in{H^1(\Omega)\times{H^1(\Omega)}}},\mathbf{u}\vert_{\Gamma_{g}}=\mathbf{g}}\}$ and $\mathcal{V}(\Omega)=\{{\mathbf{w} \, \vert \, {\mathbf{w}\in{H^1(\Omega)\times{H^1(\Omega)}}}, \allowbreak \mathbf{w}\vert_{\Gamma_{g}}=\mathbf{0}}\}$ where $\mathbf{u}=\left\{{\omega,\phi}\right\}^T$, $\mathbf{w}=\left\{{\delta\omega,\delta\phi}\right\}^T$, $\mathbf{g}$ is the prescribed Dirichlet boundary condition, and $\Gamma_g$ is the Dirichlet boundary at $x=0$, the weak form of the problem can be stated as: find $\mathbf{u}\in{\mathcal{S}(\Omega)}$ such that for all $\mathbf{w}\in{\mathcal{V}(\Omega)}$
\begin{equation}
    {\int_{0}^L\kappa(\mathbf{w})EI\kappa(\mathbf{u}) + \bar{\gamma}(\mathbf{w})sGA\bar{\gamma}(\mathbf{u})} \, dx=\int_0^L\delta\omega f dx+\delta\omega(L)Q+\delta\phi(L)M.
\end{equation}

\subsubsection{Discretization}
{We discretize $\omega$ and $\phi$ as}
\begin{align}
  {\omega = \sum_A \omega_A N_A} \\
  {\phi = \sum_A \phi N_A}
\end{align}
where {$N_A$ is a degree $p$ spline basis function and $\omega_A$ and $\phi_A$ are the corresponding control point values. The shear strain and bending strain can then be expressed as}
{\begin{align}
\gamma &= \sum_A{\begin{bmatrix}N'_A & -N_A\end{bmatrix}\begin{bmatrix}\omega_A & \phi_A\end{bmatrix}^T} \\
\kappa &= \sum_A\begin{bmatrix}0 & N'_A\end{bmatrix}\begin{bmatrix}\omega_A & \phi_A\end{bmatrix}^T
\end{align}}
The shear strain $\bar{\gamma}$ is constructed by B\'{e}zier projection of the true shear strain $\gamma$ onto a lower degree space. In other words, we project from a $p^{th}$ degree spline space with $n$ basis functions $\mathbf{N}$ defined by the knot vector
\begin{align}
\mathbf{\Xi}_{p}=\lbrace{\underbrace{0,0,\ldots,0}_\text{$p+1$ copies}},\mathbf{\Xi}_{int},{\underbrace{1,1,\ldots,1}_\text{$p+1$ copies}}\rbrace,
\label{eq:origin_knot_vector}
\end{align}
onto a $p-1^{th}$ degree spline space with $\bar{n}$ basis functions $\bar{\mathbf{N}}$ defined by the knot vector
\begin{equation}
\bar{\mathbf{\Xi}}_{p-1}=\lbrace{\underbrace{0,0,\ldots,0}_\text{$p$ copies}},\mathbf{\Xi}_{int},{\underbrace{1,1,\ldots,1}_\text{$p$ copies}}\rbrace
\label{eq:projected_knot_vector}
\end{equation}
where the internal knots, denoted by $\mathbf{\Xi}_{int}$, are the same for both spaces. The projected shear strain $\bar{\gamma}$ can then be written as
\begin{align}
  \bar{\gamma} &= \sum_A \bar{\gamma}_A \bar{N}_A.
\end{align}
The control variables $\bar{\gamma}_A$ are simply
\begin{align}
  \bar{\gamma}_A = \int_{\Omega^A} \hat{\bar{N}}_A \gamma \, d\Omega = \langle \hat{\bar{N}}_A, \gamma \rangle_{\Omega^i}
\end{align}
where $\hat{\bar{N}}_A$ is a dual basis function for the spline space of degree $p-1$ computed from \eqref{eq:local-dual-basis}.

Localizing to the B\'{e}zier element we define the strain-displacement arrays in terms of element Bernstein basis functions of degree $p$ and $p-1$ as
\begin{align}
	\mathbf{B}_e^\kappa 
    	&= \begin{bmatrix} 0 & -{B^e_{0,p}}' & \cdots & 0 & -{B^e_{p,p}}' \end{bmatrix}, \\
    \mathbf{B}_e^\gamma
    	&= \begin{bmatrix} {B^e_{0,p}}' & -B^e_{0,p} & \cdots & {B^e_{p,p}}' & -B^e_{p,p} \end{bmatrix}, \\
    \bar{\mathbf{B}}_e
    	&= \begin{bmatrix} \bar{B}^e_{0,p-1} & \cdots & \bar{B}^e_{p-1,p-1} \end{bmatrix},
\end{align}
{where $B^e_{i,p}$ is the $i^{th}$ Bernstein basis function of order $p$}. We can then compute the element arrays as
\begin{align}
\mathbf{K}_e^b &= EI\mathbf{C}^e\langle{{\mathbf{B}_e^\kappa}^T,\mathbf{B}_e^\kappa}\rangle(\mathbf{C}^e)^T, \label{eq:symmetric_timoshenko}\\
\bar{\mathbf{M}}_e &= s{GA}\bar{\mathbf{C}}^e\langle{\bar{\mathbf{B}}_e^T,\bar{\mathbf{B}}_e}\rangle(\bar{\mathbf{C}}^e)^T, \\
{\hat{\mathbf{P}}_e} &= \langle{(\hat{\bar{\mathbf{N}}}^e)^T,\mathbf{B}^\gamma_e}\rangle(\mathbf{C}^e)^T \label{eq:P-hat-elem},
\end{align}
where $\mathbf{C}^e$ is the element extraction operator for the degree $p$ spline space, $\bar{\mathbf{C}}^e$ is the element extraction operator for the degree $p-1$ spline space, and $\hat{\bar{\mathbf{N}}}^e$ are the dual basis functions restricted to the element for the degree $p-1$ spline space. The global stiffness matrix can then be written as
\begin{align}
	\mathbf{K} = \mathbf{K}^b + \bar{\mathbf{K}}^s_S \label{eqn:symmetric_b_bar_timoshenko_stiffness}
\end{align}
where
\begin{align}
	\mathbf{K}^b &=\A_e\mathbf{K}_e^b,\label{eqn:bending_stiffness} \\
	\bar{\mathbf{K}}^s&={\hat{\mathbf{P}}}^T\bar{\mathbf{M}}{\hat{\mathbf{P}}}\label{eqn:shear_stiffness} \\
	{\hat{\mathbf{P}}} &=\A_e{\hat{\mathbf{P}}_e} \label{eq:P-hat-global}\\
    \bar{\mathbf{M}} &=\A_e\bar{\mathbf{M}}_e
\end{align}
and $\A$ is the standard finite element assembly operator \protect\cite{hughes_finite_2012}. We note that the assembly of $\bar{\mathbf{K}}^s$ requires the assembly of two intermediate matrices, $\bar{\mathbf{M}}$ and ${\hat{\mathbf{P}}}$. The computation of these matrices is needed because the product of two integrals over the entire domain can not be localized to the element level.
\subsection{Non-symmetric \Bezier $\bar{B}$ projection}
\label{sec:non-symmetric}

{Simo and Hughes \protect\cite{hughes_variational_1986} have shown that $\bar{B}$ formulations and mixed formulations are equivalent. However, the development of the symmetric \Bezier $\bar{B}$ projection method presented in Section~\protect\ref{sec:symmetric-projection}, where we began by interpreting the $\bar{B}$ formulation as a strain projection method, lacks a connection to mixed formulations. In this section, we present a second method based on \Bezier projection in which we view the $\bar{B}$ formulation as a mixed formulation where, for the Timoshenko beam problem, the auxiliary variable is the shear strain.  We use the \Bezier dual basis functions as the test functions for the auxiliary variable. Once the problem has been cast as a mixed formulation we then eliminate the auxiliary variable to get a purely displacement based formulation. This approach preserves convergence rates and all assembly routines for the stiffness matrix can be performed at the element level. However, it does not produce a symmetric stiffness matrix.}

\subsubsection{The weak form}
{In the mixed formulation for the Timoshencko beam problem the shear stress, $\tau=sGA\gamma$, is taken as a new independent variable. The weak form of the mixed formulation can then be stated as: find $\mathbf{u}\in{\mathcal{S}(\Omega)}$ and $\tau\in{L}^2(\Omega)$ such that for all $\mathbf{w}\in{\mathcal{V}(\Omega)}$ and $\delta\tau\in{L}^2(\Omega)$}
\begin{align}
    \int_{0}^L\kappa(\mathbf{w})EI\kappa(\mathbf{u}) + {\gamma}(\mathbf{w})\tau \, dx&= l\langle{\mathbf{w}}\rangle\\
    \int_{0}^L\delta\tau(sGA\gamma(\mathbf{u})-\tau) \, dx &= 0.
\end{align}
\subsubsection{Discretization}
{In the finite element formulation of the mixed problem, the discretization of $\mathbf{u}$ and $\mathbf{w}$ remain the same as before. The shear strain and its variation, however, are in $L^2(\Omega)$, so their finite element approximation can consist of functions with lower regularity, such as discontinuous polynomials. When the field $\mathbf{u}$ is discretized by $p^{th}$ degree spline basis functions defined by the knot vector given in \eqref{eq:origin_knot_vector}, we use the $p-1$ degree spline basis functions, defined by the knot vector given in \eqref{eq:projected_knot_vector}, to discretize the shear stress $\tau$ and the corresponding dual basis to discretize its variation $\delta\tau$. The discrete form of the shear stress and its variation are given by}
\begin{align}
    \tau&=\sum_A\tau_A\bar{N}_A \\
    \delta\tau&=\sum_A\delta\tau_A\hat{\bar{N}}_A.
\end{align}
{The stiffness matrix for the mixed form can then be written as}
\begin{equation}
   \mathbf{K}_{mix} =  
    \begin{bmatrix}
        \mathbf{K}^b & \mathbf{P}^T \\
        sGA\hat{\mathbf{P}} & -\mathbf{I}
    \end{bmatrix}
    \label{eq:mixed_timoshenko}
\end{equation}
{where}
\begin{equation}
    \mathbf{P}=\A_e\mathbf{P}_e,
\end{equation}
{and}
\begin{equation}
    \mathbf{P}_e=\bar{\mathbf{C}}^e\langle{\bar{\mathbf{B}}_e^T,\mathbf{B}^\gamma_e}\rangle(\mathbf{C}^e)^T.
\end{equation}
{ and $\hat{\mathbf{P}}$ is given in \eqref{eq:P-hat-elem} and \eqref{eq:P-hat-global}.
We can now eliminate the control variable of the shear stress from \eqref{eq:mixed_timoshenko} to get a pure displacement formulation where the stiffness matrix can be written as}
\begin{equation}
    \mathbf{K}=\mathbf{K}^b+\bar{\mathbf{K}}^s_{NS}, 
\end{equation}
where 
\begin{equation}
    \bar{\mathbf{K}}^s_{NS}=sGA\mathbf{P}^T\hat{\mathbf{P}}.
    \label{eq:nonsymmetric_timoshenko}
\end{equation}
{We can see that the use of different function spaces for the shear strain and its variation leads to a non-symmetric stiffness matrix.}

\paragraph{Remark} { As mentioned previously, the symmetric \Bezier $\bar{B}$ formulation is not consistent with a mixed formulation. To see this, we can recover the mixed formulation of \eqref{eqn:symmetric_b_bar_timoshenko_stiffness}, which is}
\begin{equation}
    \begin{bmatrix}
        \mathbf{K}^b & \hat{\mathbf{P}}^T \\
        sGA\hat{\mathbf{P}} & -\bar{\mathbf{M}}^{-1}
    \end{bmatrix},
\end{equation}
{where both the shear stress and its variation are discretized by the dual basis functions. However, for the inner product of dual basis functions we have}
\begin{equation}
    \langle{\hat{\bar{N}}_i,\hat{\bar{N}}_j}\rangle\neq{\bar{\mathbf{M}}^{-1}_{ij}},
\end{equation}
{which shows the inconsistency between the symmetric \Bezier $\bar{B}$ formulation and the mixed formulation.}

\subsection{Bandwidth of the stiffness matrix}
\begin{figure}[!htb]
    \centering
    \begin{subfigure}[b]{0.24\linewidth}        
        \centering
        \includegraphics[width=\linewidth]{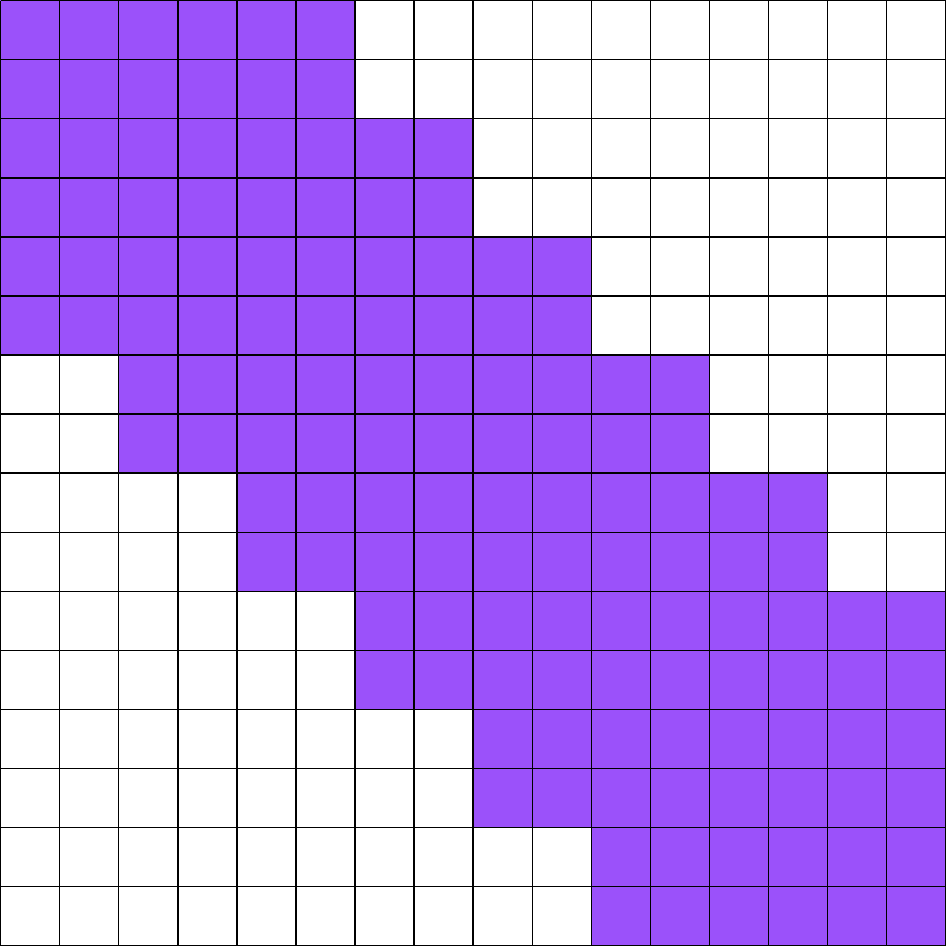}
        \caption{Standard}
    \end{subfigure}
    \begin{subfigure}[b]{0.24\linewidth}        
        \centering
        \includegraphics[width=\linewidth]{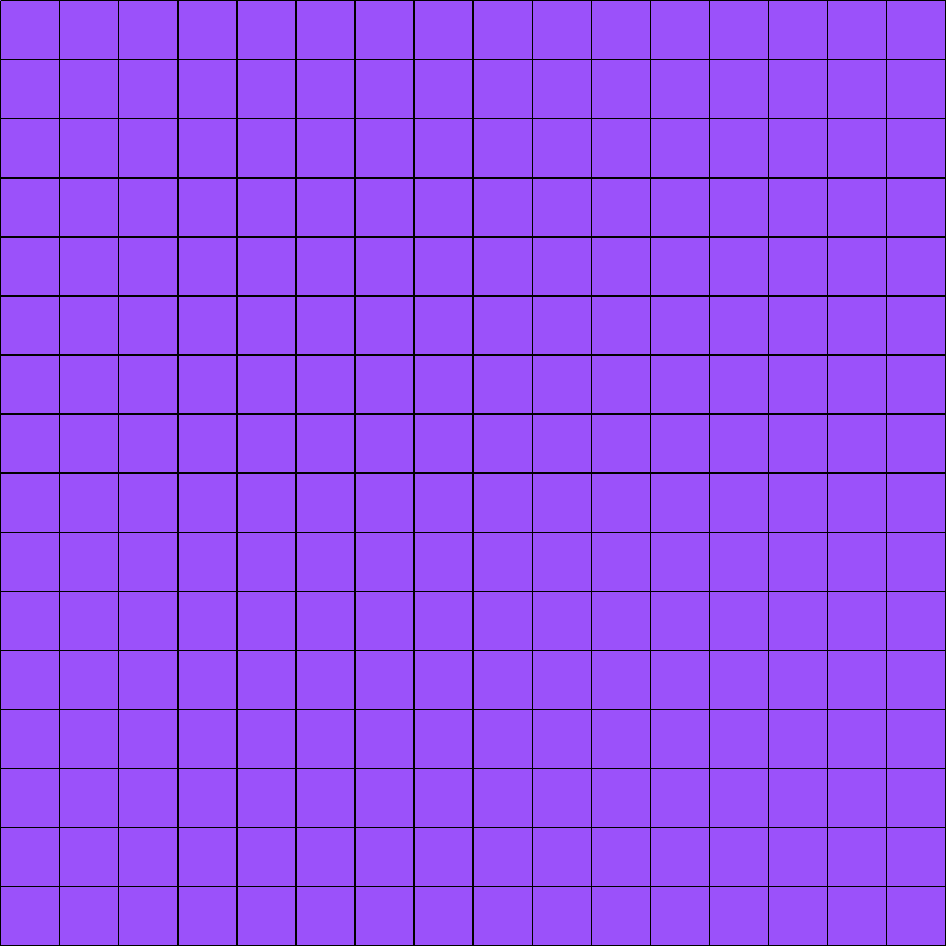}
        \caption{Global $\bar{B}$}
    \end{subfigure}
    \begin{subfigure}[b]{0.24\linewidth}        
        \centering
        \includegraphics[width=\linewidth]{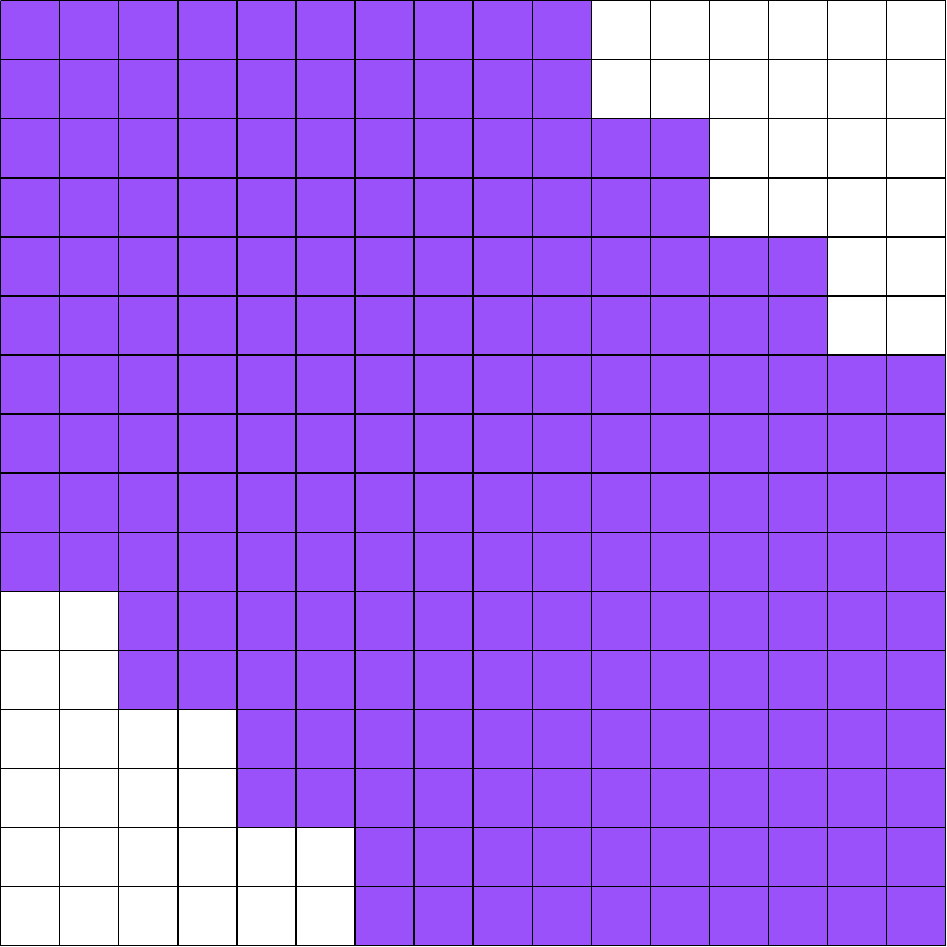}
        \caption{Symmetric B\'ezier $\bar{B}$}
    \end{subfigure}
    \begin{subfigure}[b]{0.24\linewidth}        
        \centering
        \includegraphics[width=\linewidth]{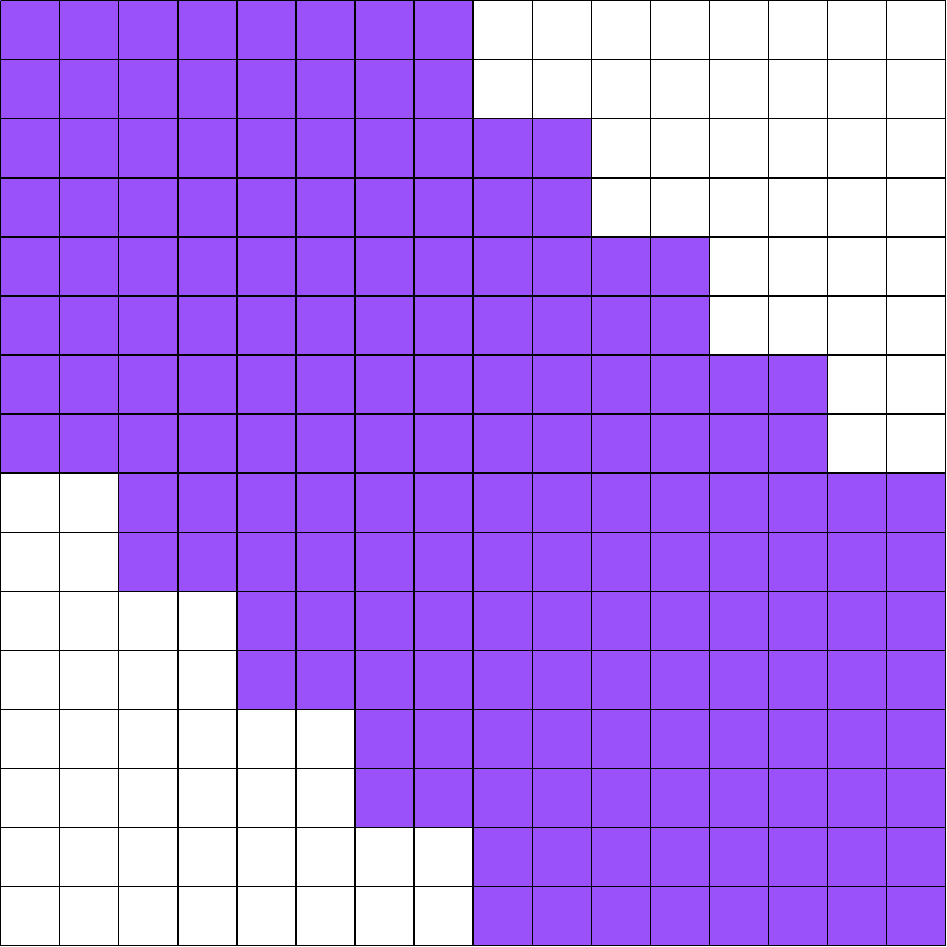}
        \caption{Non-symmetric B\'ezier $\bar{B}$}
    \end{subfigure}
    \caption{Illustrations of the structure of $2^{nd}$ order Timoshenko beam stiffness matrices for (a) a standard displacement method, (b) a global $\bar{B}$ method, (c) a symmetric B\'ezier $\bar{B}$ method and (d) a non-symmetric B\'ezier $\bar{B}$ method.}
    \label{fig:stiffness_matrix}
\end{figure}

\begin{figure}[htb!]
      \centering
      \includegraphics[width=.5\linewidth]{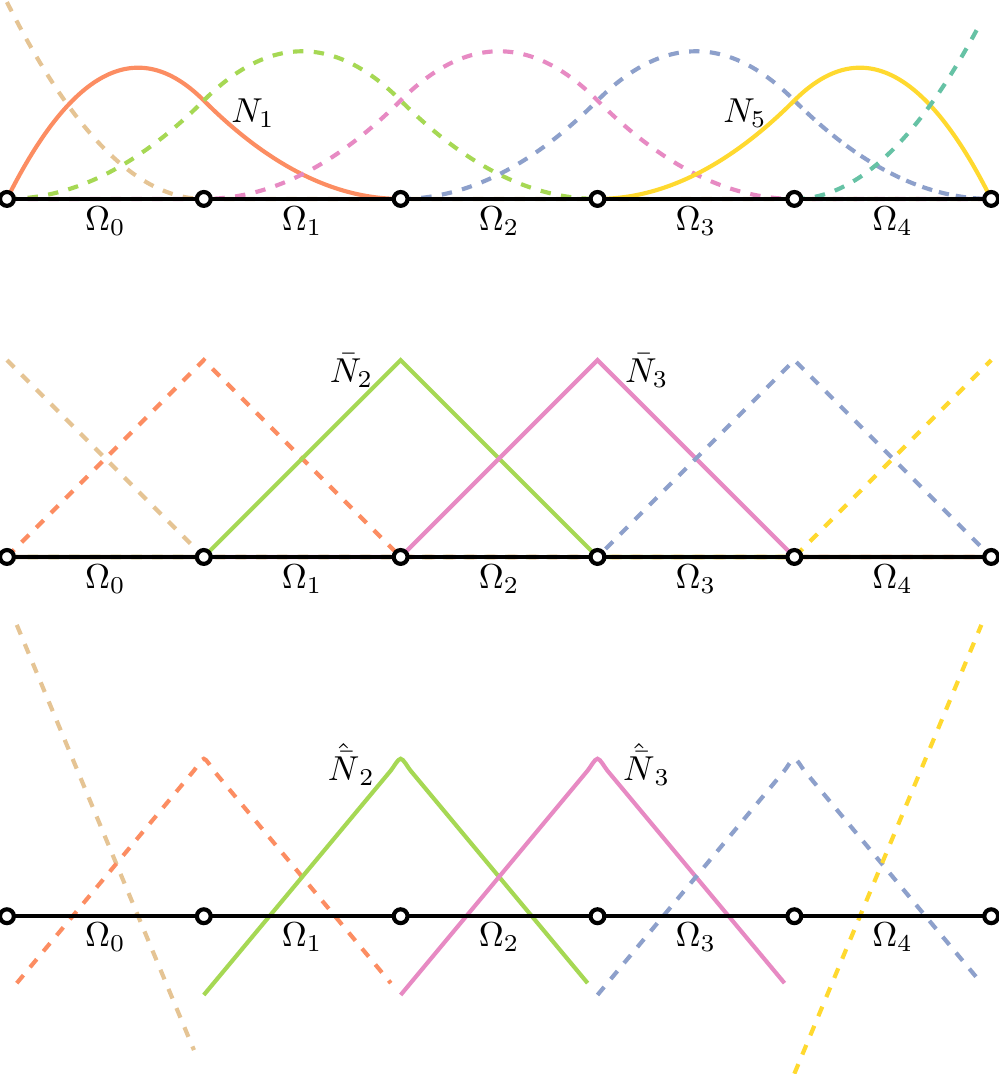}
\caption{Quadratic maximally smooth B-spline basis functions (top), associated linear basis functions (middle), and dual basis functions (bottom) for the B\'ezier $\bar{B}$ formulation.}
\label{fig:IGAelement}
\end{figure}

A global $\bar{B}$ method that utilizes a global $L^2$ projection results in a dense stiffness matrix. The B\'ezier $\bar{B}$ {methods}, on the other hand, {produce} sparse stiffness matrices. However, the coupling of the local dual basis functions does increase the bandwidth slightly. This is illustrated in Figure~\protect\ref{fig:stiffness_matrix}, which shows the structure of the stiffness matrix for the Timoshenko beam problem using the second order basis functions of maximal smoothness for a displacement-based method (Figure~\protect\ref{fig:stiffness_matrix}a), global $\bar{B}$ method (Figure~\protect\ref{fig:stiffness_matrix}b), symmetric B\'ezier $\bar{B}$ method (Figure~\protect\ref{fig:stiffness_matrix}c) and non-symmetric \Bezier $\bar{B}$ method (Figure~\protect\ref{fig:stiffness_matrix}d). The blank cells indicate zero terms in the matrix while colored cells show the location of nonzero terms.
∏
The increased bandwidth of the {symmetric} B\'ezier $\bar{B}$ method when compared to a displacement-based method can be explained by looking at the product of the integrals in \eqref{eqn:shear_stiffness}. For example, if we consider the basis functions $N_{1}$ and $N_{5}$ in Figure~\protect\ref{fig:IGAelement} we see that $\operatorname{supp}(N_{1}) \cap{} \operatorname{supp}(N_{5}) = \emptyset$, which means that the inner product of these two functions will be zero and the corresponding coefficient in the stiffness matrix will be zero in the displacement-based method. For the B\'ezier $\bar{B}$ method, however, the form of \eqref{eqn:shear_stiffness} leads to a coupling between $N_{1}$ and $N_{5}$. This can be seen by considering $\Omega_2$. Over this element, the shear stiffness can be represented as
\begin{align}
	\bar{\mathbf{K}}^s_2 = \sum_{i=1}^3 \sum_{j=1}^3 \mathbf{P}_i^T\bar{\mathbf{M}}_2\mathbf{P}_j
\end{align}
and the term of this summation that results in the coupling between $N_{1}$ and $N_{5}$ is $\mathbf{P}_1^T\bar{\mathbf{M}}_2\mathbf{P}_3$, where $\mathbf{P}_1$ is the inner product of $N_{1}$ and $\hat{\bar{N}}_{2}$, $\mathbf{P}_3$ is the inner product of $N_{5}$ and $\hat{\bar{N}}_{3}$, and $\bar{\mathbf{M}}_2$ is the inner product of $\bar{N}_{2}$ and $\bar{N}_{3}$. We can see from Figure~\protect\ref{fig:IGAelement} that $\operatorname{supp}(N_{1})\cap{}\operatorname{supp}(\hat{\bar{N}}_{2})=\Omega_1$, $\operatorname{supp}(N_{5})\cap{}\operatorname{supp}(\hat{\bar{N}}_{3})=\Omega_3$ and $\operatorname{supp}(\bar{N}_{2})\cap{}\operatorname{supp}(\bar{N}_{3})=\Omega_2$, so that $\mathbf{P}_1^T\bar{\mathbf{M}}_2\mathbf{P}_3$ is not zero. Thus we have increased the number of nonzero coefficients in the shear stiffness matrix. However, the same exercise can be used to show that there is no coupling between $N_{0}$ and $N_{6}$ for this set of basis functions so the matrix is not dense. {the bandwidth of the non-symmetric B\'ezier $\bar{B}$ method is reduced further. This is because the Gramian matrix does not appear in the this formulation.} In fact, from the formulation of the element stiffness matrix, we can show that the bandwidth of the stiffness matrix of the symmetric B\'ezier $\bar{B}$ and non-symmetric B\'ezier $\bar{B}$ methods for the Timoshenko beam are $6p-3$ and $4p-1$, respectively. 

\paragraph{Remark} In \protect\cite{bouclier_efficient_2013} a local $\bar{B}$ method for shells was proposed that was based on the local least squares method presented in \protect\cite{govindjee_convergence_2012}. This approach has a similar structure to the {symmetric B\'ezier $\bar{B}$} method presented here. However, it was shown in~\protect\cite{thomas_bezier_2015} that choosing (\protect\ref{eqn:Bezier_weight}) as the weighting provides a significant increase in the accuracy of the approximation.

\subsection{Numerical results}

\begin{figure}[h]
	\centering
	\includegraphics[width=0.4\linewidth]{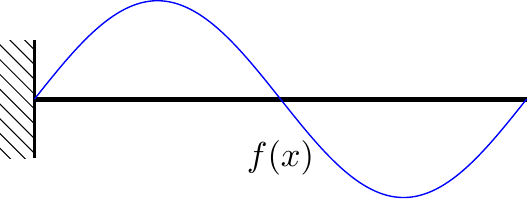}
	\caption{Straight planar cantilevered Timoshenko beam clamped at the left and loaded by a distributed load $f(x)$.}
	\label{fig:load_beam}
\end{figure}

In our study, a straight planar cantilever beam is clamped on the left end and a sinusoidal distributed load $f(x)=sin(\pi\dfrac{x}{l})$ is applied, as depicted in Figure \protect\ref{fig:load_beam}. The analytical solution for vertical displacement $w$, rotation $\phi$, bending moment $M$, and transverse shear force $Q$ are given by
\begin{align}
\begin{split}
w(x)&=\frac{EI\left(6 \pi ^2 l^2 \sin \left(\frac{\pi  x}{l}\right)+6 \pi ^3 l x\right)+sGA\left(6 l^4 \sin \left(\frac{\pi  x}{l}\right)-6 \pi l^3 x+3 \pi ^3 l^2 x^2-\pi ^3 l x^3\right)}{6 \pi ^4 sEIGA}\\
\phi(x)&=\frac{2 l^3 \cos \left(\frac{\pi  x}{l}\right)-2 l^3+2 \pi ^2 l^2 x-\pi ^2 l x^2}{2 \pi ^3 EI}\\
M(x)&=\frac{l^2 \sin \left(\frac{\pi  x}{l}\right)-\pi  l^2+\pi  l x}{\pi ^2}\\
Q(x)&=\frac{-l\cos \left(\frac{\pi  x}{l}\right)-l}{\pi}.
\end{split}
\label{eq:timoshenko_analytical}
\end{align}

\begin{figure}
    \centering
    \begin{subfigure}[b]{0.49\linewidth}        
        \centering
        \includegraphics[width=\linewidth]{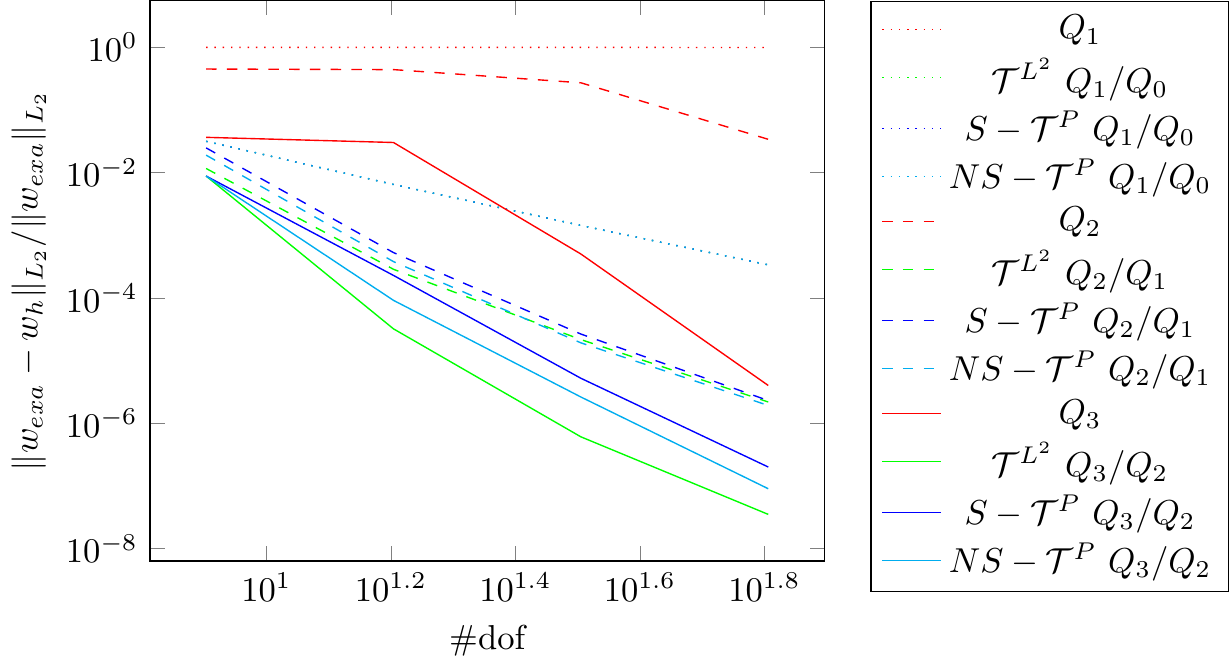}
        \caption{}
        \vspace*{2mm}
    \end{subfigure}
    \begin{subfigure}[b]{0.49\linewidth}        
        \centering
        \includegraphics[width=\linewidth]{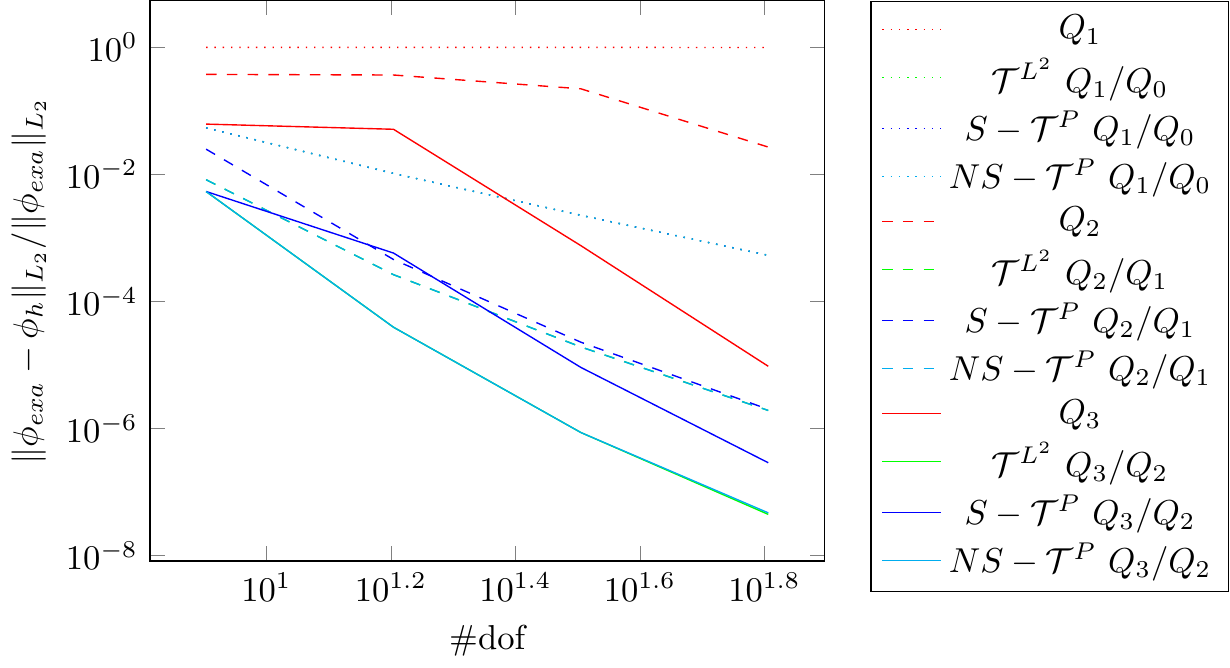}
        \caption{}
        \vspace*{2mm}
    \end{subfigure}
        \begin{subfigure}[b]{0.49\linewidth}        
        \centering
        \includegraphics[width=\linewidth]{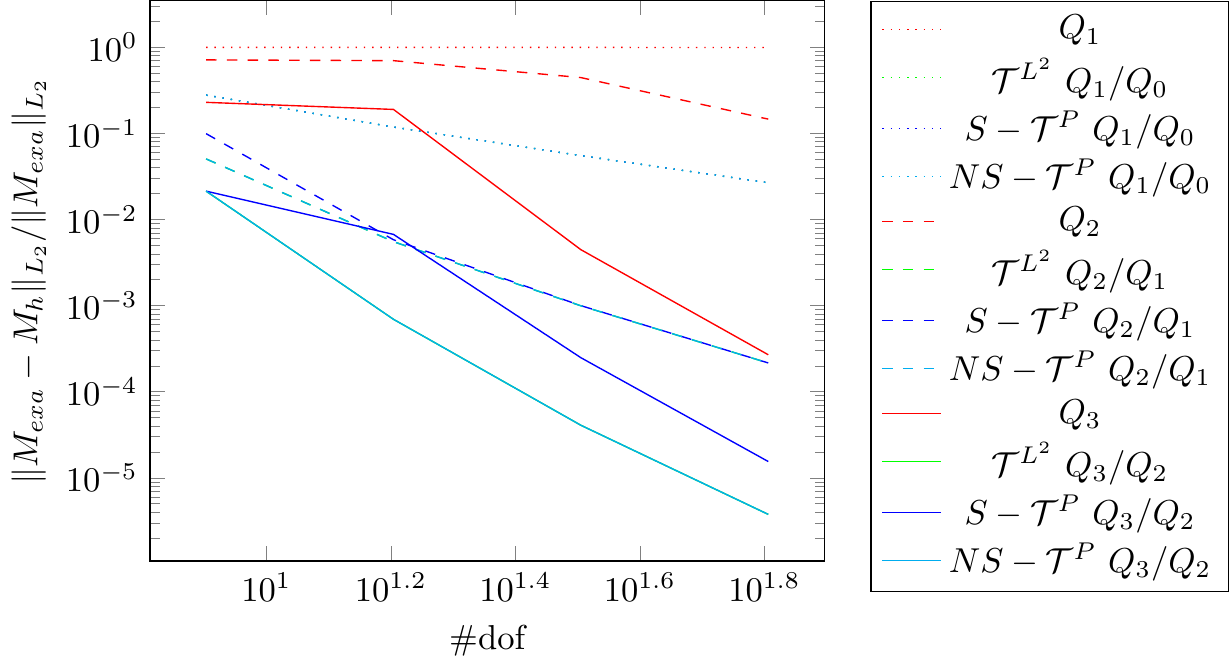}
        \caption{}
        \vspace*{2mm}
    \end{subfigure}
    \begin{subfigure}[b]{0.49\linewidth}        
        \centering
        \includegraphics[width=\linewidth]{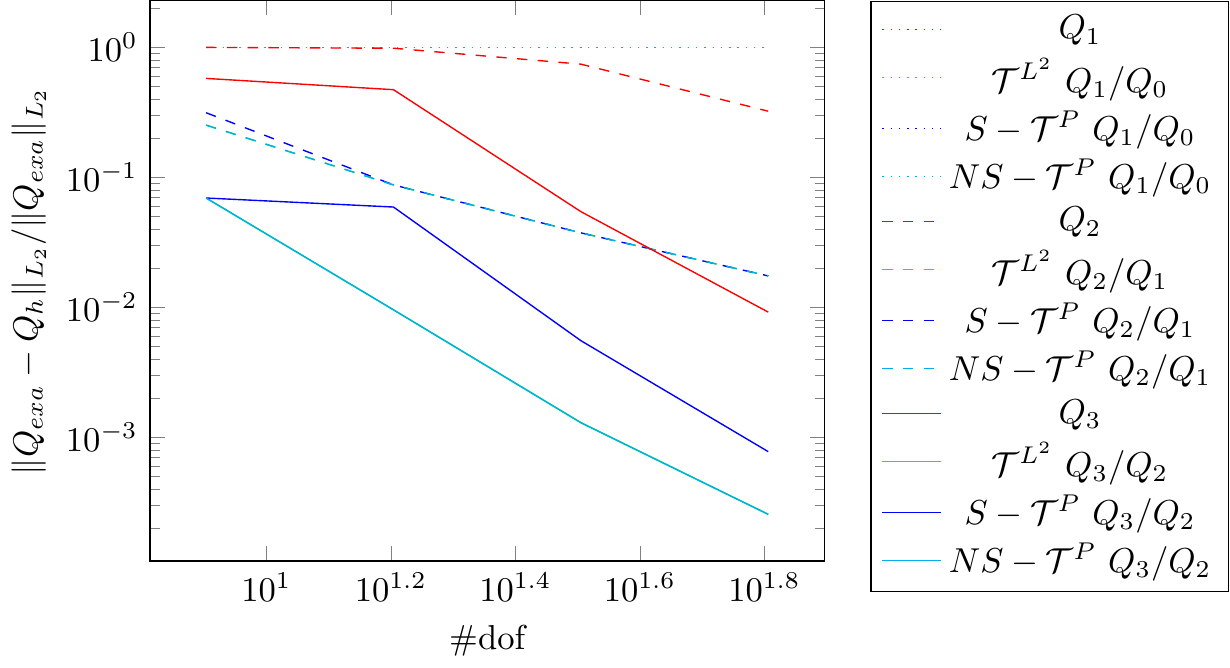}
        \caption{}
        \vspace*{2mm}
    \end{subfigure}
    \caption{Convergence studies for slenderness factor $l/t=10^{-3}$. Error in the $L^2$-norm for (a) displacement $w$, (b) rotation $\phi$, (c) bending moment $M$, and (d) shear force $Q$.}
    \label{fig:timoshenko_result}
\end{figure}

The beam has a rectangular cross-section and we use the following non-dimensional sectional and material parameters: length $l=10$, width $b=1$, thickness $t=0.01$, Young's modulus $E=10^9$, Poisson's ratio $\nu=0.3$, and a shear correction factor of $s=5/6$. A comparison of the normalized error in the $L^2$ norm for $w$, $\phi$, $M$ and $Q$ versus the number of degrees of freedom for polynomial degrees $p=1,2,3$ is shown in Figure \protect\ref{fig:timoshenko_result}. Results computed using standard finite elements are labeled $Q_1$, $Q_2$, and $Q_3$. Results computed using a global $\bar{B}$ method are labeled $\mathcal{T}^{L^2}$, those computed with the symmetric B\'ezier $\bar{B}$ method and the non-symmetric B\'ezier $\bar{B}$ method are labeled $S-\mathcal{T}^{P}$ and $NS-\mathcal{T}^{P}$, respectively. As expected, the $Q_1$ results lock and the error remains virtually unchanged as the mesh is refined. Increasing the polynomial degree does reduce the locking effect, although the reduction is minor for the $Q_2$ results. $\mathcal{T}^{L^2}$, $S-\mathcal{T}^{P}$ and $NS-\mathcal{T}^P$ are essentially locking free for all polynomial orders. The convergence rates for the $\bar{B}$ methods are at least $p+1$ for $w$, $p$ for $\phi$, $p-1$ for $M$, and $p-2$ for $Q$. These rates agree with those reported in \protect\cite{kiendl_single-variable_2015} and are optimal. To reiterate, B\'{e}zier $\bar{B}$ methods produces the same convergence rates as the global $\bar{B}$ method {and the the error plots of $NS-\mathcal{T}^{P}$ for $\phi$, $M$ and $Q$ are identical to those of $\mathcal{T}^{L^2}$}.

\begin{figure}
    \centering
    \begin{subfigure}[b]{0.49\linewidth}        
        \centering
        \includegraphics[width=\linewidth]{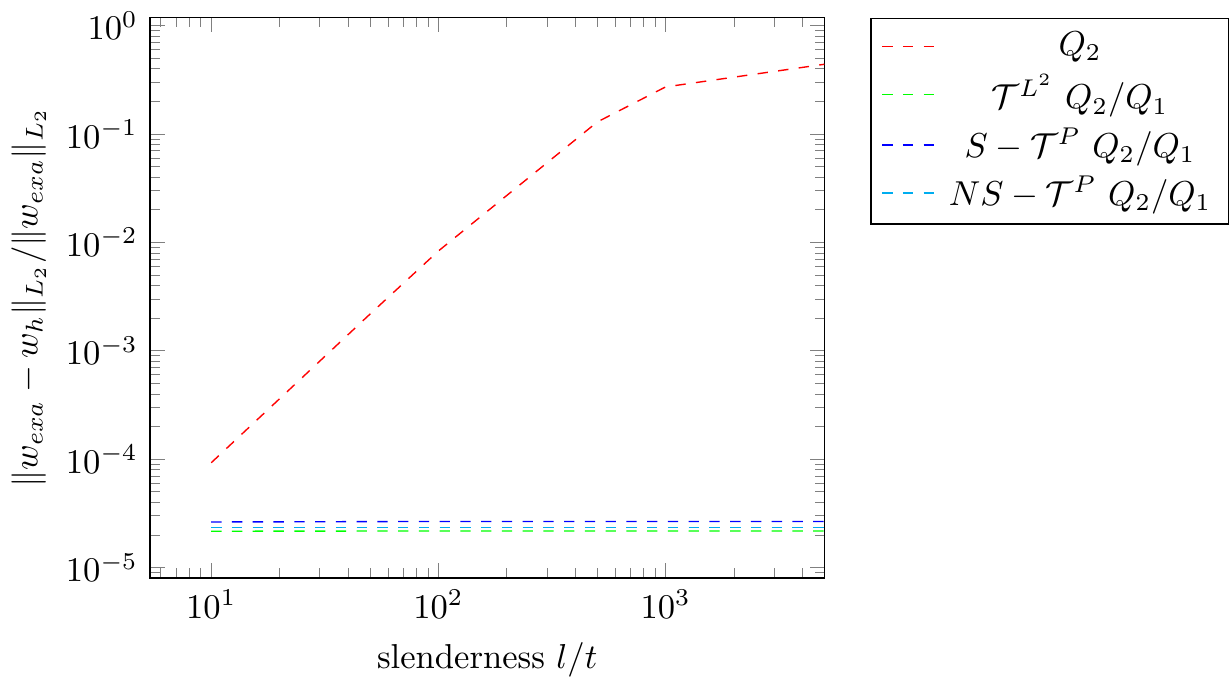}
        \caption{}
        \vspace*{2mm}
    \end{subfigure}
    \begin{subfigure}[b]{0.49\linewidth}        
        \centering
        \includegraphics[width=\linewidth]{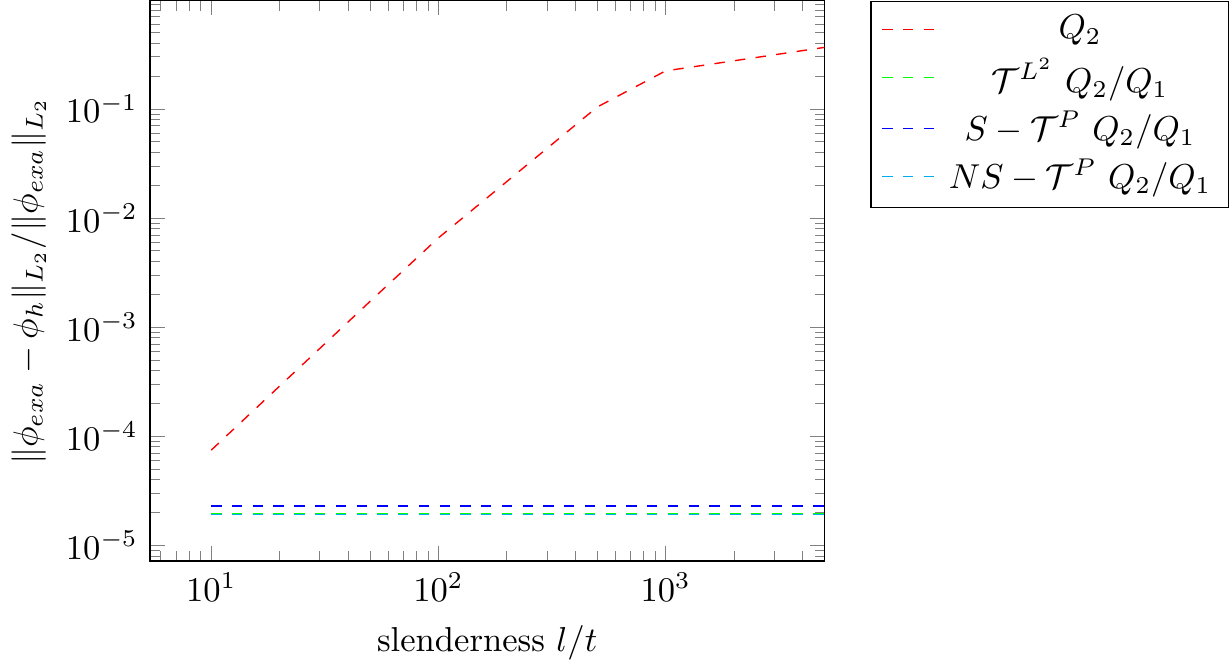}
        \caption{}
        \vspace*{2mm}
    \end{subfigure}
        \begin{subfigure}[b]{0.49\linewidth}        
        \centering
        \includegraphics[width=\linewidth]{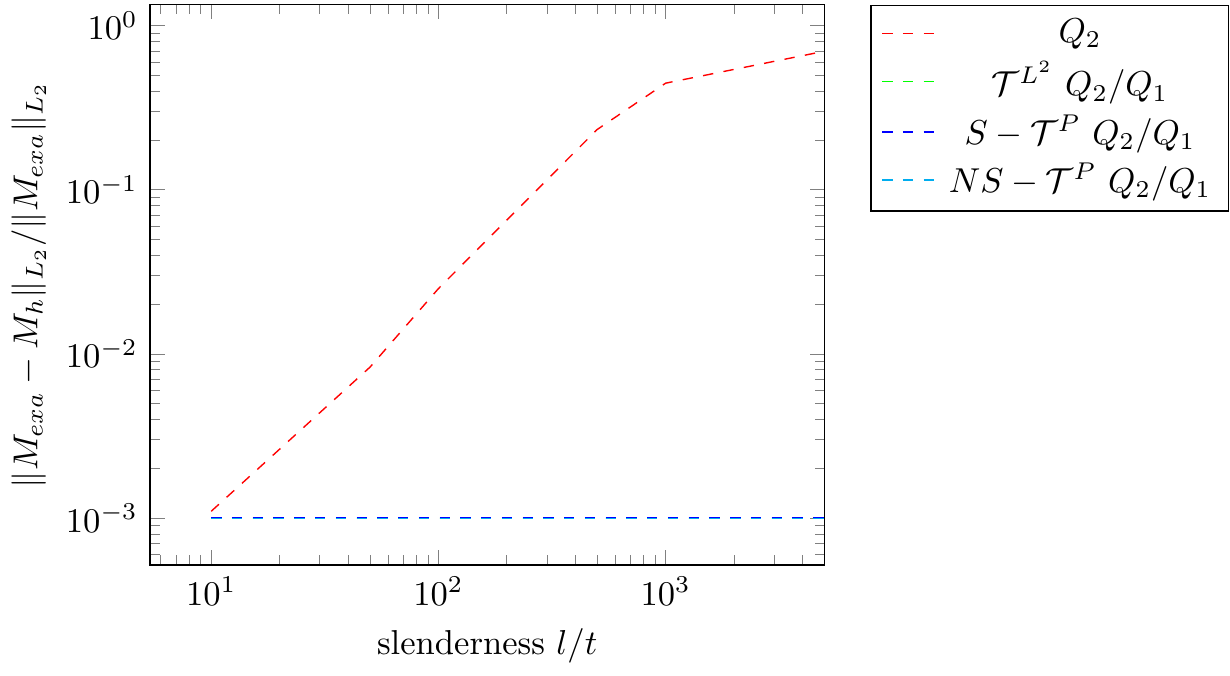}
        \caption{}
        \vspace*{2mm}
    \end{subfigure}
    \begin{subfigure}[b]{0.49\linewidth}        
        \centering
        \includegraphics[width=\linewidth]{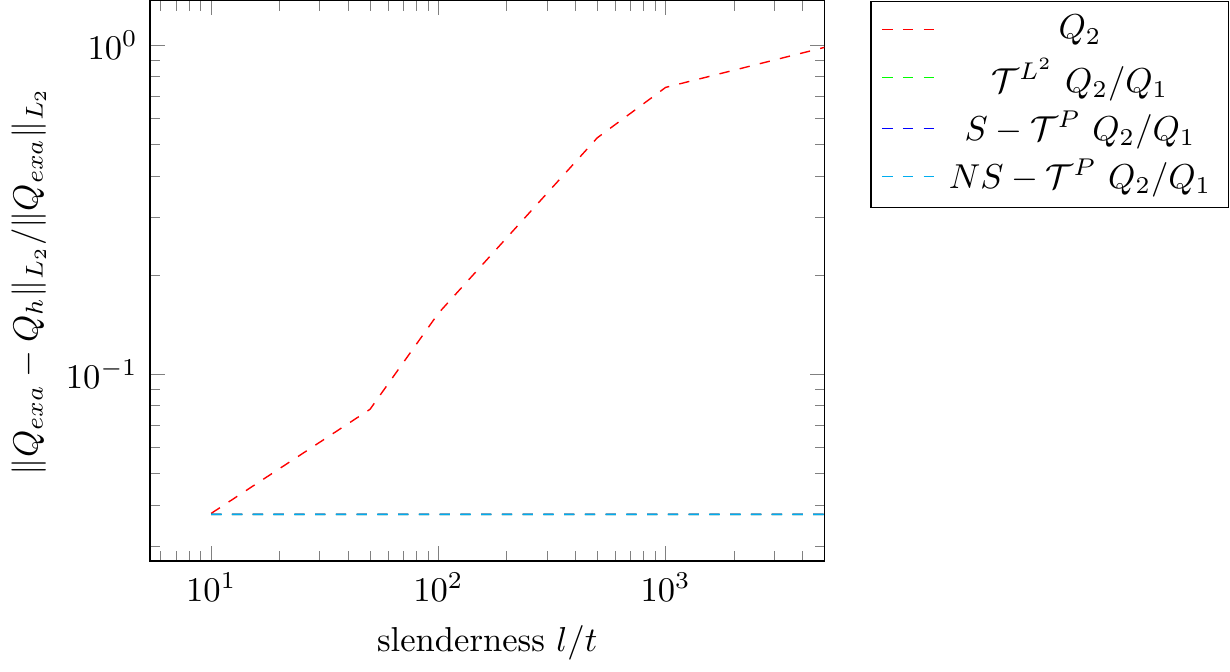}
        \caption{}
        \vspace*{2mm}
    \end{subfigure}
    \caption{Convergence study for increasing slenderness, $p=2$, and $\#{dof}=32$. Error in the $L^2$-norm for (a) displacement $w$, (b) rotation $\phi$, (c) bending moment $M$, and (d) shear force $Q$.}
    \label{fig:slenderness}
\end{figure}

We have also studied the relationship between shear locking and decreasing slenderness ratios for $p=2$. The results are shown in Figure \protect\ref{fig:slenderness}. For all three methods, the number of degrees of freedom are fixed, and the sectional and material parameters are the same as in the previous study. The slenderness ratio varies from $10$ to $5\times{10}^3$. $Q_2$ locks severely. The $\bar{B}$ methods, on the other hand, are locking free. 

\section{Volumetric locking: Nearly incompressible linear elasticity}
\label{sec:nearlyincompressible}
To demonstrate the use of \Bezier $\bar{B}$ {methods} to alleviate volumetric locking effects we study the nearly incompressible elasticity problem in two dimensions. We start with the small strain tensor $\boldsymbol{\varepsilon}$, which is defined as the symmetric part of the displacement gradient, i.e., 
\begin{equation}
\varepsilon_{ij}=\dfrac{u_{i,j}+u_{j,i}}{2}.
\end{equation}
The stress tensor is related to the strain tensor through the generalized Hooke's law
\begin{equation}
\sigma_{ij}=c_{ijkl}\varepsilon_{kl}
\end{equation}
where, for isotropic elasticity, the elastic coefficients and stress tensor can be expressed in terms of the Lam\'{e} parameters $\lambda$ and $\mu$ as
\begin{align}
c_{ijkl} &=\mu(\delta_{ik}\delta_{jl}+\delta_{il}\delta_{jk})+\lambda\delta_{ij}\delta_{kl}\\
\sigma_{ij} &= \lambda\varepsilon_{kk}\delta_{ij} + 2\mu\varepsilon_{ij}.
\end{align}
The Lam\'e parameters $\lambda$ and $\mu$ are defined in terms of Young's modulus, $E$, and Poisson's ratio, $\nu$, as
\begin{align}
\lambda &=\dfrac{\nu{E}}{(1+\nu)(1-2\nu)}\\
\mu &=\dfrac{\nu{E}}{2(1+\nu)}.
\end{align}
we can write the strong form of linear elasticity as
\begin{align}
	\sigma_{ij,j}+f_i &= 0 \text{ in $\Omega$} \\
    u_i &= g_i\text{ on $\Gamma_{g_i}$} \\
    \sigma_{ij}n_j &= h_i\text{ on $\Gamma_{h_i}$}
\end{align}
where Dirichlet boundary conditions are applied on $\Gamma_{g_i}$, Neumann boundary conditions are applied on $\Gamma_{h_i}$, and the closure of the domain $\Omega$ is $\bar{\Omega}=\Omega\cup\Gamma_{g_i}\cup\Gamma_{h_i}$. {To demonstrate the source of volumetric locking, we introduce the pressure term}
\begin{equation}
    p=-(\lambda+2\mu/3){\epsilon}_{ii}.
\end{equation}
If $\nu\rightarrow\dfrac{1}{2}$ then $\lambda$ becomes very large and {the additional constraint $\epsilon_{ii} = 0$ is applied to the volumetric strain.}

\subsection{Symmetric \Bezier $\bar{B}$ projection}
\subsubsection{The weak form}
The $\bar{B}$ approach for nearly incompressible linear elasticity splits the strain tensor $\varepsilon$ into volumetric and deviatoric strains and then replaces the volumetric strain with a projected strain. We begin with
\begin{equation}
\boldsymbol{\varepsilon}(\mathbf{u})=\boldsymbol{\varepsilon}^{vol}(\mathbf{u})+\boldsymbol{\varepsilon}^{dev}(\mathbf{u})
\end{equation}
where $\boldsymbol{\varepsilon}^{vol}=\dfrac{1}{3}\mathrm{tr}(\boldsymbol{\varepsilon})\mi$ is the volumetric strain and $\boldsymbol{\varepsilon}^{dev}=\boldsymbol{\varepsilon}-\dfrac{1}{3}\mathrm{tr}(\boldsymbol{\varepsilon})\mi$ is the deviatoric strain. The volumetric strain is then replaced by a projected volumetric strain $\bar{\boldsymbol{\varepsilon}}^{vol}$ and the new total strain becomes
\begin{equation}
\bar{\boldsymbol{\varepsilon}}=\bar{\boldsymbol{\varepsilon}}^{vol}+\boldsymbol{\varepsilon}^{dev}.
\end{equation}
The weak form can then be written as: find $\mathbf{u}\in{\mathcal{S}(\Omega)}$ such that for all $\mathbf{w}\in{\mathcal{V}(\Omega)}$
\begin{equation}
    \bar{a}\langle{\mathbf{w},\mathbf{u}}\rangle=\bar{l}\langle{\mathbf{w}}\rangle
\end{equation}
where
\begin{align}
    \bar{a}\langle{\mathbf{w},\mathbf{u}}\rangle&=\int_{\Omega}\bar{\epsilon}_{ij}(\mathbf{w})c_{ijkl}\bar{\epsilon}_{kl}(\mathbf{u})d\Omega, \\
    \bar{l}\langle{\mathbf{w}}\rangle&=\int_{\Omega}\mathbf{w}\cdot\mathbf{f}d\Omega+\int_{\Gamma_{h}}\mathbf{w}\cdot\mathbf{h}d\Gamma.
\end{align}

\subsubsection{Discretization}
Following the same approach as was described for Timoshenko beams in Section~\protect\ref{sec:timoshenkobeam} we define element level strain-displacement matrices in terms of the Bernstein basis 
\begin{align}
{\mathbf{B}^{dev}_e}&=
\begin{bmatrix}
\dfrac{2}{3}\dfrac{\partial{B}^e_{0,p}}{\partial{x}} & -\dfrac{1}{3}\dfrac{\partial{B}^e_{0,p}}{\partial{y}} & -\dfrac{1}{3}\dfrac{\partial{B}^e_{0,p}}{\partial{z}} & \cdots & \dfrac{2}{3}\dfrac{\partial{B}^e_{p,p}}{\partial{x}} & -\dfrac{1}{3}\dfrac{\partial{B}^e_{p,p}}{\partial{y}} & -\dfrac{1}{3}\dfrac{\partial{B}^e_{p,p}}{\partial{z}} \\
-\dfrac{1}{3}\dfrac{\partial{B}^e_{0,p}}{\partial{x}}  & \dfrac{2}{3}\dfrac{\partial{B}^e_{0,p}}{\partial{y}} & -\dfrac{1}{3}\dfrac{\partial{B}^e_{0,p}}{\partial{z}} & \cdots & -\dfrac{1}{3}\dfrac{\partial{B}^e_{p,p}}{\partial{x}}  & \dfrac{2}{3}\dfrac{\partial{B}^e_{p,p}}{\partial{y}} & -\dfrac{1}{3}\dfrac{\partial{B}^e_{p,p}}{\partial{z}}\\
-\dfrac{1}{3}\dfrac{\partial{B}^e_{0,p}}{\partial{x}} & -\dfrac{1}{3}\dfrac{\partial{B}^e_{0,p}}{\partial{y}} & \dfrac{2}{3}\dfrac{\partial{B}^e_{0,p}}{\partial{z}} & \cdots & -\dfrac{1}{3}\dfrac{\partial{B}^e_{p,p}}{\partial{x}} & -\dfrac{1}{3}\dfrac{\partial{B}^e_{p,p}}{\partial{y}} & \dfrac{2}{3}\dfrac{\partial{B}^e_{p,p}}{\partial{z}}\\
0 & \dfrac{\partial{B}^e_{0,p}}{\partial{z}} & \dfrac{\partial{B}^e_{0,p}}{\partial{y}} & \cdots & 0 & \dfrac{\partial{B}^e_{p,p}}{\partial{z}} & \dfrac{\partial{B}^e_{p,p}}{\partial{y}}\\
\dfrac{\partial{B}^e_{0,p}}{\partial{z}} & 0 & \dfrac{\partial{B}^e_{0,p}}{\partial{x}} & \cdots & \dfrac{\partial{B}^e_{p,p}}{\partial{z}} & 0 & \dfrac{\partial{B}^e_{p,p}}{\partial{x}}\\
\dfrac{\partial{B}^e_{0,p}}{\partial{y}} & \dfrac{\partial{B}^e_{0,p}}{\partial{x}} & 0 & \cdots & \dfrac{\partial{B}^e_{p,p}}{\partial{y}} & \dfrac{\partial{B}^e_{p,p}}{\partial{x}} & 0\\
\end{bmatrix},
\end{align}
\begin{align}
{\mathbf{B}_e^{vol}}&=
\begin{bmatrix}
\dfrac{\partial{B}^e_{0,p}}{\partial{x}} & \dfrac{\partial{B}^e_{0,p}}{\partial{y}} & \dfrac{\partial{B}^e_{0,p}}{\partial{z}} & \cdots & \dfrac{\partial{B}^e_{p,p}}{\partial{x}} & \dfrac{\partial{B}^e_{p,p}}{\partial{y}} & \dfrac{\partial{B}^e_{p,p}}{\partial{z}}
\end{bmatrix}
\end{align}
The deviatoric part of the element stiffness matrix can then be computed from the corresponding strain-displacement matrices as
\begin{align}
	\mathbf{K}_e^{dev}
    	=\mathbf{C}^e\langle {\mathbf{B}_{e}^{dev}}^T\mathbf{D}{\mathbf{B}}_{e}^{dev}\rangle (\mathbf{C}^e)^T.
\end{align}
where $\mathbf{C}^e$ is the element extraction operator for the degree $p$ spline space. The volumetric part of the stiffness matrix is computed using Bezier projection. The intermediate element matrices are
\begin{align}
\bar{\mathbf{M}}^{vol}_e&=\bar{\mathbf{C}}^e\langle{\bar{\mathbf{B}}_e^T,\dfrac{1}{3}(3\lambda+2\mu)\bar{\mathbf{B}}_e}\rangle(\bar{\mathbf{C}}^e)^T\\
\hat{\mathbf{P}}^{vol}_e&=\langle{(\hat{\bar{\mathbf{N}}}^e)^T, \mathbf{B}_e^{vol}}\rangle(\mathbf{C}^e)^T
\end{align}
where $\bar{\mathbf{C}}^e$ is the element extraction operator for the degree $p-1$ spline space, $\hat{\bar{\mathbf{N}}}^e$ are the dual basis functions restriced to the element, and \begin{equation}
\bar{\mathbf{B}}_e=
\begin{bmatrix}
{B}^e_{0,p-1} & \cdots & {B}^e_{p-1,p-1}
\end{bmatrix}.
\end{equation}The global stiffness matrix can then be assembled as
\begin{align}
	\mathbf{K} = \mathbf{K}^{dev} + \bar{\mathbf{K}}^{vol}_{S}
\end{align}
where
\begin{align}
	\mathbf{K}^{dev} &=\A_e\mathbf{K}_e^{dev}, \\
	\bar{\mathbf{K}}_S^{vol}&=\hat{\mathbf{P}}^T\bar{\mathbf{M}}\hat{\mathbf{P}} \\
	\hat{\mathbf{P}} &=\A_e\hat{{\mathbf{P}}}^{vol}_e \\
    \bar{\mathbf{M}} &=\A_e\bar{\mathbf{M}}^{vol}_e.
\end{align}
\subsection{Non-symmetric \Bezier $\bar{B}$ projection}
\subsubsection{The weak form}
{A mixed formulation for nearly incompressible elasticity can be developed by considering the pressure $p$ as an independent variable. The weak statement of the problem is then given as: find $\mathbf{u}\in{\mathcal{S}(\Omega)}$ and $p\in{L}^2(\Omega)$ such that for all $\mathbf{w}\in{\mathcal{V}(\Omega)}$ and $\delta{p}\in{L}^2(\Omega)$}
\begin{align}
    \hat{a}\langle{\mathbf{w},\mathbf{u}}\rangle-b\langle{\mathbf{w},p}\rangle&=l\langle{\mathbf{w}}\rangle, \\
    - b\langle{\mathbf{u},\delta{p}}\rangle-\dfrac{1}{(\lambda+2\mu/3)}\int_{\Omega}\delta{p}pd\Omega&=0,
\end{align}
where
\begin{align}
    \hat{a}\langle{\mathbf{w},\mathbf{u}}\rangle&=\int_{\Omega}\epsilon_{ij}(\mathbf{w})\hat{c}_{ijkl}\epsilon_{kl}(\mathbf{u})d\Omega,\\
    \hat{c}_{ijkl}&=\mu\left(\delta_{ik}\delta_{jl}+\delta_{il}\delta_{jk}-2/3\delta_{ij}\delta_{kl}\right),\\
    b\langle{\mathbf{w},p}\rangle&=\int_{\Omega}\epsilon_{ii}(\mathbf{w})pd\Omega.
\end{align}

\subsubsection{{Discretization}}

Following the same pattern as the non-symmetric \Bezier $\bar{B}$ method for the Timoshenko beam problem, we use the same discretization for $\mathbf{u}$ and $\mathbf{w}$ and use lower order spline basis functions and corresponding dual basis functions for the discretization of the pressure $p$ and its variation $\delta{p}$, respectively. The discretized stiffness matrix in mixed form can be written as
\begin{equation}
    \mathbf{K}_{mix}=
    \begin{bmatrix}
        \mathbf{K}^{dev} & -\mathbf{P}^{T}\\
        -\hat{\mathbf{P}} & -\dfrac{1}{(\lambda+2\mu/3)}\mathbf{I}
    \end{bmatrix},
\end{equation}
{where}
\begin{equation}
    \mathbf{P}=\A_e\mathbf{P}_e^{vol},
\end{equation}
with
\begin{equation}
    \mathbf{P}_e^{vol}=\bar{\mathbf{C}}^e\langle{\bar{\mathbf{B}}_e^T,\mathbf{B}^{vol}_e}\rangle(\mathbf{C}^e)^T.
\end{equation}
By eliminating the pressure control variables, we obtain a pure displacement formulation with the stiffness matrix taking the form
\begin{equation}
	\mathbf{K} = \mathbf{K}^{dev} + \bar{\mathbf{K}}^{vol}_{NS},
\end{equation}
where
\begin{equation}
    \bar{\mathbf{K}}^{vol}_{NS}=(\lambda+2\mu/3){\mathbf{P}}^T\hat{\mathbf{P}}.
\end{equation}
We note again, that in contrast to the symmetric \Bezier $\bar{B}$ method, the assembly of the stiffness matrix in this case can be performed in the standard way with no need for global matrix operations.
 
\subsection{Numerical results}
{We begin this section by numerically evaluating the inf-sup constant for the global $\bar{B}$ and non-symmetric \Bezier $\bar{B}$ methods.} 
We then investigate the performance of the B\'{e}zier $\bar{B}$ method for two nearly incompressible linear elasticity problems under plane strain conditions. {For the numerical examples,} we first study the Cook's membrane problem, which is discretized with B-spline basis functions, and in the second problem we model an infinite plate with a circular hole using NURBS. Results computed using standard finite elements are labeled $Q_1$, $Q_2$, $Q_3$, and $Q_4$. Results computed using a global $\bar{B}$ method are labeled $\mathcal{T}^{L^2}$ and those computed with the {symmetric} B\'ezier $\bar{B}$ method and the {non-symmetric B\'ezier $\bar{B}$ method} are labeled $S-\mathcal{T}^{P}$ and $NS-\mathcal{T}^{P}$, respectively.

\subsubsection{{Numerical evaluation of the inf-sup condition}}

The inf-sup condition is also refered to as the Ladyzhenskaya-Babuska-Brezzi condition (or simply LBB)~\protect\cite{babuvska1973finite, brezzi1974existence, ladyzhenskaya1969mathematical}. It is a crucial condition to ensure the solvability, stability and optimality of a mixed problem. For the nearly incompressible elasticity problem the inf-sup condition is stated as: for $\delta{p}\neq{0}$ and $\mathbf{u}\neq{0}$
\begin{equation}
    \inf_{\delta{p}\in{L^2(\Omega)}}\sup_{\mathbf{u}\in\mathcal{S}(\Omega)}\dfrac{\vert{b\langle{\mathbf{u},\delta{p}}\rangle}\vert}{\|{\delta{p}}\|_{L^2(\Omega)}\|{\mathbf{u}}\|_{H^1(\Omega)}}\geq\beta>0.
\end{equation}
In a discretized problem, the inf-sup condition requires the variable $\beta$ to be a constant that is independent of the mesh size.\par

\begin{figure}[htb!]
	\centering
	\includegraphics[width=0.5\linewidth]{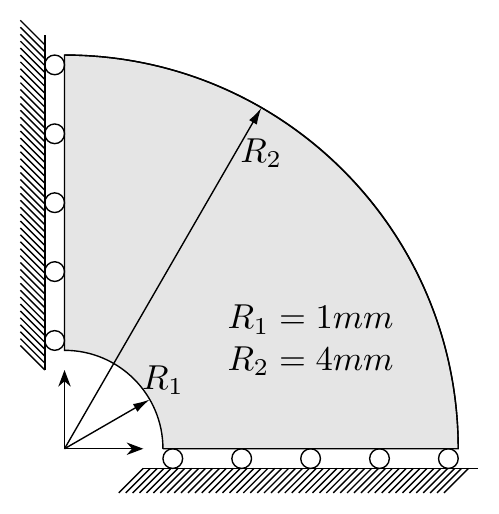}
	\caption{Geometry and boundary conditions for the inf-sup test.}
	\label{fig:quarter_annulus}
\end{figure}

\begin{figure}
    \centering
    \begin{subfigure}[b]{0.18\linewidth}        
        \centering
        \includegraphics[width=\linewidth]{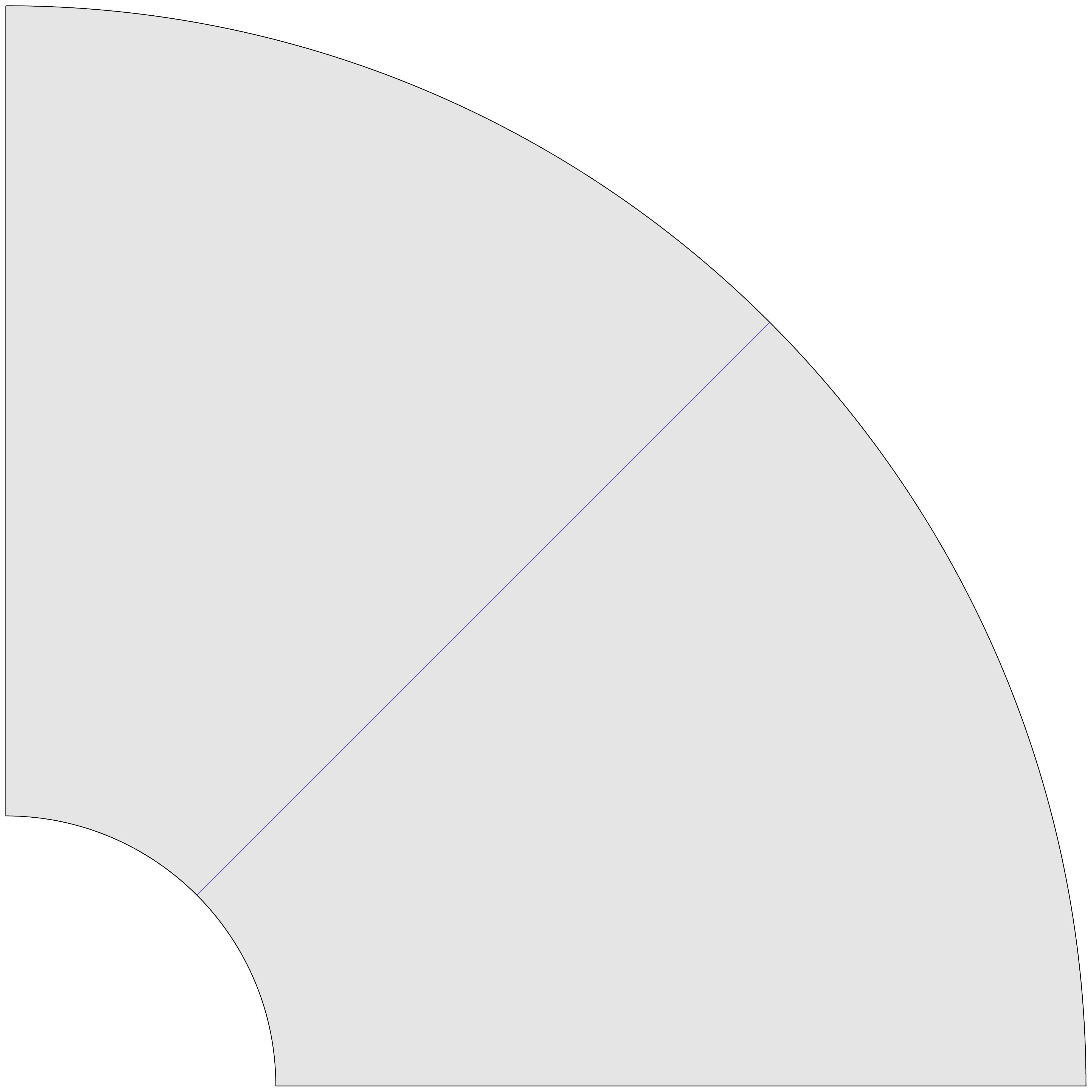}
    \end{subfigure}
    \begin{subfigure}[b]{0.18\linewidth}        
        \centering
        \includegraphics[width=\linewidth]{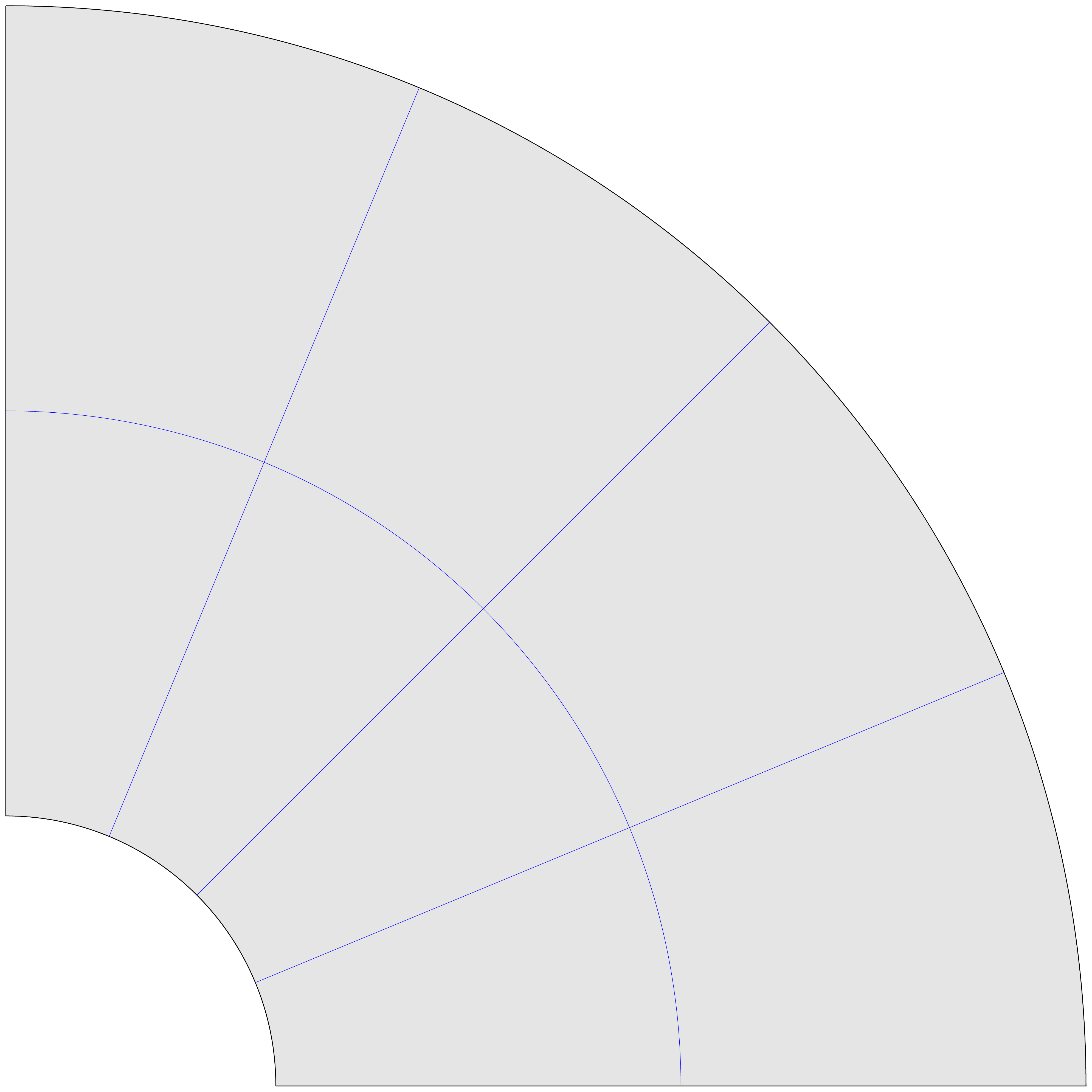}
    \end{subfigure}
    \begin{subfigure}[b]{0.18\linewidth}        
        \centering
        \includegraphics[width=\linewidth]{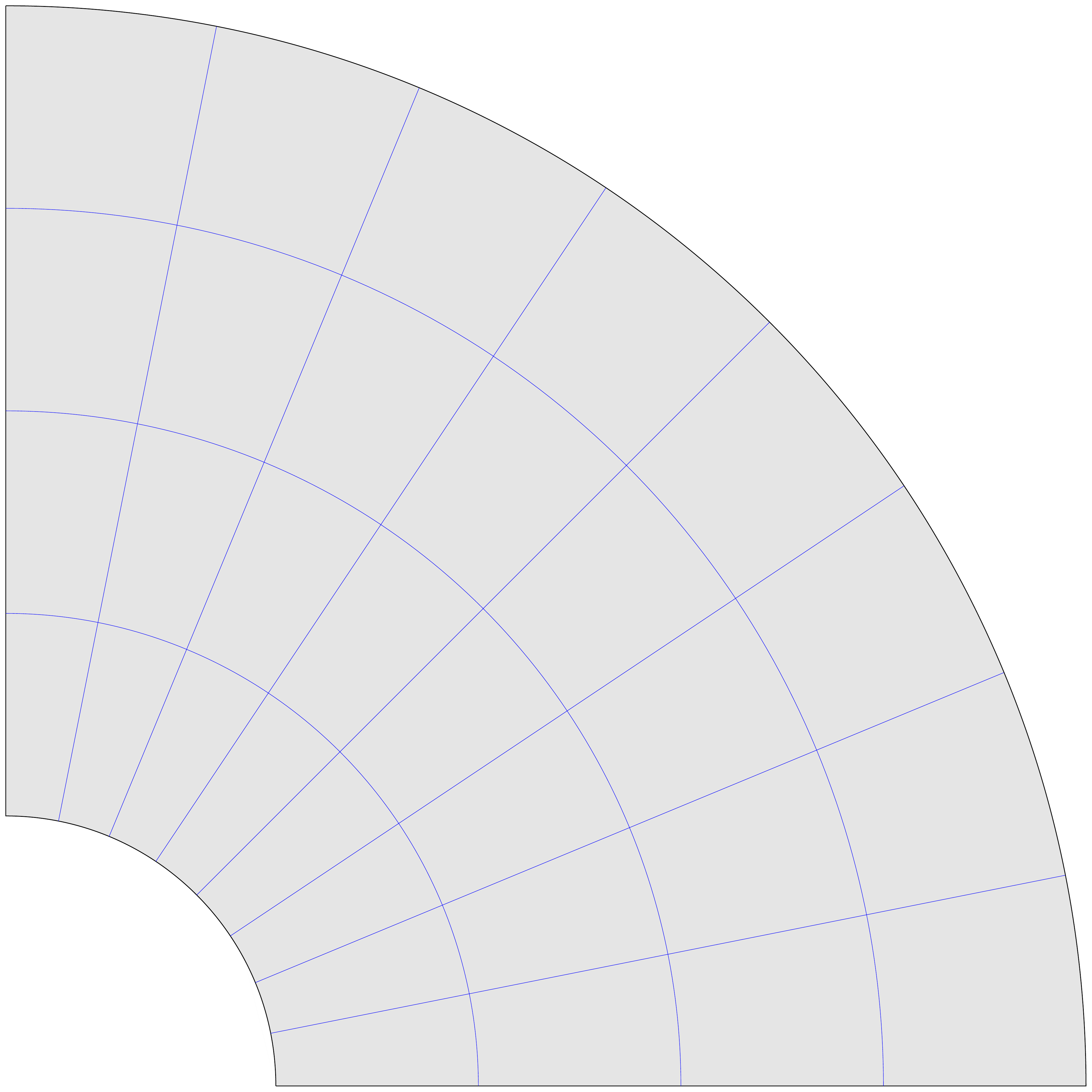}
    \end{subfigure}
    \begin{subfigure}[b]{0.18\linewidth}        
        \centering
        \includegraphics[width=\linewidth]{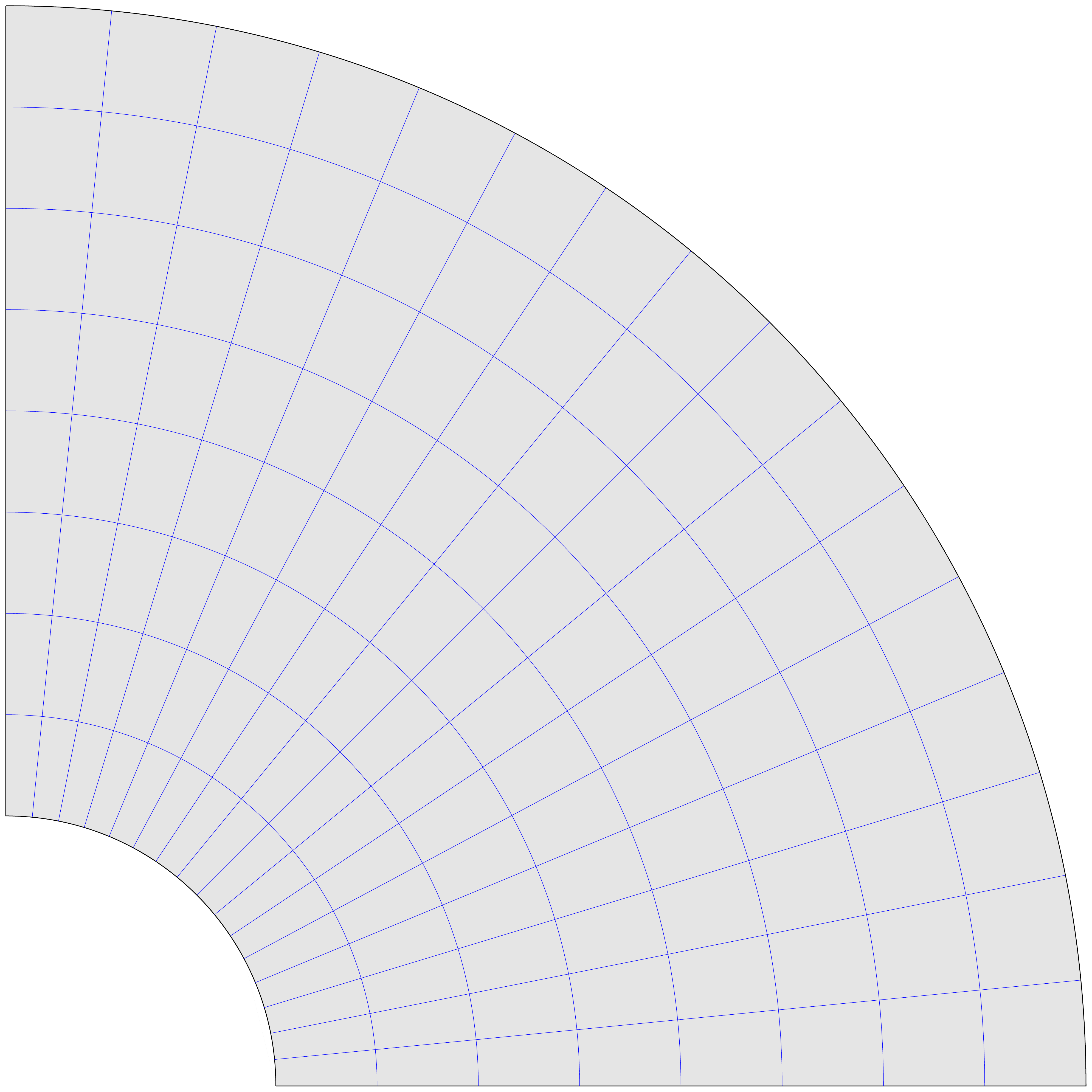}
    \end{subfigure}
    \begin{subfigure}[b]{0.18\linewidth}        
        \centering
        \includegraphics[width=\linewidth]{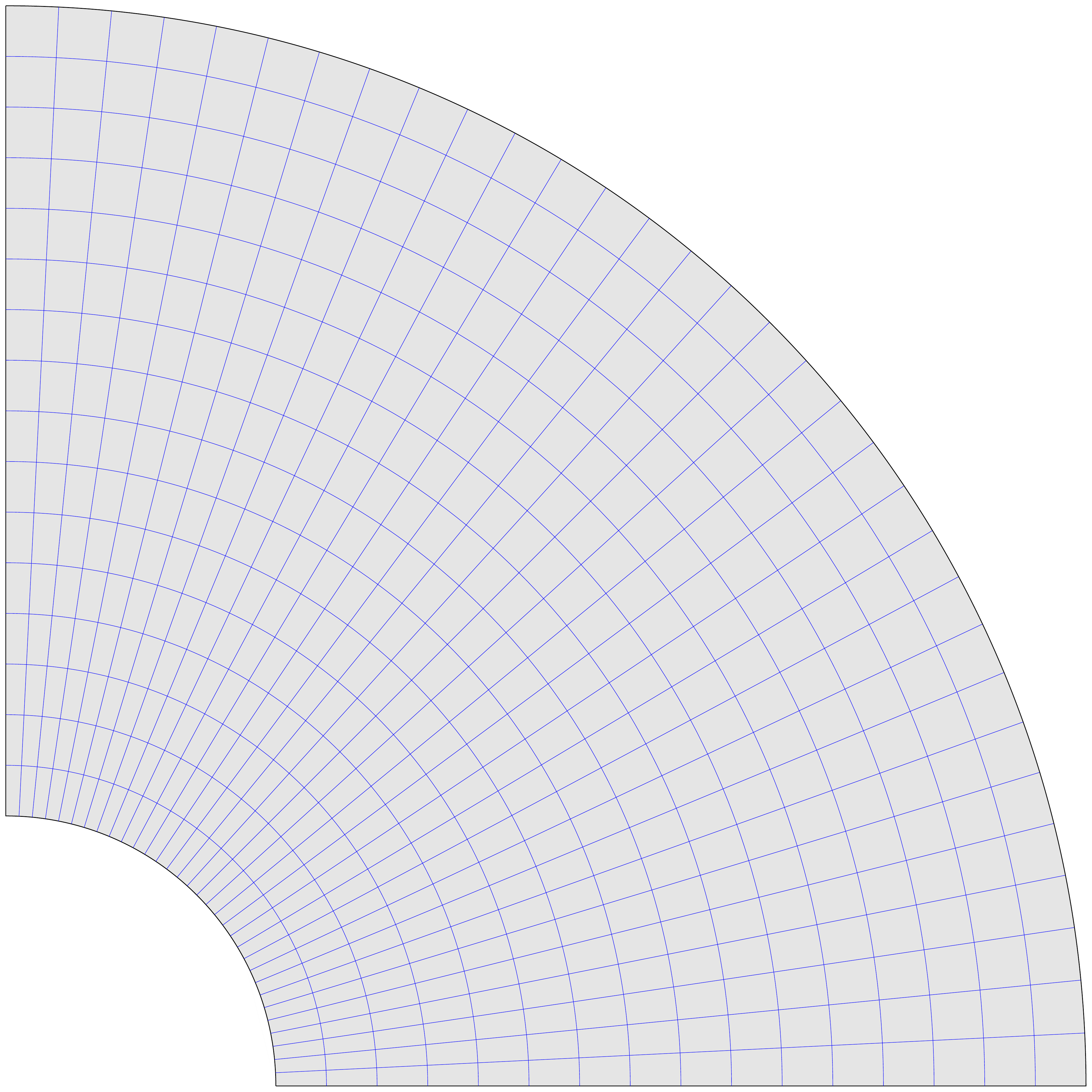}
    \end{subfigure}
    \caption{Sequence of meshes for inf-sup test.}
    \label{fig:mesh_hole}
\end{figure}

Here, we consider the inf-sup condition of a uniformly refined quarter annulus. The geometry and boundary conditions are shown in Figure \protect\ref{fig:quarter_annulus}. The geometry of the quarter annulus can be exactly represented using a biquadratic NURBS basis. The knot vector for the coarsest discretization is given by
\begin{align}
\begin{split}
\Xi_\xi\times\Xi_\eta=\lbrace{0,0,0,1,1,1}\rbrace\times\lbrace{0,0,0,1,1,1}\rbrace
\end{split}
\end{align}
and the corresponding weights and control points associated with each basis function are given in Table \protect\ref{table:weights} and \protect\ref{table:control_points}. For higher-order elements and finer discretizations the weights and corresponding control points are identified by an order elevation and knot insertion algorithm, respectively. The B\'ezier mesh representation for the discretizations are shown in Figure \protect\ref{fig:mesh_hole}.

\begin{table}
\begin{center}
\caption{Weights for the plate with a circular hole}\label{table:weights}
\begin{tabular}{l@{\hskip 1cm}l@{\hskip 1cm}l@{\hskip 1cm}l}
\hline
i    & $w_{i,1} $ & $w_{i,2}$ & $w_{i,3}$\\
\hline
1    & 1    & ${1}/{\sqrt{2}}$ & 1 \\
2    & 1    & ${1}/{\sqrt{2}}$ & 1 \\
3    & 1    & ${1}/{\sqrt{2}}$ & 1 \\
\hline
\end{tabular}
\end{center}
\end{table}

\begin{table}
\begin{center}
\caption{Control points for the plate with a circular hole}\label{table:control_points}
\begin{tabular}{l@{\hskip 1cm}l@{\hskip 1cm}l@{\hskip 1cm}l}
\hline
i    & $B_{i,1} $ & $B_{i,2}$ & $B_{i,3}$\\
\hline
1    & (0,1)    & (1,1)     & (1,0)   \\
2    & (0,2.5)  & (2.5,2.5) & (2.5,0) \\
3    & (0,4)    & (4,4)     & (4,0)   \\
\hline
\end{tabular}
\end{center}
\end{table}

Only the global $\bar{B}$ method and the non-symmetric \Bezier $\bar{B}$ method are considered here, as the symmetric \Bezier $\bar{B}$ method lacks a connection to a mixed formulation. As a counter example, the well-known pair $Q_p/Q_p$ of the global $\bar{B}$ method that violates the inf-sup condition is also tested here. Our tests follow the procedure proposed by Chapelle and Bathe in \protect\cite{chapelle_inf-sup_1993}.\par
\begin{figure}[htb!]
    \centering
    \begin{subfigure}[b]{0.31\linewidth}        
        \centering
        \includegraphics[width=\linewidth]{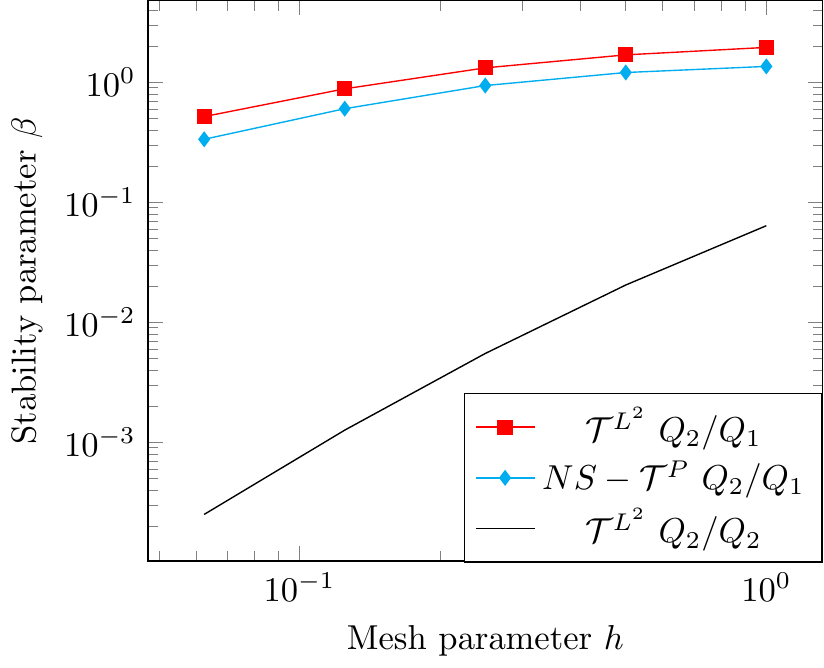}
    \end{subfigure}
    \begin{subfigure}[b]{0.31\linewidth}        
        \centering
        \includegraphics[width=\linewidth]{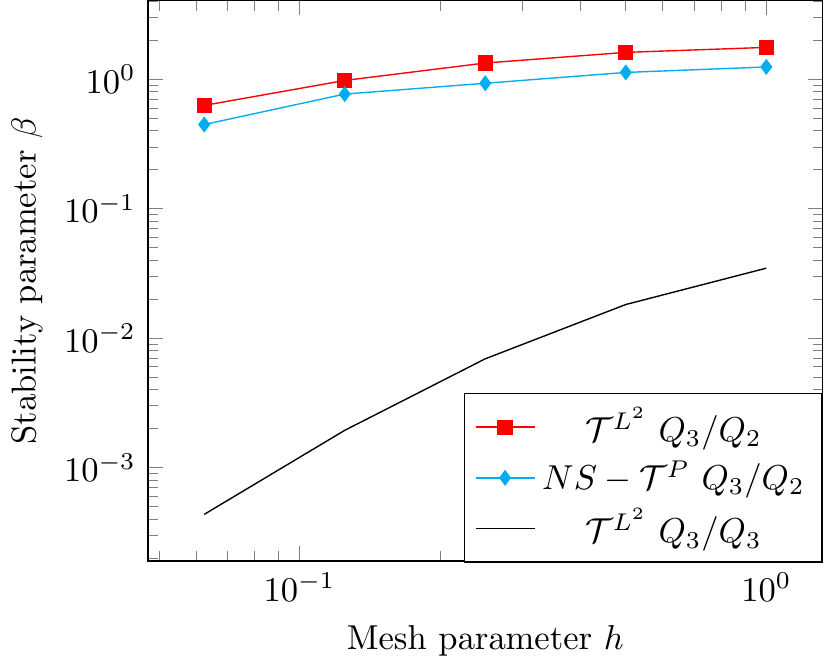}
    \end{subfigure}
    \begin{subfigure}[b]{0.31\linewidth}        
        \centering
        \includegraphics[width=\linewidth]{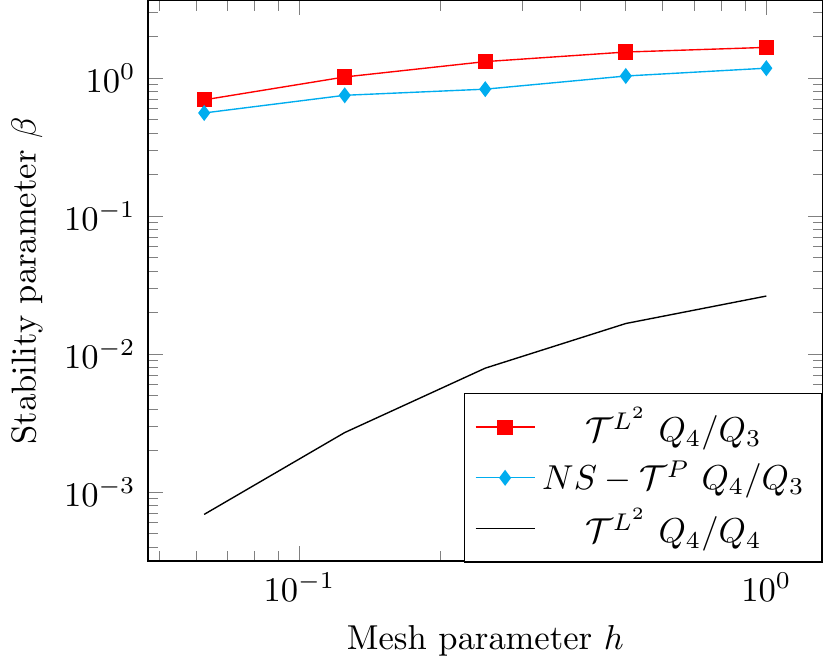}
    \end{subfigure}

    \caption{Inf-sup test results for nearly incompressible elasticity. The global $\bar{B}$ method, $\mathcal{T}^{L^2} Q4/Q3$, and the non-symmetric \Bezier $\bar{B}$ method, $NS - \mathcal{T}^P Q4/Q3$, do not strictly satisfy the LBB condition, but compared to the $\mathcal{T}^{L^2} Q_4/Q_4$ method, both methods reduce constraints to a favorable level. }
    \label{fig:inf_sup}
\end{figure}

Figure \protect\ref{fig:inf_sup} shows the numerical results. As can be seen, the global $\bar{B}$ method and the non-symmetric \Bezier $\bar{B}$ method do not strictly satisfy the LBB condition since $\beta$ is not independent of the mesh size. This result is consistent with the statement made in \protect\cite{elguedj:hal-00457010} that the global $\bar{B}$ method does not reduce the constraints sufficiently to satisfy the LBB condition. However, compared to the $Q_p/Q_p$ pair, both methods reduce constraints to a more favorable level. If we compare the results for the global $\bar{B}$ method with the non-symmetric \Bezier $\bar{B}$ method we see that their stability parameter $\beta$ decreases at the same rate and the stability parameter for the non-symmetric \Bezier $\bar{B}$ method is slightly lower than that for the global $\bar{B}$ method. These results indicate a similar optimality in convergence for both methods and a slightly higher error for the non-symmetric \Bezier $\bar{B}$ method.

\subsubsection{Cook's membrane problem}

\begin{figure}[htb!]
	\centering
	\includegraphics[width=0.5\linewidth]{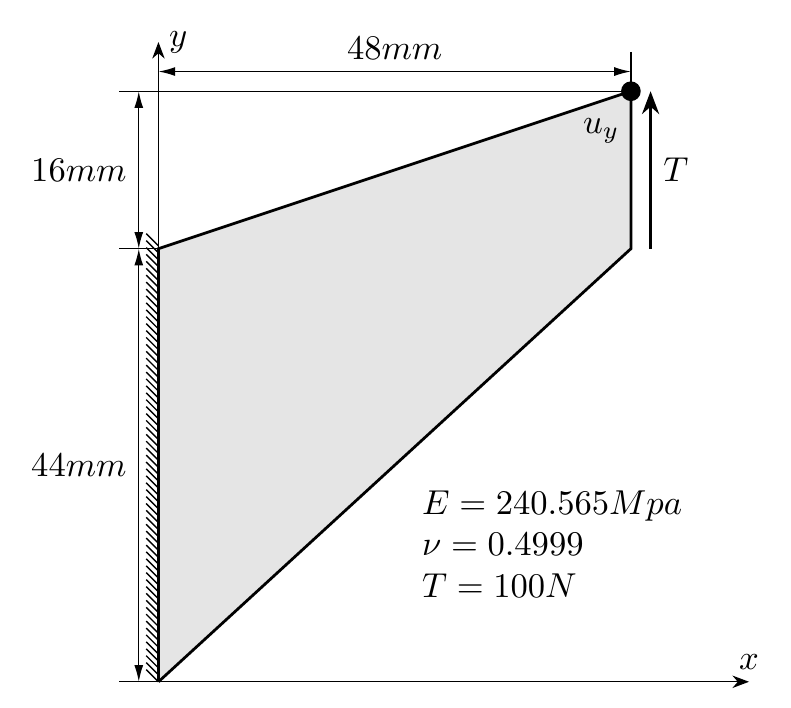}
	\caption{Geometry, boundary conditions, and material properties for the Cook's membrane problem.}
	\label{fig:Cook's}
\end{figure}

\begin{figure}[htb!]
    \centering
    \begin{subfigure}[b]{0.18\linewidth}        
        \centering
        \includegraphics[width=\linewidth]{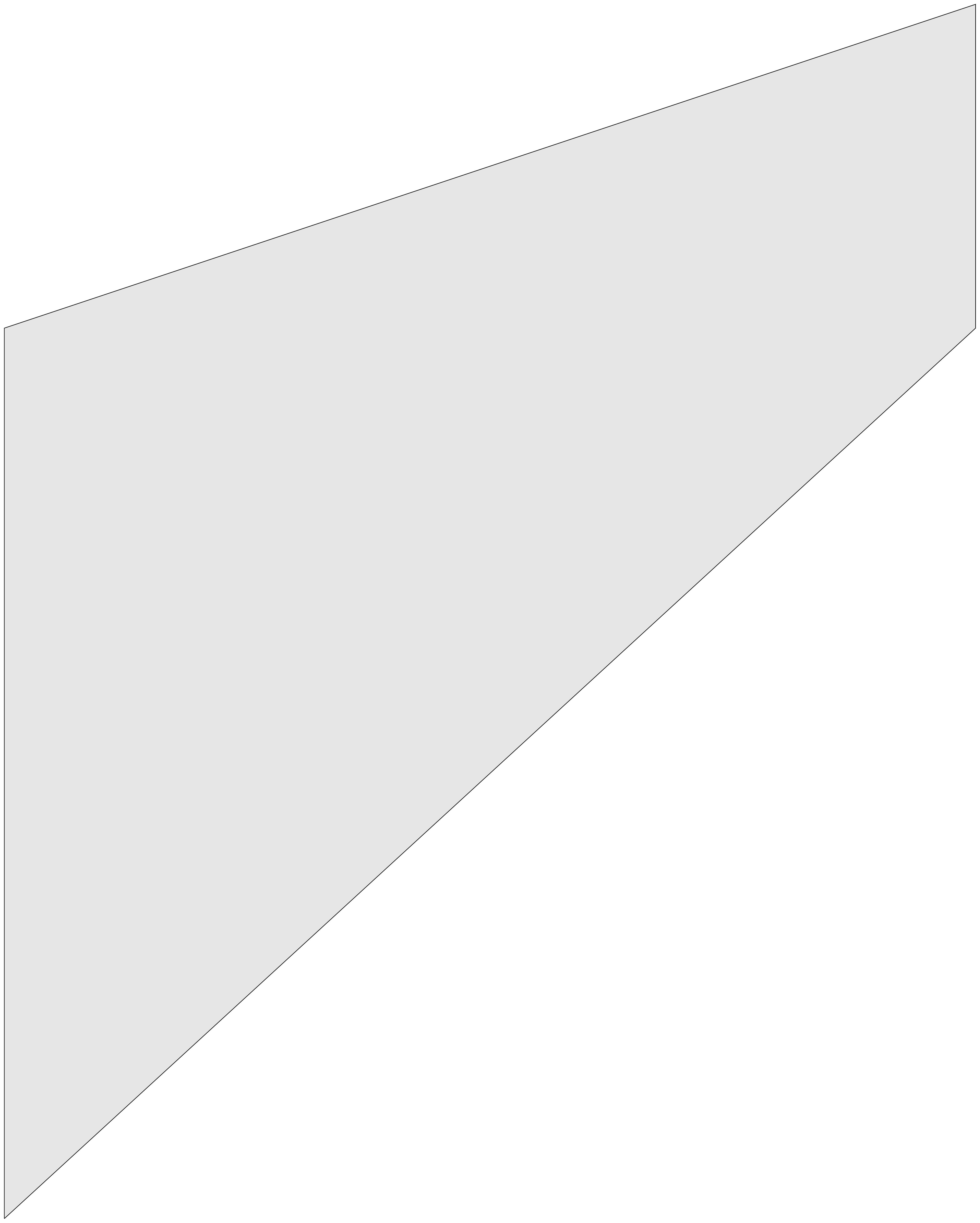}
    \end{subfigure}
    \begin{subfigure}[b]{0.18\linewidth}        
        \centering
        \includegraphics[width=\linewidth]{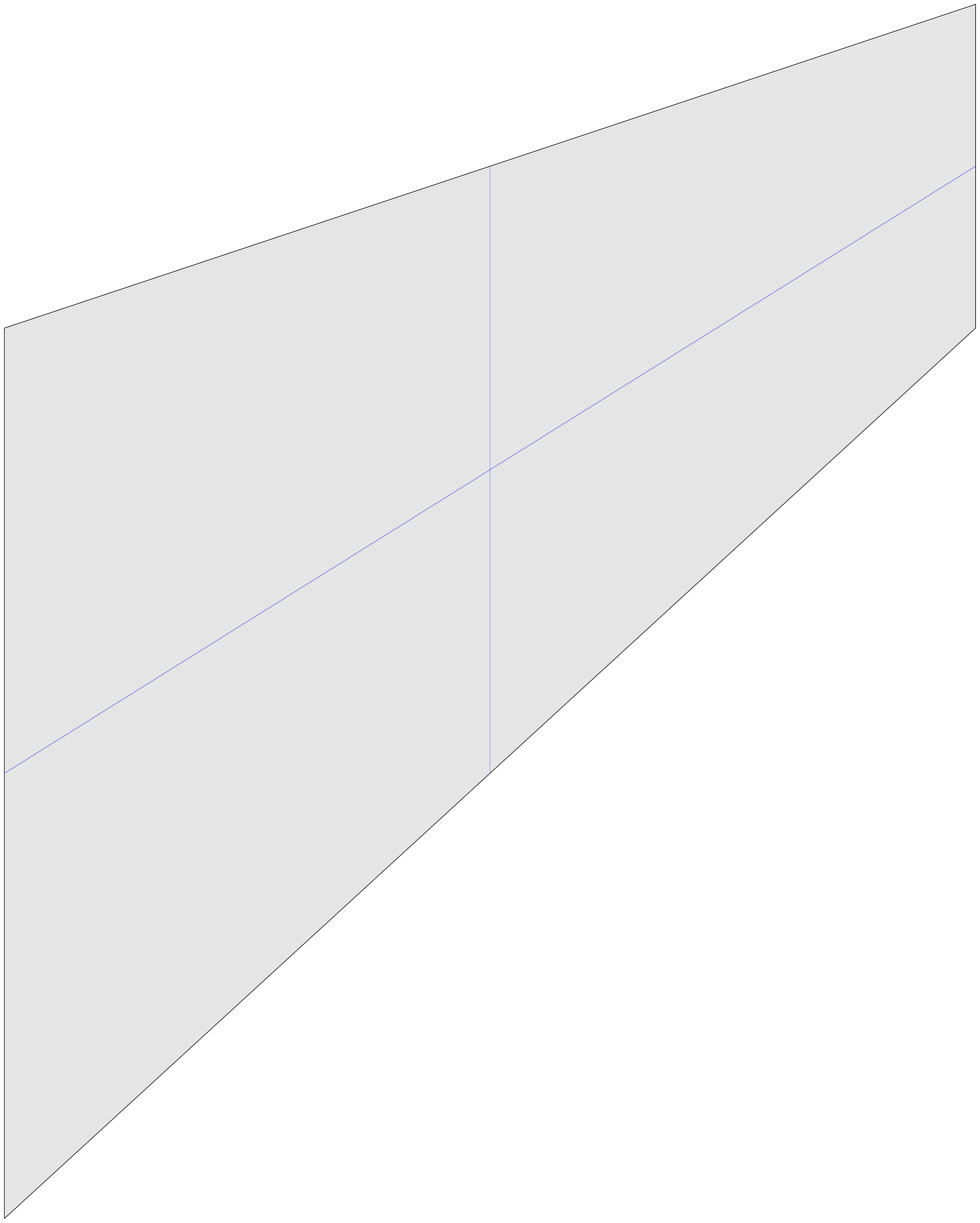}
    \end{subfigure}
    \begin{subfigure}[b]{0.18\linewidth}        
        \centering
        \includegraphics[width=\linewidth]{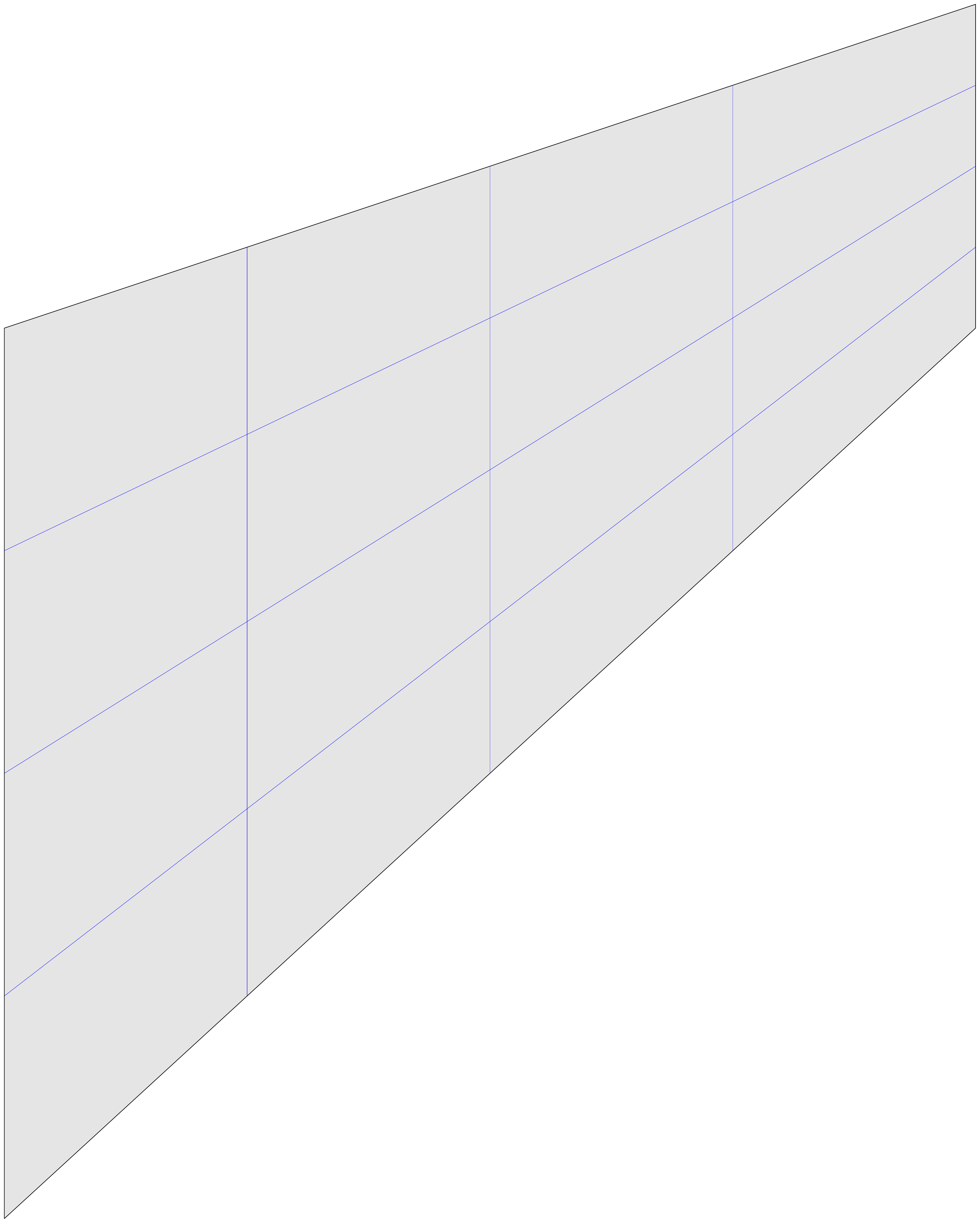}
    \end{subfigure}
    \begin{subfigure}[b]{0.18\linewidth}        
        \centering
        \includegraphics[width=\linewidth]{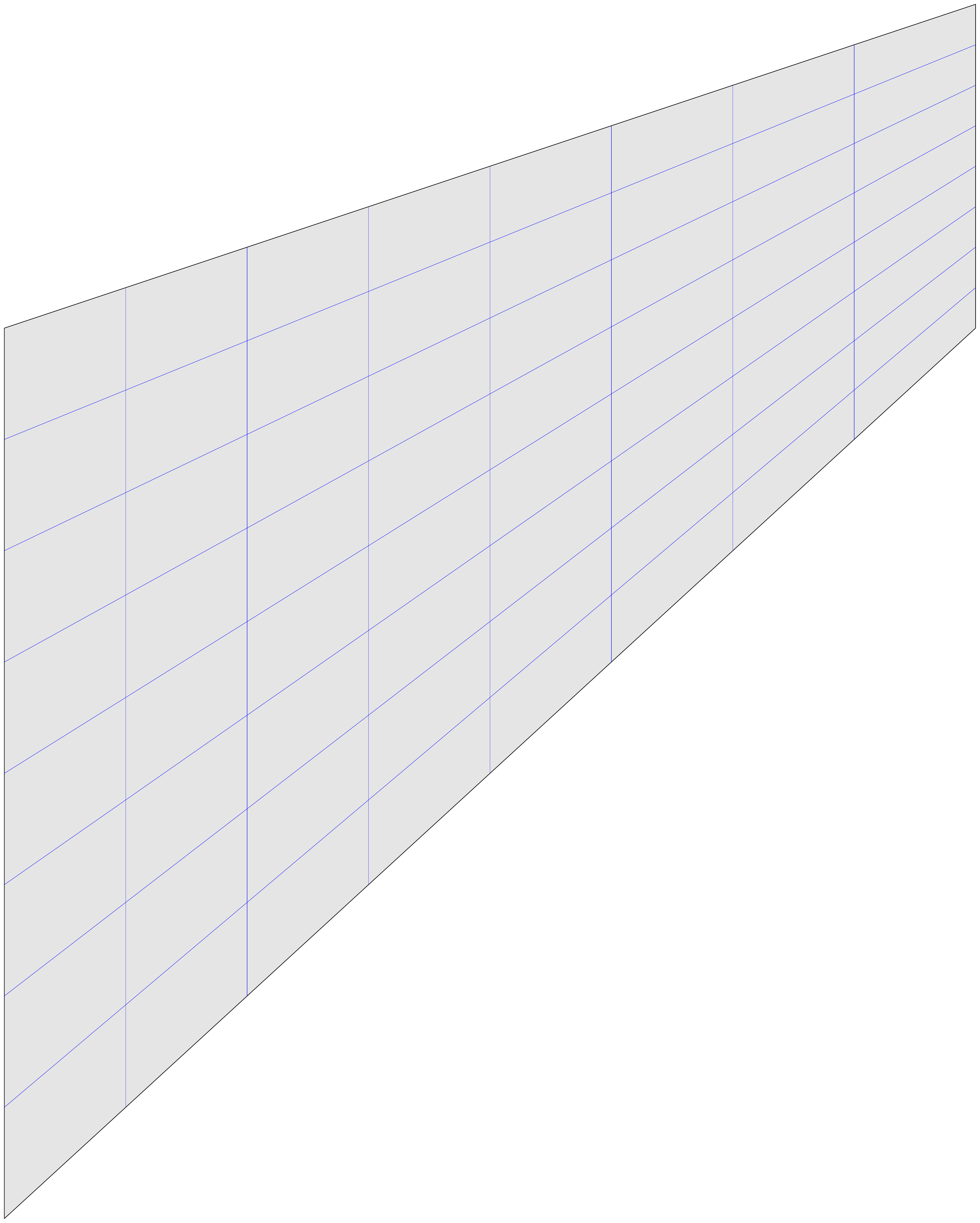}
    \end{subfigure}
    \begin{subfigure}[b]{0.18\linewidth}        
        \centering
        \includegraphics[width=\linewidth]{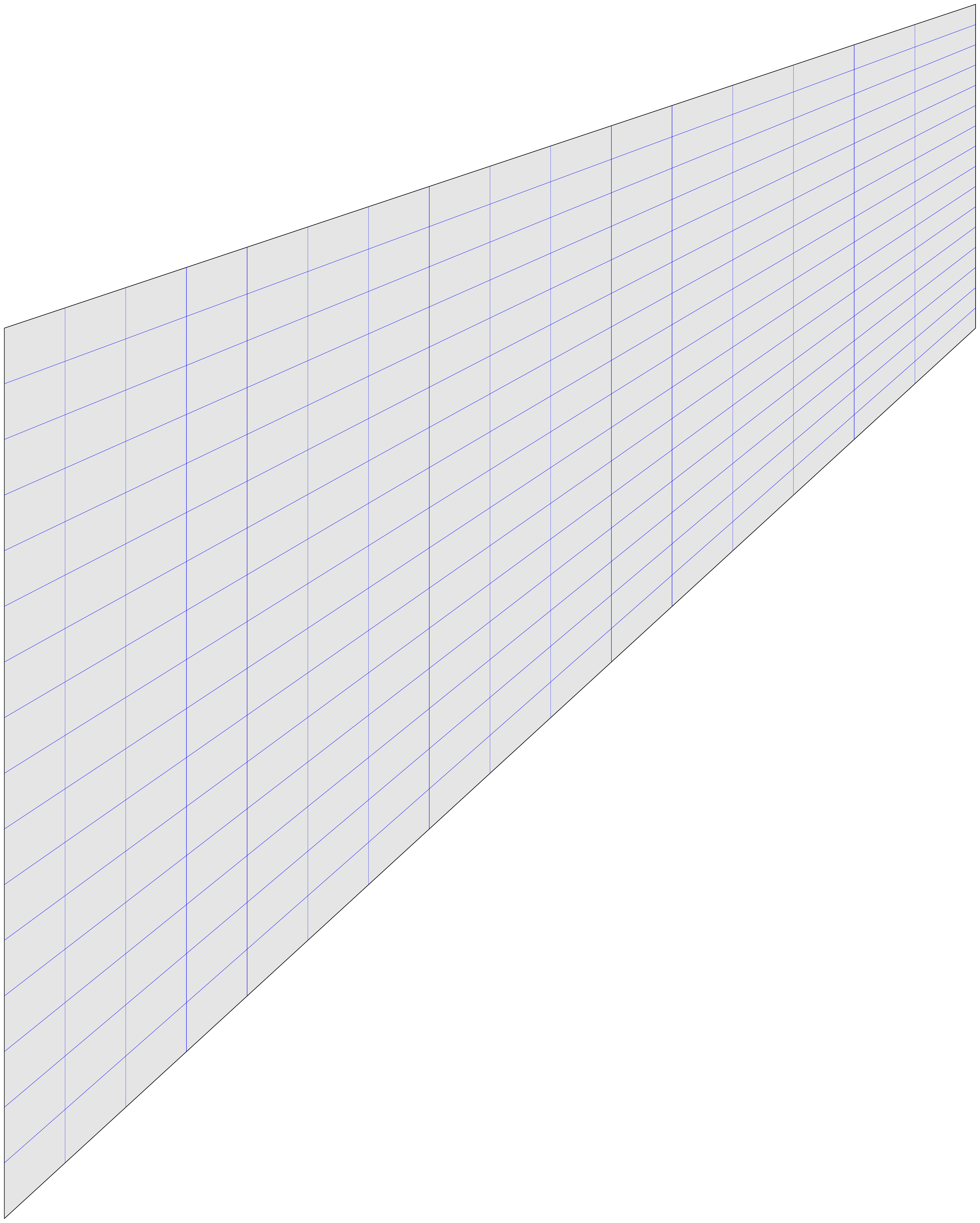}
    \end{subfigure}
    \caption{Sequence of meshes for Cook's membrane problem.}
    \label{fig:mesh_cook}
\end{figure}

This benchmark problem is a standard test for combined bending and shearing response. The geometry, boundary conditions, and material properties are shown in Figure \protect\ref{fig:Cook's}. The left boundary of the tapered panel is clamped, the top and bottom edges are free with zero traction boundary conditions, and the right boundary is subjected to a uniformly distributed traction load in the $y$-direction as shown. The meshes used are shown in Figure \protect\ref{fig:mesh_cook}.

A comparison of the displacement of the top right corner with respect to the number of elements per side is shown in Figure \protect\ref{fig:Cook's_result}. $Q_1$ locks and mesh refinement has little impact. Locking is somewhat reduced for the higher-order elements $Q_p$, $p > 1$. The $\bar{B}$ methods perform very well for all degrees.

\begin{figure}[htb!]
    \centering
	\begin{subfigure}[b]{\textwidth}
        \includegraphics[width=.8\linewidth]{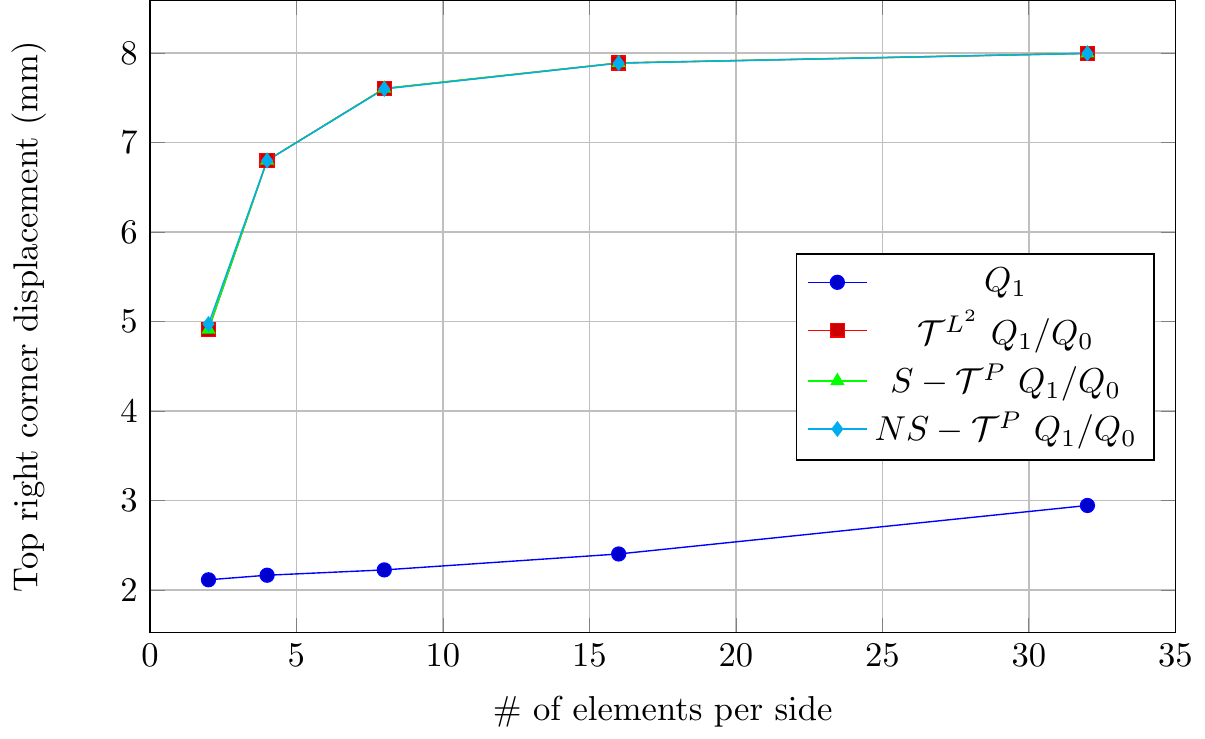}
    \end{subfigure}
    \centering
	\begin{subfigure}[b]{\textwidth}
        \includegraphics[width=.8\linewidth]{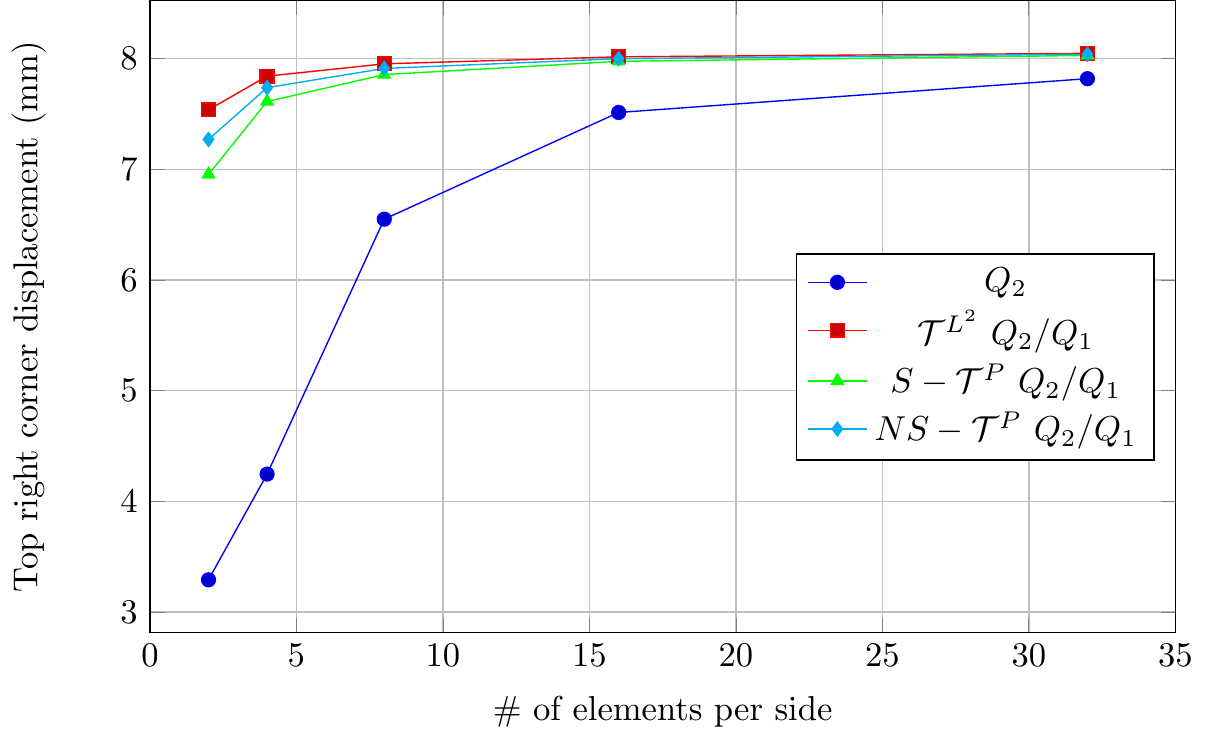}
    \end{subfigure}
\end{figure}

\begin{figure}[htb!]\ContinuedFloat
    \centering
    \begin{subfigure}[b]{\textwidth}
        \includegraphics[width=.8\linewidth]{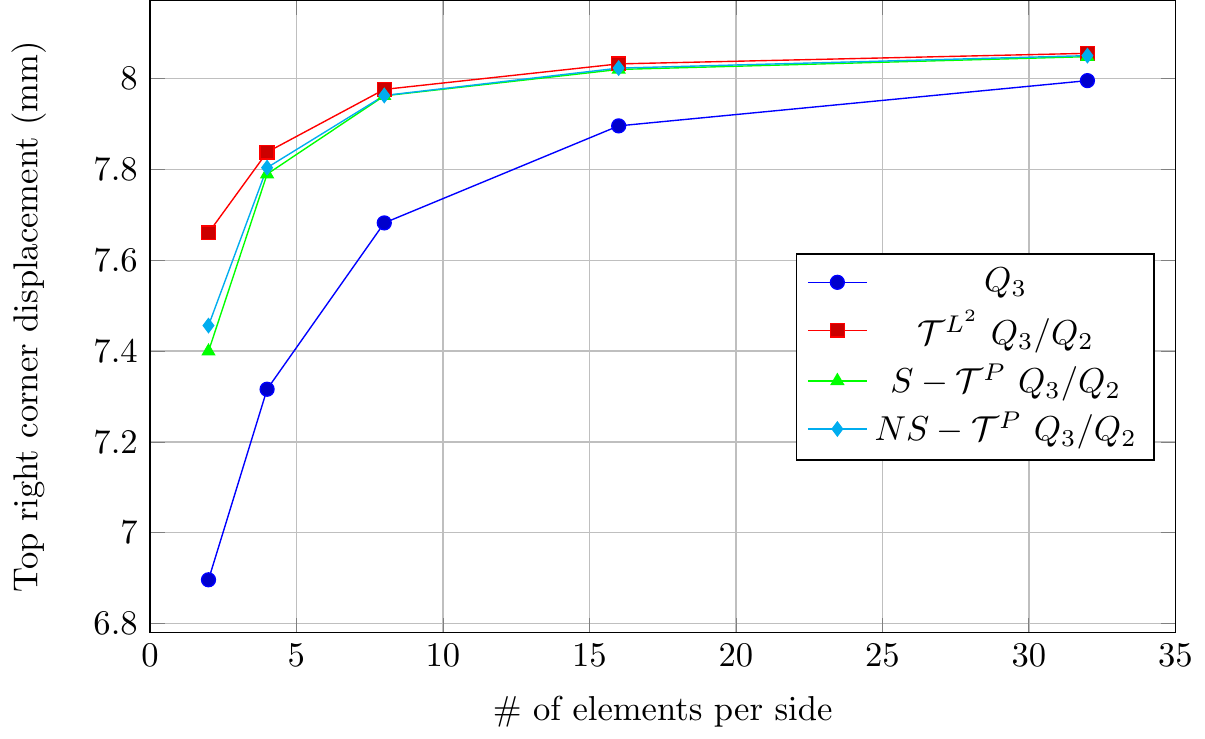}
    \end{subfigure}
    \centering
     \begin{subfigure}[b]{\textwidth}
        \includegraphics[width=.8\linewidth]{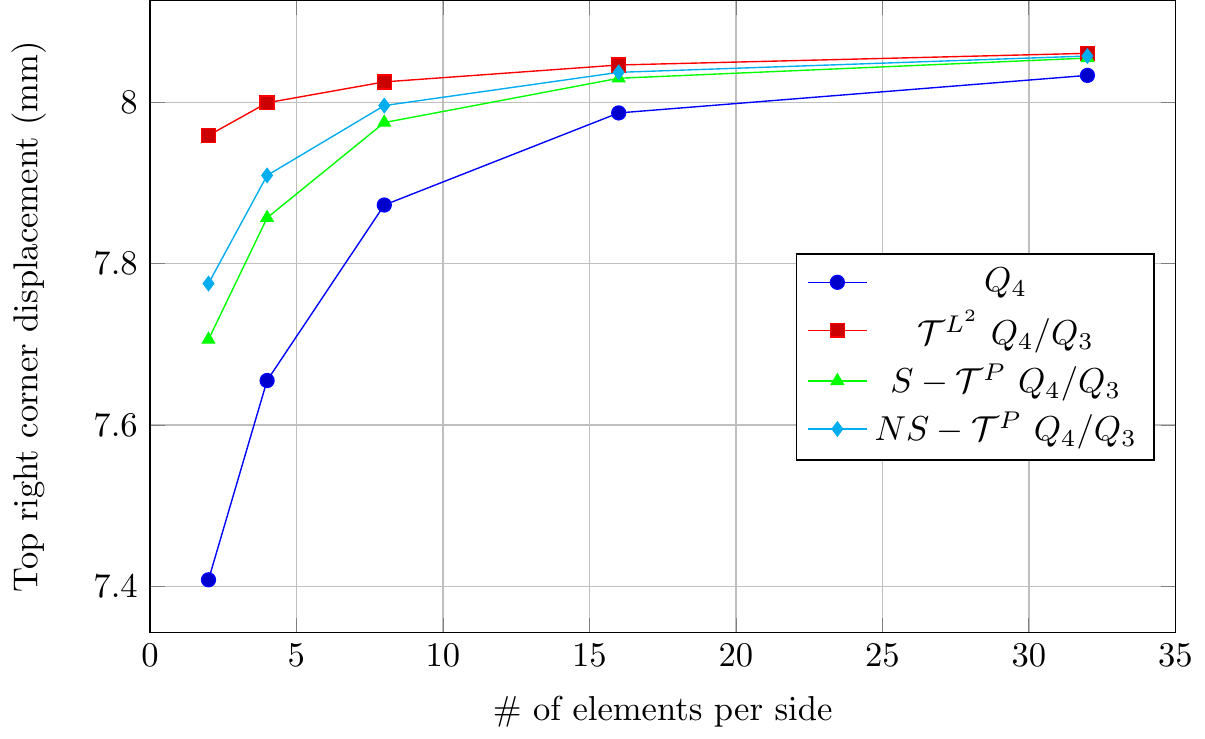}
    \end{subfigure}
	\caption{Cook's membrane: comparison of the vertical displacement at the top right corner for the different methods and degrees.}
	\label{fig:Cook's_result}
\end{figure}

\clearpage

\subsubsection{Infinite plate with a circular hole}

\begin{figure}[htb!]
	\centering
	\begin{subfigure}[t]{0.5\textwidth}
    \centering
    \includegraphics[width=1\linewidth]{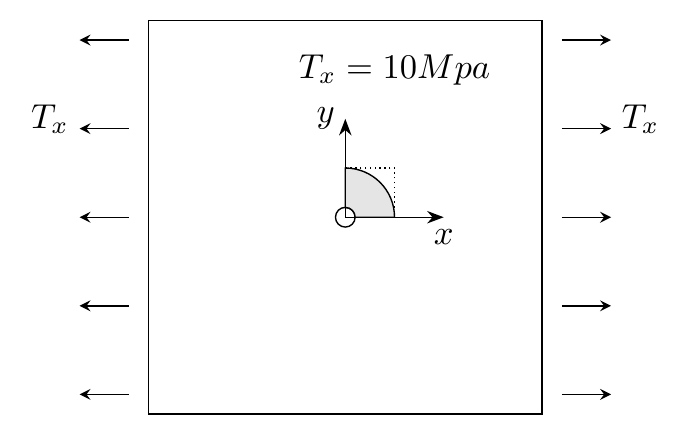}
    \caption{Infinite plate with a hole subjected to uniaxial tension at $x=\pm\infty$.}
    \label{fig:platewithhole_geometry_a}
    \end{subfigure}%
    ~ 
    \begin{subfigure}[t]{0.5\textwidth}
    \centering
    \includegraphics[width=1\linewidth]{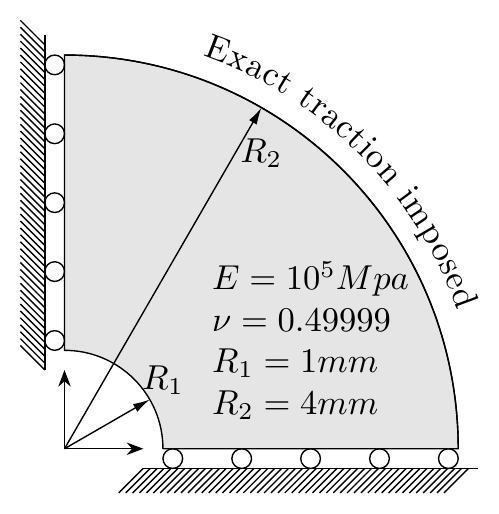}
    \caption{A representation of the computational model.}
    \label{fig:platewithhole_geometry_b}
    \end{subfigure}
    \caption{Geometry, boundary conditions, and material properties for the infinite plate with a hole.}
	\label{fig:platewithhole_geometry}
\end{figure}

The setup for the infinite plate with a circular hole problem is shown in Figure \protect\ref{fig:platewithhole_geometry} and the discretizations are shown in Figure \protect\ref{fig:mesh_hole}. The traction along the outer edge is evaluated from the exact solution which is given by
\begin{align}
\begin{split}
\sigma_{rr}(r,\theta)&=\dfrac{T_x}{2}(1-\dfrac{R_1^2}{r^2})+\dfrac{T_x}{2}(1-4\dfrac{R_1^2}{r^2}+3\dfrac{R_1^4}{r^4})cos(2\theta)\\
\sigma_{\theta\theta}(r,\theta)&=\dfrac{T_x}{2}(1+\dfrac{R_1^2}{r^2})-\dfrac{T_x}{2}(1+3\dfrac{R_1^4}{r^4})cos(2\theta)\\
\sigma_{r\theta}(r,\theta)&=-\dfrac{T_x}{2}(1+2\dfrac{R_1^2}{r^2}-3\dfrac{R_1^4}{r^4})sin(2\theta).
\end{split}
\end{align}

\begin{figure}[h]
    \center
\begin{tabular}{ccc}
\includegraphics[width=.31\linewidth]{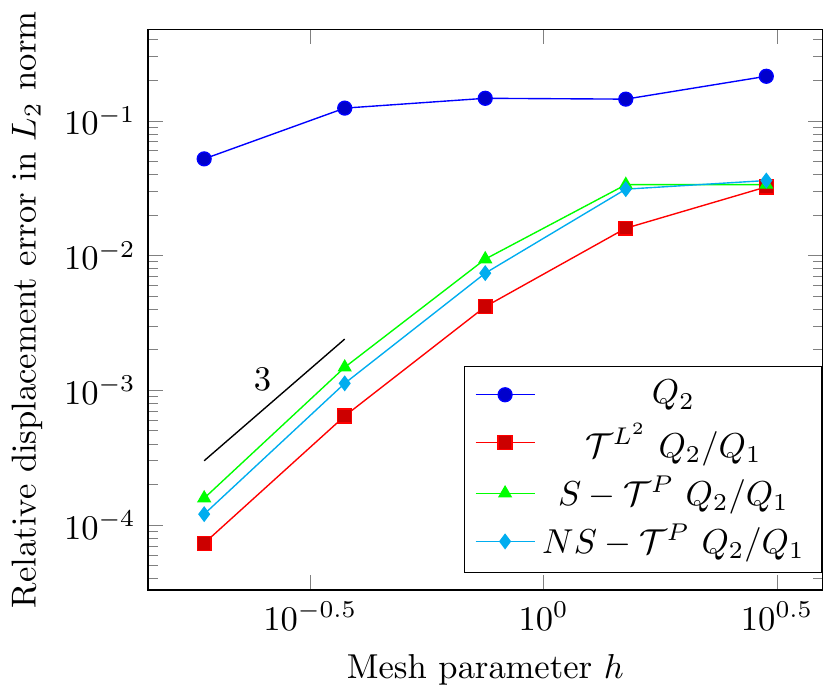} & \includegraphics[width=.31\linewidth]{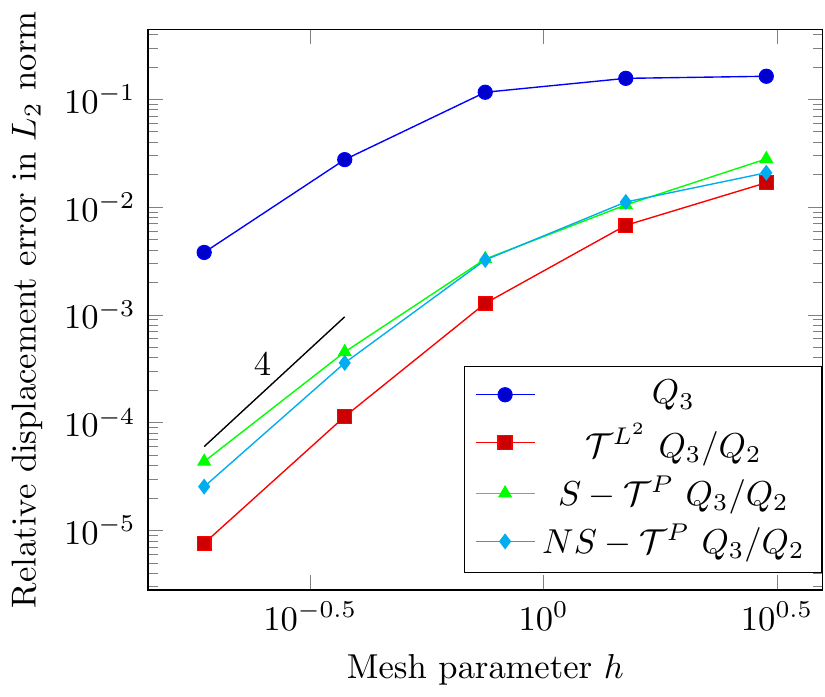} & \includegraphics[width=.31\linewidth]{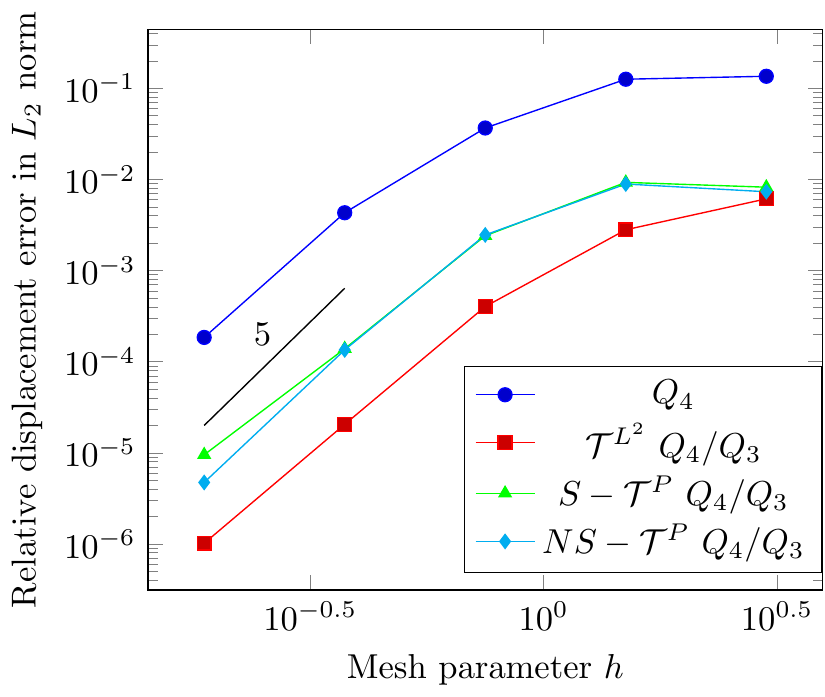}\\
\includegraphics[width=.31\linewidth]{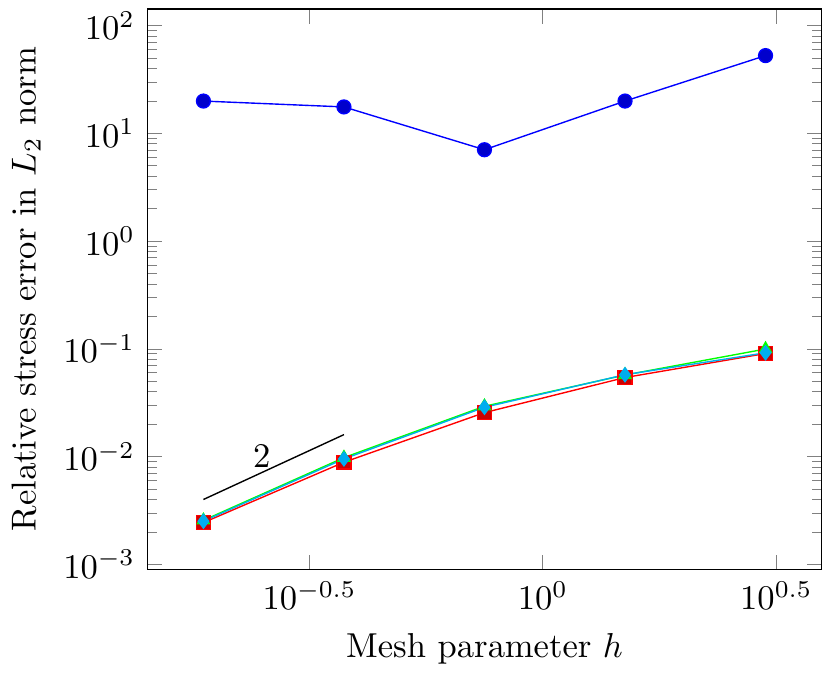} & \includegraphics[width=.31\linewidth]{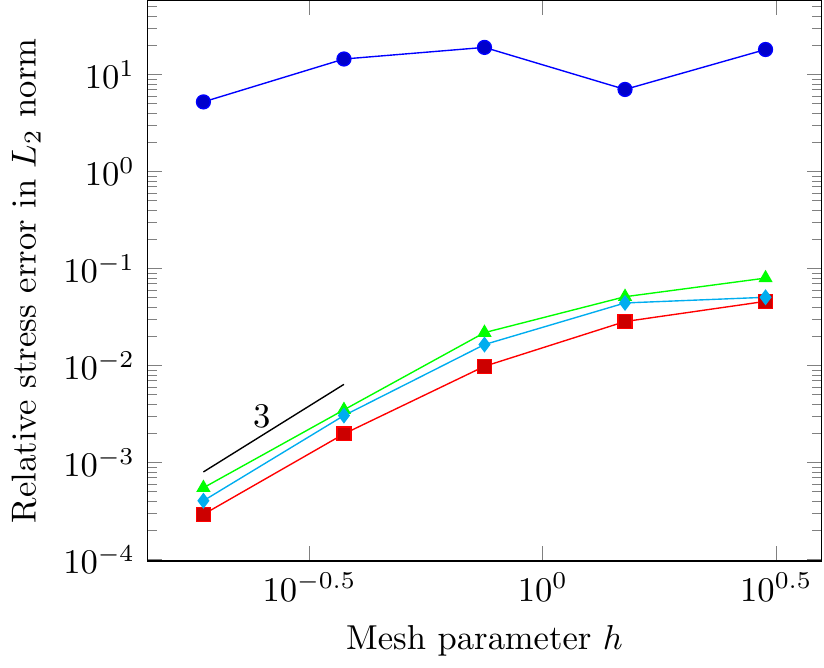} & \includegraphics[width=.31\linewidth]{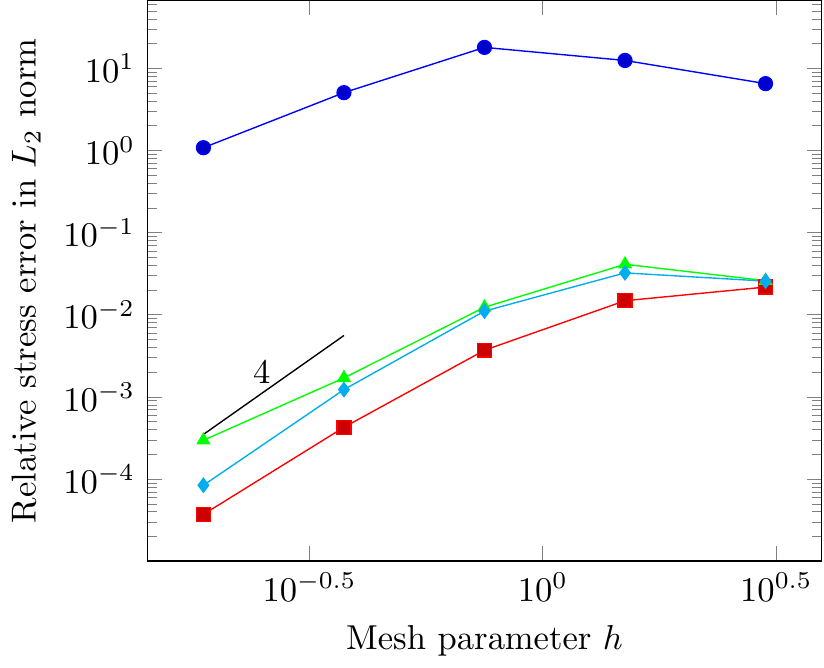}\\
\includegraphics[width=.31\linewidth]{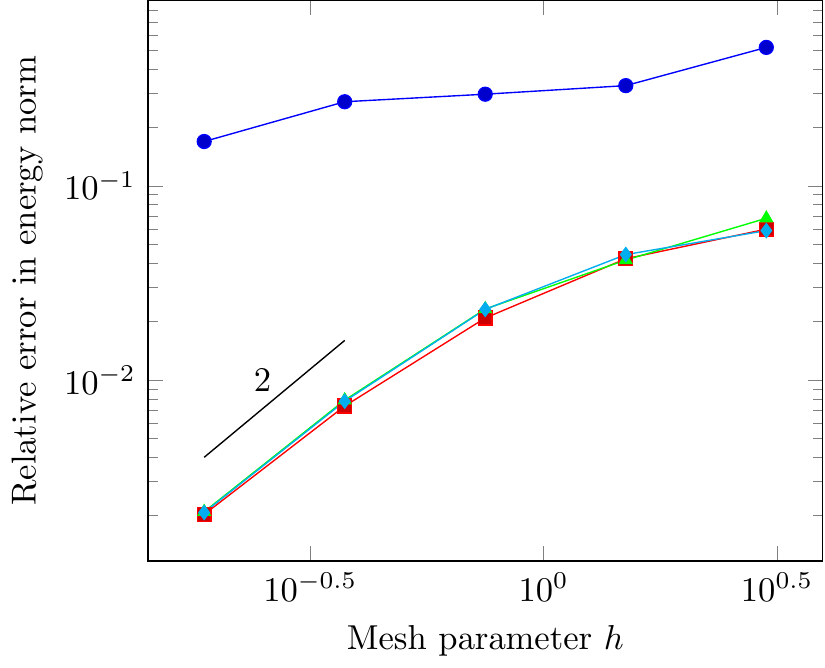} & \includegraphics[width=.31\linewidth]{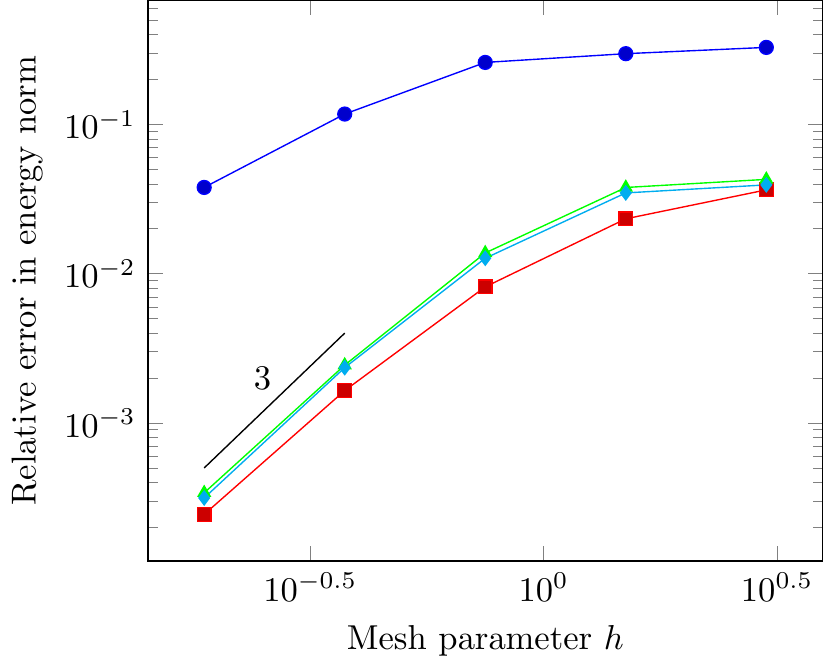} & \includegraphics[width=.31\linewidth]{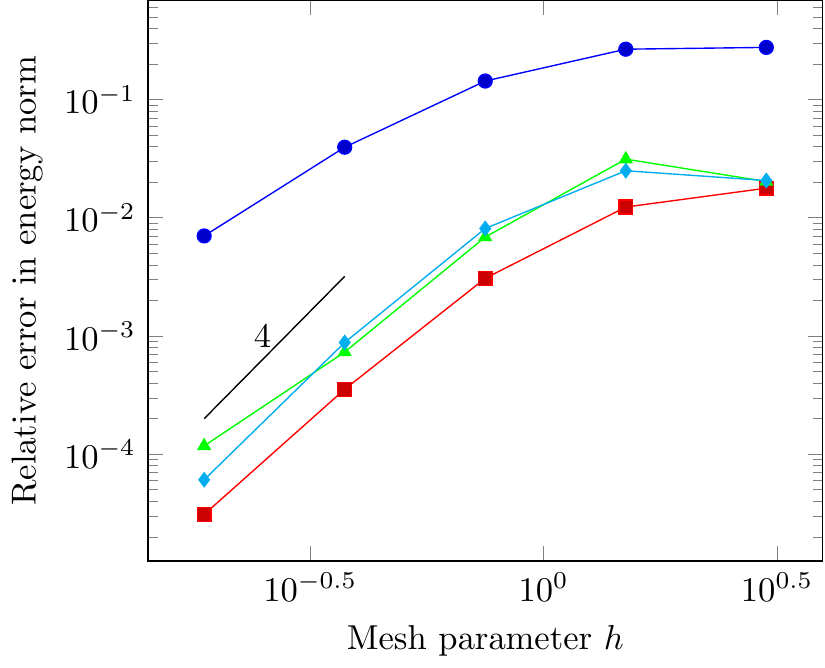}
\end{tabular}
    \caption{Convergence study of the plate with a circular hole. The relative $L^2$ error of displacement, stress and the relative error in energy norm with respect to mesh refinement.}
	\label{fig:platewithhole_convergence}
\end{figure}

Convergence plots for the relative error of the displacement and energy in the $L^2$ norm are shown in Figure \protect\ref{fig:platewithhole_convergence}. As can be seen, the standard $Q_p$ approximations suffer from severe volumetric locking for all orders while, on the other hand, the $\bar{B}$ methods remedy locking for all cases. {For the symmetric \Bezier $\bar{B}$ method, optimal rates are achieved in all three measures for biquadratic elements and the optimal energy convergence has been achieved for bicubic elements, but convergence has degraded in all three measures for the biquartic elements. This reduction in convergence rates results from the fact that the derivation of the symmetric \Bezier $\bar{B}$ method is purely based on the engineering analogy between the $L^2$ projection and \Bezier projection operations. The non-symmetric \Bezier $\bar{B}$ method, on the other hand, achieves optimal convergence in the displacement, stress, and energy norms for all elements with slightly higher errors than those of the global $\bar{B}$ method.}

\begin{figure}[htb!]
    \center
    \begin{tabular}{cc}
    \includegraphics[width=.4\linewidth]{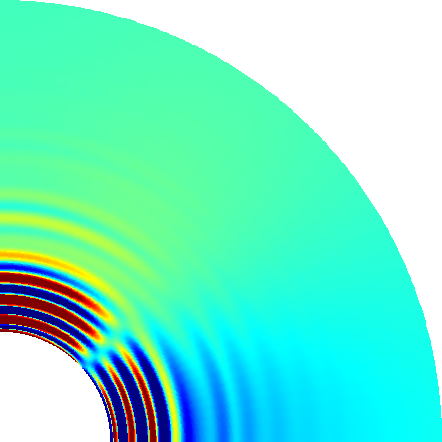} & \includegraphics[width=.4\linewidth]{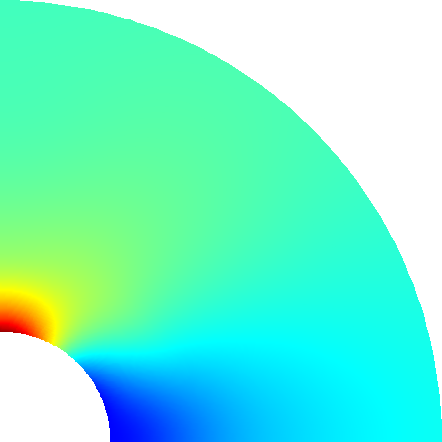}\\
    $Q_4$ & $\mathcal{T}^{L^2} Q_4/Q_3$\\
    \includegraphics[width=.4\linewidth]{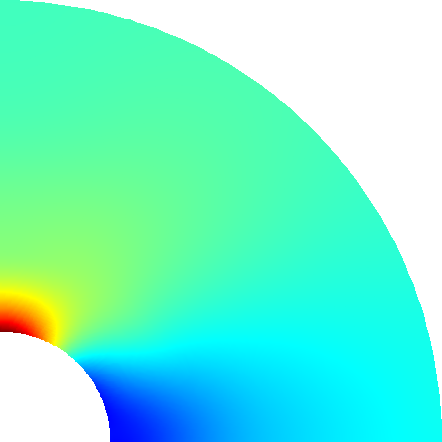} & \includegraphics[width=.4\linewidth]{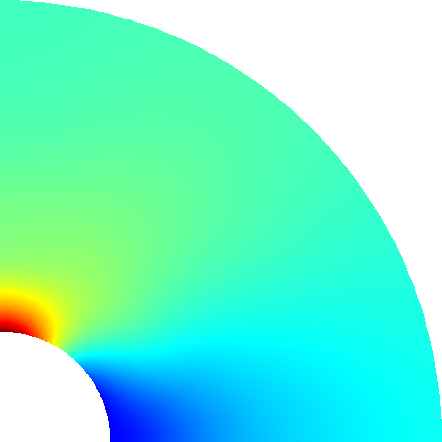}\\
    $S-\mathcal{T}^{p} Q_4/Q_3$ & $NS-\mathcal{T}^{p} Q_4/Q_3$\\
    \includegraphics[width=.4\linewidth]{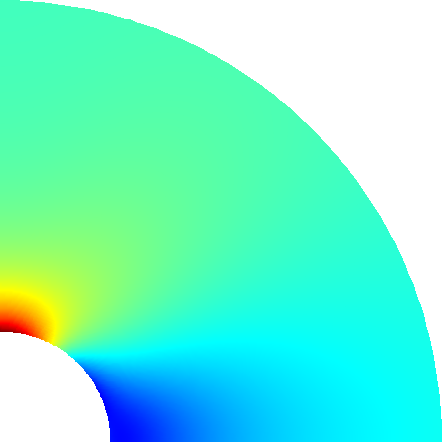} & \includegraphics[width=.14\linewidth]{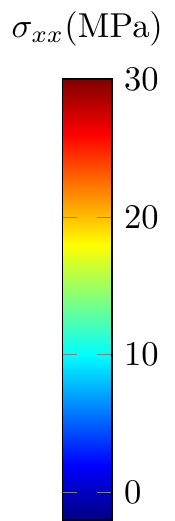}\\
    Reference & Legend
    \end{tabular}
        \caption{ Contour plots of $\sigma_{xx}^h$ for the plate with a circular hole ($p=4$, and the finest mesh is used). For reference the analytical solution is also plotted. }
        \label{fig:platewithhole_contour}
\end{figure}

\begin{figure}[htb!]
    \center
    \begin{tabular}{cc}
    \includegraphics[width=.44\linewidth]{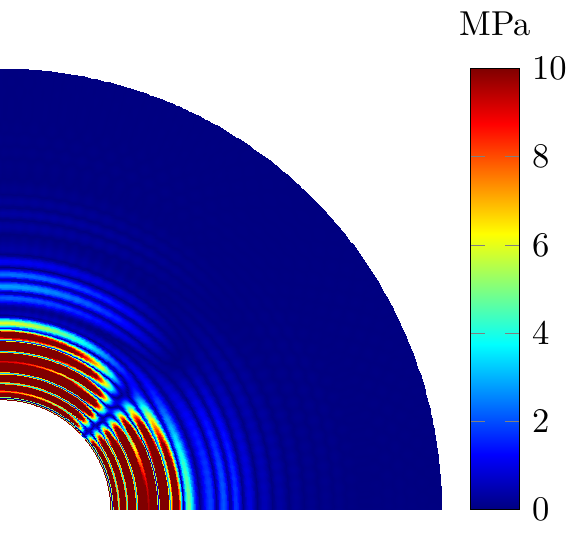} & \includegraphics[width=.5\linewidth]{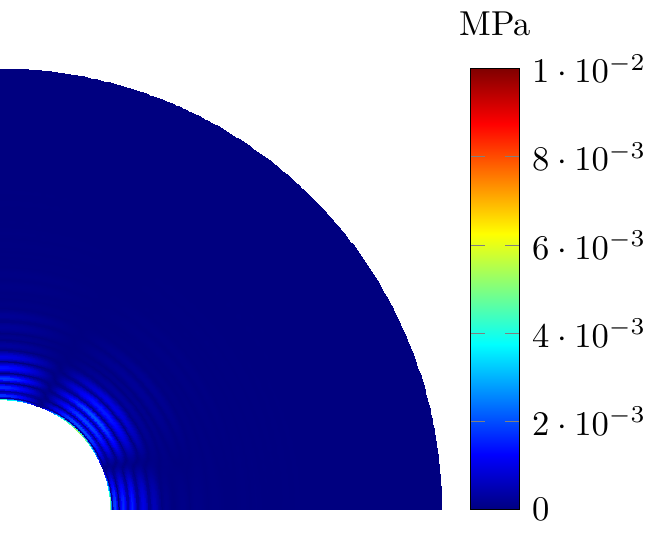}\\
    $Q_4$ & $\mathcal{T}^{L^2} Q_4/Q_3$\\
    \includegraphics[width=.5\linewidth]{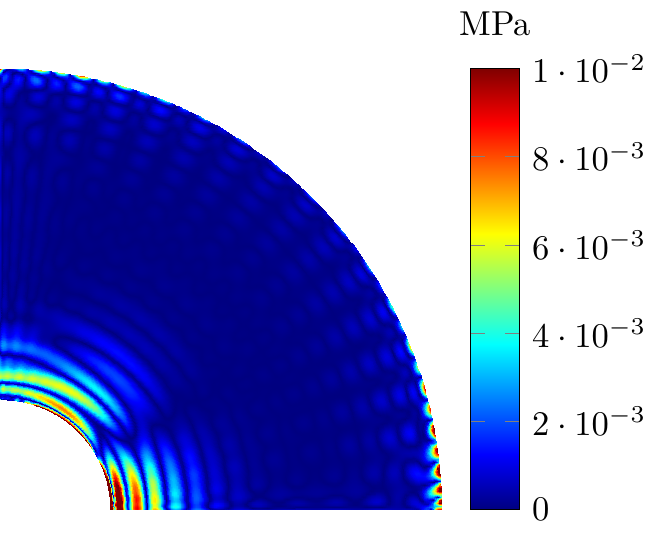} & \includegraphics[width=.5\linewidth]{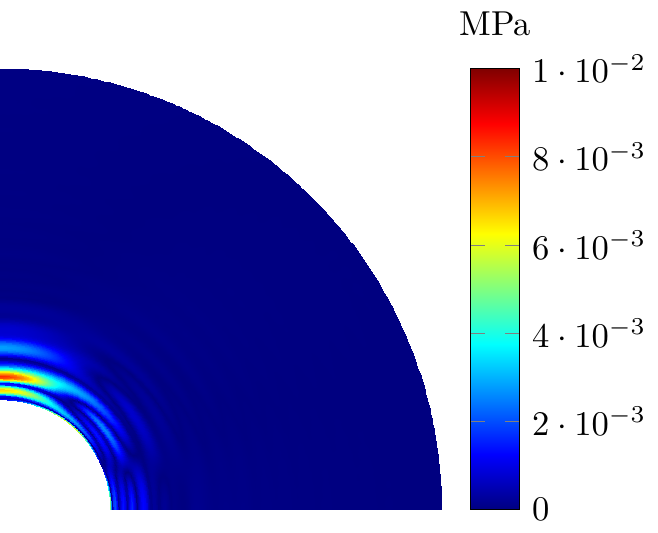}\\
    $S-\mathcal{T}^{p} Q_4/Q_3$ & $NS-\mathcal{T}^{p} Q_4/Q_3$
    \end{tabular}
        \caption{Contour plots of $\vert{\sigma_{xx}-\sigma_{xx}^h}\vert$ for the plate with a circular hole ($p=4$, and the finest mesh is used).}
        \label{fig:platewithhole_error_contour}
\end{figure}

{Contour plots of  $\sigma_{xx}$ from the finest biquartic discretization are shown in Figure \protect\ref{fig:platewithhole_contour}. We can see that results from all $\bar{B}$ methods are consistent with the reference solution, but using the standard finite element approach results in meaningless stresses. Figure \protect\ref{fig:platewithhole_error_contour} shows the absolute error of $\sigma_{xx}$ from the same discretization. We can see that the projection methods produce an error of less than $.1\%$ of the maximum $\sigma_{xx}$, while the error for the standard $Q_4$ element is of the same order as the maximum $\sigma_{xx}$. We can also see that the non-symmetric method provides a slight improvement when compared to the symmetric method.}

\section{Conclusions}
\label{sec:conclusion}

We have presented two B\'{e}zier $\bar{B}$ projection methods, which we have called symmetric and non-symmetric \Bezier $\bar{B}$ projection, as an approach to overcome locking phenomena in structural mechanics applications of isogeometric analysis. Each approach maintains the sparsity of the resulting linear system. The methods utilize B\'ezier extraction and projection, which makes it simple to implement them in an existing finite element framework and makes it applicable to any spline representation which can be written in \Bezier form. In contrast to global $\bar{B}$ methods, which produce dense stiffness matrices, the B\'ezier $\bar{B}$ approach results in a sparse stiffness matrix while still benefiting from higher-order convergence rates. {We have made the connection between the non-symmetric method and a mixed formulation and shown that, although this method does not strictly satisfy the inf-sup condition, it reduces constraints sufficiently to provide optimal convergence rates for the problems studied here.}

We have demonstrated the performance of the approach in the context of shear deformable beams (to alleviate transverse shear locking) and nearly incompressible elasticity problems (to alleviate volumetric locking). The proposed method reduces locking errors and achieves nearly optimal convergence rates { for the symmetric method and optimal rates for the non-symmetric method}. The cases where optimal rates were not achieved {when using the the symmetric formulation are a symptom of the fact that the symmetric formulation is not directly related to a variational principle}.

{The two methods presented here provide a choice between a formulation that results in a symmetric stiffness matrix but requires matrix operations at the global level and potentially less accuracy and a formulation that results in a non-symmetric stiffness matrix that can be assembled in the standard element routine approach and achieves optimal convergence rates. The trade-offs are between higher costs in assembly for the symmetric formulation versus potentially higher costs in solving a non-symmetric system. In either case, however, the cost is less than using the standard $\bar{B}$ formulation.}
\bibliographystyle{elsarticle-num}

\end{document}